\documentclass[12pt]{article}
\usepackage{amsmath,amsthm,amsfonts,amssymb,amscd}
\def\qed{{\unskip\nobreak\hfil\penalty50
\hskip2em\hbox{}\nobreak\hfil$\square$
\parfillskip=0pt \finalhyphendemerits=0\par}\medskip}
\def\proof{\trivlist \item[\hskip \labelsep{\bf Proof\ }]}
\def\endproof{\null\hfill\qed\endtrivlist}

\def\Ad{{\mathrm {Ad}}}

\def\Hom{{\mathrm {Hom}}}

\def\o{{\mathrm {opp}}}

\def\Tr{{\mathrm {Tr}}}

\def\dim{{\mathrm {dim}}}
\def\a{\alpha}
\def\b{\beta}

\def\e{\varepsilon}

\def\l{\lambda}

\def\phi{\varphi}

\def\o{\omega}
\def\Om{\Omega}

\def\s{{\sigma}}

\def\r{{\rho}}

\def\emptyset{\varnothing}
\def\setminus{\smallsetminus}

\def\Diff{{\mathrm {Diff}}}
\def\Mob{{\rm\textsf{M\"ob}}}

\def\Ad{{\mathrm {Ad}}}

\def\Hom{{\mathrm {Hom}}}

\def\o{{\mathrm {opp}}}

\def\Tr{{\mathrm {Tr}}}

\def\dim{{\mathrm {dim}}}
\def\a{\alpha}
\def\b{\beta}

\def\e{\varepsilon}

\def\l{\lambda}

\def\phi{\varphi}

\def\o{\omega}
\def\Om{\Omega}

\def\s{{\sigma}}

\def\r{{\rho}}

\newtheorem{theorem}{Theorem}[section]
\newtheorem{lemma}[theorem]{Lemma}
\newtheorem{conjecture}[theorem]{Conjecture}
\newtheorem{corollary}[theorem]{Corollary}
\newtheorem{definition}[theorem]{Definition}

\newtheorem{proposition}[theorem]{Proposition}
\newtheorem{remark}[theorem]{Remark}

\def\emptyset{\varnothing}
\def\setminus{\smallsetminus}

\def\Diff{{\mathrm {Diff}}}
\def\Mob{{\rm\textsf{M\"ob}}}

\def\res{\!\restriction\!}
\def\A{{\cal A}}
\def\g{{\bold g}}
\def\B{{\cal B}}

\def\D{{\cal D}}

\def\K{{\cal K}}
\def\I{{\cal I}}

\def\L{{\cal L}}

\def\M{{\cal M}}

\def\H{{\cal H}}

\def\1#1{{\bf #1}}
\def\2#1{{\mathcal #1}}
\def\3#1{{\sl #1}}
\def\4#1{{\tt #1}}
\def\5#1{{\sf #1}}
\def\6#1{{\mathfrak #1}}
\def\7#1{{\mathbb #1}}

\renewcommand{\qed}{\ \hfill $\blacksquare$}

\newcommand{\bdefin}{\begin{definition}}
\newcommand{\blemma}{\begin{lemma}}
\newcommand{\bprop}{\begin{proposition}}
\newcommand{\btheor}{\begin{theorem}}
\newcommand{\bcoro}{\begin{corollary}}
\newcommand{\bconj}{\begin{conjecture}}
\newcommand{\edefin}{\end{definition}}
\newcommand{\elemma}{\end{lemma}}
\newcommand{\eprop}{\end{proposition}}
\newcommand{\etheor}{\end{theorem}}
\newcommand{\ecoro}{\end{corollary}}
\newcommand{\econj}{\end{conjecture}}
\newcommand{\brem}{\begin{remark}}
\newcommand{\erem}{\end{remark}}

\newcommand{\ba}{\begin{array}}
\newcommand{\ea}{\end{array}}
\newcommand{\bea}{\begin{eqnarray}}
\newcommand{\eea}{\end{eqnarray}}
\newcommand{\bean}{\begin{eqnarray*}}
\newcommand{\eean}{\end{eqnarray*}}

\title{\huge Solitons in Affine and Permutation Orbifolds\\}
\author{
{\sc Victor G. Kac}\footnote{Supported in part by NSF.}\\
Department of Mathematics\\
MIT\\
Cambridge, MA 02139\\
E-mail: {\tt kac@math.mit.edu}\\
{}\\
{\sc Roberto Longo}\footnote{Supported in part by GNAMPA-INDAM and 
MIUR.}\\
Dipartimento di Matematica\\
Universit\`a di Roma ``Tor Vergata''\\
Via della Ricerca Scientifica, 1, I-00133 Roma, Italy\\
E-mail: {\tt longo@mat.uniroma2.it}\\
{}\\
{\sc Feng Xu}\footnote{Supported in part by NSF.}\\
Department of Mathematics\\
University of California at Riverside\\
Riverside, CA 92521\\
E-mail: {\tt xufeng@math.ucr.edu}}
\begin{document}
\date{}
\maketitle

\begin{abstract}
We consider properties of solitons 
in general orbifolds in the algebraic quantum 
field theory framework and constructions of solitons in affine and 
permutation orbifolds. Under general conditions we show that our 
construction gives all the twisted representations of the fixed point
subnet. This allows us to prove a number of conjectures: in the affine
orbifold case we clarify the issue of ``fixed point resolutions'';
in the permutation orbifold case we determine
all irreducible representations of the orbifold,
and we  also determine the fusion rules in a nontrivial case, 
which imply an integral 
property of chiral data for any completely rational conformal net.
\end{abstract}

\newpage

\section{Introduction}
Let $\A$ be a completely rational conformal net (cf. \S \ref{complete 
rationality} and def. \ref{abr} following \cite{KLM}).  
Let $\Gamma$ be a finite group acting properly on $\A$ (cf. definition
(\ref{p'})). The starting point of this paper is Th. \ref{orb} proved
in \cite{Xu3} which states that the fixed point subnet (the orbifold)
$\A^\Gamma$ is also  completely rational, and by \cite{KLM} 
$\A^\Gamma$ has finitely many irreducible representations which are 
divided into two classes: the ones that are obtained from the 
restrictions of a representation of $\A$ to $\A^\Gamma$ which 
are called untwisted representations, and the ones which are
twisted (cf. definition after Th. \ref{orb}). It follows from 
 Th. \ref{orb} that twisted representation of $\A^\Gamma$
always exists if  $\A^\Gamma\neq \A.$ The motivating question for
this paper is how to construct these   twisted representations of 
$\A^\Gamma.$\par
It turns out that all  representations of 
$\A^\Gamma $ are closely related to the solitons of $\A$ (cf.
\S3.3 and Prop. \ref{sg}). Solitons are representations of $\A_0$, 
the restriction of $\A$ to the real line identified with a circle with 
one point removed. Every representation of $\A$ restricts to a soliton
of $\A_0$, but not every soliton of $\A_0$ can be extended to a 
representation of $\A$. In \S4 we develop general theory of solitons
in the case of orbifolds with two main results: Th. \ref{solitonind}
gives a formula for the index of of solitons obtained from restrictions,
and Th. \ref{projrep} clarifies the general structure of the restriction
of a soliton. These results are natural extensions of similar results
in \cite{LX} and \cite{Xu5} in special cases. \par
The construction of solitons depends on the net $\A$
and the action of $\Gamma$. In the case of affine orbifold, our 
construction (cf. def. \ref{affinesl}) 
is partially inspired by the ``twisted representations'' of 
\cite{KT}, and in fact can be viewed as an ``exponentiated  version''
of   the ``twisted representations'' of 
\cite{KT} (cf. \S5.2.1). Combined with the general properties 
of solitons described above, 
this construction allows us to clarify the issue of ``fixed point''
problem in \cite{KT} in Th. \ref{fix}, and we also show that our
construction gives all the irreducible representations of the 
fixed point subnet under general conditions in Th. \ref{allaffine}
and Cor. \ref{allaa}, thus answering our motivating question in this
case. \par
In the case of permutation orbifolds (cf. \S6), our construction of
solitons in (\ref{sp}) 
is a simple generalizations of the construction of solitons
in \cite{LX} for the case of cyclic orbifolds. Note that the     
construction of solitons
in \cite{LX} also leads to structure results such as a dichotomy for
any split local conformal net. In Th. \ref{allcyclic} (resp. Th. \ref{allpermu}) 
we show that our construction 
gives all the irreducible representations of the cyclic orbifold (resp.
the permutation orbifold), and in Th. \ref{decomposecy} (resp. 
Th. \ref{decomposepermu}) we list 
all the irreducible representations of the cyclic orbifold (resp.
the permutation orbifold). These results generalize the results of
\cite{LX} and prove a claim in \cite{Ba} which is based on heuristic 
arguments. Using theses results in  \S9 
we determine the fusion rules for the first nontrivial 
case when $n=2$ in Th. \ref{n=2f} which implies an integral property of 
the chiral data for any completely rational net (cf. Cor. \ref{integer}),
proving a conjecture in the paper \cite{BHS}, that contained the first 
computations leading to correct fusion rules. 
\par
The rest of this paper is organized as follows: \S2 and \S3 are 
preliminaries on the algebraic quantum field theory framework where 
orbifold construction is considered. 
In these sections we have collected some basic notions that
appear in this paper for the convenience of the reader who may
not have an operator algebra background.
The results in \S2,\S3 are
known except Prop. \ref{normal} on extensions of solitons 
which plays an important role in \S7. In \S4 we apply the results in
\S2 and \S3 to obtain general properties of solitons under inductions
and restrictions, and in particular we prove Th. \ref{solitonind} and
Th. \ref{projrep}. In \S5, after recalling basic definitions 
and properties in affine orbifold from \cite{KT}, we give
the constructions of solitons in \S5.2 and in \S5.2.1 compare it with the 
twisted representations in \cite{KT}. Th. \ref{allaffine} and its
Cor. \ref{allaa} are proved in \S5.3. In \S5.5 we clarify the issue of
fixed point resolutions in Th. \ref{fix}, Cor. \ref{cyclic2}. In
\S5.6 we illustrate the results of \S5.5 in an example 
considered in \cite{KT}.   In \S6 we first recall the construction of
solitons from \cite{LX} in the cyclic permutation case, and in \S6.3
give the general construction of solitons for permutation orbifolds. 
We prove in \S7 the important property of these solitons (Th. \ref{permu}) 
which so far has no direct proof. In \S8 we apply the results of
previous sections to prove four theorems which are briefly described
above. In \S9 after proving some simple properties of $S$ matrix
(cf. Lemma \ref{Smatrix1}), we determine the fusions of solitons
in cyclic orbifold in a special case in Prop. \ref{degreen}. In \S9.3
we determine the fusion rules for the case $n=2$ in Th. \ref{n=2f}
which implies an integral property in Cor. \ref{integer}.  

\section{Elements of Operator Algebras and Conformal QFT}

For the convenience of the reader we collect here some basic notions 
that appear in this paper. This is only a guideline and the reader 
should look at the references for a more complete treatement.

\subsection{von Neumann algebras}
Let $\H$ be a Hilbert space that  we always assume  to be separable to
simplify the exposition. With $B(\H)$ the algebra of all bounded linear
operator on $\H$ a {\em von Neumann algebra} $M$ is a $^*$-subalgebra
of $B(\H)$ containing the identity operator such that $M = M^{-}$
(weak closure or, equivalently, strong closure).

Equivalently $M = M^{\prime\prime}$
(von Neumann density theorem), where the prime denotes the 
commutant: $M'\equiv\{a\in B(\H): xa=ax\ \forall x\in M\}$.

A linear map $\eta$ from a von Neumann algebra $M$ to a von Neumann 
algebra $N$ is  {\em positive} if $\eta(M_+)\subset N_+$, where 
$M_+\equiv\{x\in M:x>0\}$ denotes the cone of positive elements of $M$. 
$\eta$ is {\em normal} if commutes with the sup operation,
namely $\sup\eta(x_i)=\eta(\sup x_i)$ for any bounded increasing 
net of elements in $M_+$; $\eta$ is normal iff it is 
weakly (equivalently strongly) continuous on the unit ball of $M$. $\eta$ 
is {\em faithful} if $\eta(x)=0$, $x\in M_+$, implies $x=0$.
By a {\em homomorphism} of a von Neumann algebra we shall always mean an 
identity preserving homomorphism commuting with the $^*$-operation, 
and analogously for {\em isomorphisms} and {\em endomorphisms}. Isomorphisms 
between von Neumann algebras are automatically normal. 
By a {\em representation} of 
$M$ on a Hilbert space $\K$ we mean a homomorphism of $M$ into $B(\K)$.

A {\em state} $\o$ on von Neumann algebra $M$ is a positive linear 
functional on $M$ with the normalization $\o(1) = 1$. The relevant 
states for a von Neumann algebras are the normal states. By the GNS 
construction, every normal state of $M$ is given by 
$\o(x)=(\pi(x)\Omega,\Omega)$, where $\pi$ is a normal representation 
of $M$ on a Hilbert space $\K$ and $\Omega\in \K$ is cyclic (i.e. 
$\pi(M)\Omega$ is dense in $\K$, see below). Given $\o$, the triple 
$(\K,\pi,\Omega)$ is unique up to unitary equivalence. 

A {\em factor} is a von Neumann algebra with trivial center, namely 
$M\cap M' =\mathbb C$. We note that a factor is a simple algebra,
i.e., the only weakly closed ideal of the factor is either trival
or equal to the factor itself.
If $M$ is a factor (and $\K$ is separable), 
a representation of $M$ on $\K$ is automatically normal. 

A factor $M$ is {\em finite} if there exists a tracial state $\o$ on 
$M$, namely $\o(xy)=\o(yx)$, $x,y\in M$
(automatically normal and unique). Otherwise $M$ is called an 
{\em infinite factor}. For a factor $M$, the following are equivalent:
\begin{itemize}
\item $M$ is infinite;
\item $M$ is isomorphic to $M\otimes B(\K)$, with $\K$ a 
separable infinite dimensional Hilbert space;
\item $M$ contains a non-unitary isometry (an isometry $v$ is an 
operator with the property $v^*v=1$);
\item $M$ contains a non degenerate Hilbert $H$ space of isometries 
with arbitrary dimension (but separable).
\end{itemize}
Here a {\em Hilbert space of isometries $H$ in $M$} we mean a norm 
closed linear subspace $H\subset M$ such that $x^* y\in\mathbb C$ for 
all $x,y\in M$. Thus $x,y\mapsto y^* x$ is scalar product on $H$.
 Then, if $L$ is a set with $\{v_{\ell}, \ell\in L\}$ 
an orthonormal basis for $H$, we have $v^*_{\ell} v_{\ell'} 
=\delta_{\ell \ell'}$, namely the $v_i$'s are isometries $H$ with 
pairwise orthogonal range projections.  $H$ is non-degenerate if the 
left support of $H$ is 1,  that is the final projections form a 
partition of the identity: $\sum_{\ell\in L} v_{\ell} v^*_{\ell} = 1$.

A factor $M$ is {\em of type III} (or purely infinite) if every 
non-zero projection $e\in M$ is equivalent to the identity, namely 
there exists an isometry $v\in  M$ wit $vv^* = 1$. As we shall see, 
factors appearing 
in CFT as local algebras are of type III and the reader may focus on 
this case for the need of this paper.

A {\em semifinite} factor is a factor $M$ isomorphic to $M_0 \otimes 
B(\K)$ with $M_0$ a finite factor and $\K$ a Hilbert space. 
Semifinite factor are characterized by the existence of a normal, 
possible unbounded, trace (that we do not define here). A factor is 
either semifinite or of type III.

A factor $M$ of type III has only one representation (on a separable 
Hilbert space) up to unitary equivalence. Namely, if $\pi: M\to B(\K)$ 
is a representation, there exists a unitary $U:\H\to\K$ such that
$\pi(x)=U x U^*$, $x\in M$.

{\sc Refs:} \cite{T}.

\subsection{Tomita-Takesaki modular theory} Let $M$ be a von Neumann
algebra and $\omega $ a normal faithful state on $M$.
By the GNS construction, we may assume that $\omega =
(\,\cdot\,\Omega,\Omega)$ with $\Omega$ a cyclic and separating vector
($\M$ acts standardly).

Here a vector $\Omega$ is {\em cyclic} if $\overline{M\Omega}=\H$ and 
{\em separating} if $x\in M$, $x\Omega = 0$ implies $x=0$. A vector is 
cyclic for $M$ iff it is separating for $M'$.

The anti-linear operator $x\Omega\mapsto x^*\Omega$, $x\in M$, is 
closable and its closure is denoted by $S$. The polar decomposition 
$S=J\Delta^{1/2}$ gives a antiunitary involution $J$, the {\em modular 
involution}, and a positive 
non-singular linear operator $\Delta\equiv S^*S$, the {\em modular 
operator}.

We have 
\begin{gather}
\Delta^{it}M\Delta^{-it}= M\\
JM J= M^\prime ,
\end{gather}
in other words the
modular theory associates with $\Omega $ a canonical ``evolution", i.e. a
one-parameter group of {\it modular automorphisms} of $M$
$\sigma^\omega_t\equiv \Ad\Delta^{it}$ and an anti-isomorphism 
$\Ad J$ of $\M$ with $\M'$.

Let $N\subset M$ be an inclusion of
von Neumann algebras. We always assume that $N$ and $M$ have the same 
identity. A {\em conditional expectation} $\e:M\to N$ is a 
positive, unital map from $M$ onto $N$ such that 
$\e(n_1 xn_2)=n_1\e(x)n_2$, $x\in M$, $n_1, n_2\in N$.

If $\o$ is a faithful normal state of $M$, by Takesaki theorem there 
exists a normal conditional expectation $\e:M\to N$ preserving $\o$ 
(i.e. $\o\cdot\e =\o$) if and only if $N$ is globally invariant under 
the modular group $\s^{\o}$ of $M$.

If $\r$ is an endomorphism of $M$ and $\e: M\to \r(M)$ is a 
conditional expectation, the map $\varphi\equiv \r^{-1}\cdot\e$ 
satisfies $\varphi\cdot\r = {\rm id}$ and is called a {\em left 
inverse} of $\r$. 

{\sc Refs:} \cite{T}.

\subsection{Jones index}

Let $N\subset M$ an inclusion of factors. The index of $N$ in $M$ can 
be defined by different point of views: analytic, probabilistic or tensor 
categorical.

{\em Analytic definition.} The index was originally considered by Jones in 
the setting of finite factors.
Assume $M$ to be
finite and let $\o$ be the faithful tracial state $\omega$ on $M$. 
As above we may assume that $\o$ is the vector state given by the 
vector $\Omega$. With $e$ the projection onto 
$\overline{N\Omega}$, the von Neumann algebra generated by $M$ and $e$
\[
M_1 = \{M,e\}'' = J_{M} N'J_{M}
\]
is a semifinite factor. $N\subset M$ has finite index iff $M_1$ is finite 
and the {\em index} is then defined by $\lambda = \omega(e)^{-1}$ with $\omega$ also
denoting the tracial state of $M_1$. Jones theorem shows the possible values
for the index:
\[
\lambda \in \left\{4\text{cos}^2\frac{\pi}n, n\geq 3\right\}\cup\left
[4,\infty\right].
\]
If $N\subset M$ is an inclusion of finite factor, there exists a 
unique trace-preserving conditional expectation $\e:M\to N$ ($\s^{\o}$ 
is trivial in this case).

A definition for the index $[M:N]_\e$ of an arbitrary
inclusion of factors $N\subset M$ with a faithful normal conditional expectation
$\e : M\to N$ was given by Kosaki using Connes-Haagerup 
dual weights. It depends on the choice of $\varepsilon$. Given $\e$, 
choose a normal faithful state $\o$ of $M$ with $\o\cdot\e = \o$ and 
$\Omega$ a cyclic vector implementing $\o$.
If $[M:N]_\e 
<\infty$, it is possible to define a canonical expectation $\e':M_1\to 
M$ and then $[M:N]_\e = \e'(e)^{-1}$, with $e$ the projection onto 
$\overline{N\Omega}$. Jones restriction on the index values holds for $[M:N]_\e$
as well.

The good properties are shared by  the {\it minimal index} 
\[
[M:N] = \underset{\varepsilon}{\inf}
[M:N]_\varepsilon = [M:N]_{\varepsilon_0}
\]
where $\varepsilon_0$ is the unique {\it minimal  conditional
expectation}. 

The analytic point of view will not play an explicit role in this 
paper.

{\em Probabilistic definition.}
Pimsner and Popa inequality, and its extension to the infinite factor 
case, shows that $\l\equiv [M:N]^{-1}_\e$ is the 
best constant such that
\[
\varepsilon(x)\geq \lambda x, \quad x\in M^+,
\]
where $\varepsilon: M\to N$ a normal conditional
expectation (if $M$ is finite-dimensional $\l$ is not an optimal bound).

This gives a general way to define the index and a powerful 
tool to check whether a given inclusion has finite index.

{\em Tensor categorical definition.} We shall get to this point in a 
moment.

{\sc Refs:} \cite{J,Ko,L1,L5,LR,PP,T} and references therein.

\subsection{Joint modular structure. Sectors}
 Let $N\subset M$ be an
inclusion of infinite factors.  We may
assume that $N'$ and $M'$ are infinite so $M$ and $N$ have a cyclic 
and separating vector. 
 With $J_{N}$ and $J_{M}$  modular conjugations of $N$ and
$M$, the unitary $\Gamma = J_{N} J_{M}$ implements a {\it canonical
endomorphism} of $M$ into $N$ 
\[
\gamma(x) = \Gamma
x\Gamma^*,\qquad x\in M.
\]
$\gamma$ depends on the choice of
$J_{N}$ and $J_{M}$ only up to perturbations by an inner automorphism
of $M$ associated with a unitary in $N$. The restriction $\gamma | N$ 
is called the {\em dual canonical endomorphism} (it is the canonical 
endomorphism associated with $\gamma(M)\subset N$). $\gamma$ is canonical as a
sector of $M$ as we define now.

Given the infinite factor $M$, the {\it sectors of $M$}  are given by
$$\text{Sect}(M) = \text{End}(M)/\text{Inn}(M)$$
namely $\text{Sect}(M)$ is the quotient of the semigroup of the
endomorphisms of $M$ modulo the equivalence relation: $\rho,\rho'\in
\text{End}(M),\, \rho\thicksim\rho'$ iff there is a unitary $u\in M$
such that $\rho'(x)=u\rho(x)u^*$ for all $x\in M$.

$\text{Sect}(M)$ is a $^*$-semiring (there is an addition, a product and
an involution) equivalent to the Connes correspondences (bimodules) on
$M$ up to unitary equivalence. If $\r$ is
an element of $\text{End}(M)$ we shall denote by $[\r]$
its class in $\text{Sect}(M)$. The operations are:

{\it Addition} (direct sum): Let $\rho_1,\rho_2,\dots\r_n\in\text{End}(M)$.
Choose a non-degenerate $n$-dimensional Hilbert $H$ space of isometries in $M$
and a basis $v_1,\dots v_n$ for $H$. Then
\[
\r(x)\equiv\sum_{i=1}^n v_i\r_i(x)v^*_i,\quad x\in M,
\]
is an endomorphism of $M$. The definition of the direct sum 
endomorphism $\r$ does not depend on the choice of $H$ or on the 
basis, up to inner automorphism of $M$,
namely $\r$ is a well-defined sector of $M$.

{\it Composition} (monoidal product). The usual composition of maps
$$\rho_1\cdot\rho_2(x) = \rho_1(\rho_2(x)), \qquad x\in M,$$ defined on
$\text{End}(M)$ passes to the quotient $\text{Sect}(M)$.

{\it Conjugation.} With $\rho\in\text{End}(M)$, choose a canonical
endomorphism $\gamma_\rho : M\to\rho(M)$. Then 
\[
\bar\rho = \rho^{-1}\cdot\gamma_\rho
\]
 well-defines a
conjugation in $\text{Sect}(M)$. By definition we thus have
\begin{equation}\label{canr}
\gamma_\rho = \rho\cdot\bar\rho
\end{equation}

{\sc Refs:} \cite{ILP,Ko,L1'} and references therein.

\subsection{The tensor category $\text{End}(M)$}
With $M$ an infinite factor, then  $\text{End}(M)$ is a
strict {\it tensor $C^*$-category}, as is already implicit in the 
previous section. 

More precisely define a category $\text{End}(M)$ 
whose objects are the elements of  $\text{End}(M)$ and the arrows 
$\text{Hom}(\r,\r')$ between the objects $\r,\r'$ are 
\[
\text{Hom}(\r,\r')\equiv\{a\in M: a\r(x)=\r'(x)a \ \forall x\in M\}.
\]
The composition of intertwiners (arrows) is the operator product.
Clearly $\text{Hom}(\r,\r')$ is a Banach space and there is a 
$^*$-operation $a\in \text{Hom}(\r,\r')\mapsto a^*\in 
\text{Hom}(\r',\r)$ with the usual properties and the $C^*$-norm 
equality $||a^* a|| = ||a||^2$. Thus $\text{End}(M)$ is a  
{\em $C^*$-category}.

Moreover there is a tensor (or monoidal) product in $\text{End}(M)$. 
The tensor product $\r\otimes\r'$ is simply the composition $\r\r'$. 
For simplicity the symbol $\otimes$ is thus omitted in this case: 
$\r\otimes\r'=\r\r'$. If $\s,\s'\in \text{End}(M)$, and $t\in 
\text{Hom}(\r,\r')$, $s\in\text{Hom}(\s,\s')$, the tensor product 
arrow $t\otimes s$ is the element of $\text{Hom}(\r\otimes\s,\r'\otimes\s')$ 
given by
\[
t\otimes s\equiv t\r(s) = \r'(s)t\ .
\]
As usual, there is a natural compatibility between tensor product and 
composition, thus $\text{End}(M)$ is a {\em $C^*$-tensor category}. Moreover there 
is an identy object $\iota$ for the tensor procuct (the identity automorphism).

So far we have not made much use that $M$ is an infinite factor. This 
enters crucially for the conjugation in $\text{End}(M)$.

If $\rho$ is irreducible (i.e. $\rho(M)'\cap M =
\mathbb C$) and has finite index, then $\bar\rho$ is the unique  sector such
that $\rho\bar\rho$ contains the identity sector.
More generally the objects $\rho,\bar\rho\in\text{End}(M)$ are conjugate
according to the  analytic definition and have finite index if and only
if there exist isometries $v\in \Hom(\iota,\rho\bar\rho)$ and $\bar v\in
\Hom(\iota,\bar\rho\rho)$ such that
\[
\bar v^*\otimes 1_{\bar\r} \cdot 1_{\bar\r}\otimes v \equiv  \bar v^*\bar\rho(v)=\frac1d, 
\qquad v^*\otimes 1_{\r} \cdot 1_{\r} \otimes \bar v
\equiv v^*\rho(\bar v)=\frac1d,
\]
for some $d>0$.
 
The minimal possible value of $d$ in the above formulas is the {\it
dimension} $d(\rho)$ of $\rho$; it is related to the minimal index by 
\[
[M:\r(M)]=d(\rho)^2
\]
(tensor categorical definion of the index) 
and satisfies the dimension properties
\smallskip

$d(\rho_1\oplus\rho_2)=d(\rho_1)+d(\rho_2)$

$d(\rho_1\rho_2)=d(\rho_1)d(\rho_2)$

$d(\bar\rho)=d(\rho).$
\smallskip 

It follows that the the subcategory of $\text{End}(M)$ having 
finite-index objects is a {\em $C^*$-tensor category with conjugates} 
and direct sums.

Formula (\ref{canr}) shows that given $\gamma\in\text{End}(M)$ the problem
of deciding whether it is a canonical endomorphism with respect to some
subfactor  is essentially the problem of finding a ``square root"
$\rho$.
$\gamma$ is canonical and has finite index iff 
there exist isometries $t\in \Hom(\iota,\gamma)$,  
$s\in \Hom(\gamma,\gamma^2)$
satisfying the algebraic relations
\begin{gather}
s^*s^*=s^*\gamma(s^*)\\
s^*\gamma(t)\in\mathbb C\backslash\{0\}\ ,\quad
s^*t\in\mathbb C\backslash\{0\}. 
\end{gather}
It is immediate to generalize the notion of $\text{Sect}(M)$ to 
$\text{Sect}(M,N)$, for a pair of factors $M,N$. They are the 
homomorphisms of $M$ into $N$ up to unitary equivalence given by a 
unitary in $N$. If $N\subset M$ is an inclusion of infinite factors, 
the canonical endomorphism $\gamma: M\to N$ is a well defined element of 
$\text{Sect}(M,N)$; if $[M:N]<\infty$, the above formula show that 
$\gamma$ is the conjugate sector of the inclusion homomorphism $\iota_N 
:N\to M$:
\[
\gamma= \bar\iota_N \iota_N,\qquad \gamma\restriction N = \iota_N \bar 
\iota_N\ .
\]

We use $\langle  \lambda , \mu \rangle$ to denote
the dimension of   $\text{\rm Hom}(\lambda , \mu )$; it can be 
$\infty$, but it is finite if $\l,\mu$ have finite index.  
$\langle  \lambda , \mu \rangle$
depends only on $[\lambda ]$ and $[\mu ]$. Moreover we have
if $\nu$ have finite dimension, then 
$\langle \nu \lambda , \mu \rangle =
\langle \lambda , \bar \nu \mu \rangle $,
$\langle \lambda\nu , \mu \rangle
= \langle \lambda , \mu \bar \nu \rangle $ which follows from Frobenius
duality.  
$\mu $ is a subsector of $\lambda $ if there is an isometry $v\in M$ 
such that $\mu(x)= v^* \lambda(x)v, \forall x\in M.$
We will also use the following
notation: if $\mu $ is a subsector of $\lambda $, we will write as
$\mu \prec \lambda $  or $\lambda \succ \mu $.  A sector
is said to be irreducible if it has only one subsector. 

{\sc Refs:} \cite{DHR,L2,LRo} and references therein.

\section{Conformal nets on $S^1$}
\label{nets}

By an interval of the circle we mean an open connected
non-empty subset $I$ of $S^1$ such that the interior of its 
complement $I'$ is not empty. 
We denote by $\I$ the family of all intervals of $S^1$.

A {\it net} $\A$ of von Neumann algebras on $S^1$ is a map 
\[
I\in\I\to\A(I)\subset B(\H)
\]
from $\I$ to von Neumann algebras on a fixed Hilbert space $\H$
that satisfies:
\begin{itemize}
\item[{\bf A.}] {\it Isotony}. If $I_{1}\subset I_{2}$ belong to 
$\I$, then
\begin{equation*}
 \A(I_{1})\subset\A(I_{2}).
\end{equation*}
\end{itemize}
If $E\subset S^1$ is any region, we shall put 
$\A(E)\equiv\bigvee_{E\supset I\in\I}\A(I)$ with $\A(E)=\mathbb C$ 
if $E$ has empty interior (the symbol $\vee$ denotes the von Neumann 
algebra generated). 

The net $\A$ is called {\it local} if it satisfies:
\begin{itemize}
\item[{\bf B.}] {\it Locality}. If $I_{1},I_{2}\in\I$ and $I_1\cap 
I_2=\emptyset$ then 
\begin{equation*}
 [\A(I_{1}),\A(I_{2})]=\{0\},
 \end{equation*}
where brackets denote the commutator.
\end{itemize}
The net $\A$ is called {\it M\"{o}bius covariant} if in addition 
satisfies
the following properties {\bf C,D,E,F}:
\begin{itemize}
\item[{\bf C.}] {\it M\"{o}bius covariance}. 
There exists a strongly 
continuous unitary representation $U$ of the M\"{o}bius group 
$\Mob$ (isomorphic to $PSU(1,1)$) on $\H$ such that
\begin{equation*}
 U(g)\A(I) U(g)^*\ =\ \A(gI),\quad g\in \Mob,\ I\in\I.
\end{equation*}
Note that this implies $\A(\bar I)=\A(I)$, $I\in\I$ (consider a 
sequence 
of elements $g_n\in\Mob$ converging to the identity such that 
$g_n\bar I \nearrow I$).
\item[{\bf D.}] {\it Positivity of the energy}. The generator of the 
one-parameter rotation subgroup of $U$ (conformal Hamiltonian) is 
positive. 
\item[{\bf E.}] {\it Existence of the vacuum}. There exists a unit 
$U$-invariant vector $\Omega\in\H$ (vacuum vector), and $\Omega$ 
is cyclic for the von Neumann algebra $\bigvee_{I\in\I}\A(I)$.
\end{itemize}
By the Reeh-Schlieder theorem $\Omega$ is cyclic and separating for 
every fixed $\A(I)$. The modular objects associated with 
$(\A(I),\Omega)$ have a geometric meaning
\[
\Delta^{it}_I = U(\Lambda_I(2\pi t)),\qquad J_I = U(r_I)\ .
\]
Here $\Lambda_I$ is a canonical one-parameter subgroup of $\Mob$ and $U(r_I)$ is a 
antiunitary acting geometrically on $\A$ as a reflection $r_I$ on $S^1$. 

This implies {\em Haag duality}: 
\[
\A(I)'=\A(I'),\quad I\in\I\ ,
\]
where $I'$ is the interior of $S^1\setminus I$.

\begin{itemize}
\item[{\bf F.}] {\it Irreducibility}. $\bigvee_{I\in\I}\A(I)=B(\H)$. 
Indeed $\A$ is irreducible iff
$\Om$ is the unique $U$-invariant vector (up to scalar multiples). 
Also  $\A$ is irreducible
iff the local von Neumann 
algebras $\A(I)$ are factors. In this case they are III$_1$-factors in 
Connes classification of type III factors
(unless $\A(I)=\mathbb C$ for all $I$).
\end{itemize}
By a {\it conformal net} (or diffeomorphism covariant net)  
$\A$ we shall mean a M\"{o}bius covariant net such that the following 
holds:
\begin{itemize}
\item[{\bf G.}] {\it Conformal covariance}. There exists a projective 
unitary representation $U$ of $\Diff(S^1)$ on $\H$ extending the unitary 
representation of $\Mob$ such that for all $I\in\I$ we have
\begin{gather*}
 U(g)\A(I) U(g)^*\ =\ \A(gI),\quad  g\in\Diff(S^1), \\
 U(g)xU(g)^*\ =\ x,\quad x\in\A(I),\ g\in\Diff(I'),
\end{gather*}
\end{itemize}
where $\Diff(S^1)$ denotes the group of smooth, positively oriented 
diffeomorphism of $S^1$ and $\Diff(I)$ the subgroup of 
diffeomorphisms $g$ such that $g(z)=z$ for all $z\in I'$.
\par
Let $G$ be a simply connected  compact Lie group. By Th. 3.2
of \cite{FG}, 
the vacuum positive energy representation of the loop group
$LG$ (cf. \cite{PS}) at level $k$ 
gives rise to an irreducible conformal net 
denoted by {\it ${\A}_{G_k}$}. By Th. 3.3 of \cite{FG}, every 
irreducible positive energy representation of the loop group
$LG$ at level $k$ gives rise to  an irreducible covariant representation 
of ${\A}_{G_k}$. 
 
\subsection{Doplicher-Haag-Roberts superselection sectors in CQFT}
The DHR theory was originally made on the 4-dimensional Minkowski 
spacetime, but can be generalized to our setting. There are however 
several important structure differences in the low dimensional case.

A (DHR) representation $\pi$ of $\A$ on a Hilbert space $\H$ is a map 
$I\in\I\mapsto  \pi_I$ that associates to each $I$ a normal 
representation of $\A(I)$ on $B(\H)$ such that
\[
\pi_{\tilde I}\res\A(I)=\pi_I,\quad I\subset\tilde I, \quad 
I,\tilde I\subset\I\ .
\]
$\pi$ is said to be M\"obius (resp. diffeomorphism) covariant if 
there is a projective unitary representation $U_{\pi}$ of $\Mob$ (resp. 
$\Diff^{(\infty)}(S^1)$, the infinite cover of $\Diff(S^1)$ ) on $\H$ such that
\[
\pi_{gI}(U(g)xU(g)^*) =U_{\pi}(g)\pi_{I}(x)U_{\pi}(g)^*
\]
for all $I\in\I$, $x\in\A(I)$ and $g\in \Mob$ (resp. 
$g\in\Diff^{(\infty)}(S^1)$). Note that if $\pi$ is irreducible and 
diffeomorphism covariant then $U$ is indeed a projective unitary 
representation of $\Diff(S^1)$.

By definition the irreducible conformal net is in fact an irreducible 
representation of itself and we will call this representation the {\it 
vacuum representation}.

Given an interval $I$ and a representation $\pi$ 
of $\A$, there is an {\em endomorphism of $\A$ localized in $I$} equivalent 
to $\pi$; namely $\r$ is a representation of $\A$ on the vacuum Hilbert 
space $\H$, unitarily equivalent 
to $\pi$, such that $\r_{I'}=\text{id}\restriction\A(I')$.  

Fix an interval $I_0$ and  endomorphisms $\r,\r'$ of $\A$ localized in $I_0$.
Then the {\em composition} 
(tensor product) $\r\r'$ is defined by
\[
(\r\r')_I=\r_I \r'_I
\]
with $I$ an interval containing $I$. One can indeed define $(\r\r')_I$ 
for an arbitrary interval $I$ of $S^1$ (by using covariance) and get a 
well defined endomorphism of $\A$ localized in $I_0$.
Indeed the endomorphisms of $\A$ localized in a given interval form a tensor 
$C^*$-category. For our needs $\r,\r'$ will be always localized in a 
common interval $I$. 

If $\pi$ and $\pi'$ are representations of $\A$, fix an interval $I_0$ 
and choose endomorphisms $\r,\r'$ localized in $I_0$ with $\r$ equivalent 
to $\pi$ and $\r'$ equivalent to $\pi'$. Then $\pi\cdot\pi'$ is 
defined (up to unitary equivalence)  to be $\r\r'$. The class of a 
DHR representation modulo unitary equivalence is a {\em 
superselection sectors} (or simply a sector).

Indeed the localized endomorphisms of $\A$ for a tensor 
$C^*$-category. For our needs, $\r,\r'$ will be always localized in a 
common interval $I$. 

We now define  the statistics. Given the endomorphism $\r$ of $\A$
localized in $I\in\I$, choose an equivalent endomorphism $\r_0$
localized in an interval $I_0\in\I$ with $\bar I_0\cap\bar I
=\emptyset$ and let $u$ be a local intertwiner in $\Hom(\r,\r_0)$ as
above, namely $u\in \Hom(\r_{\tilde I},\r_{0,\tilde I})$ with $I_0$ following
clockwise $I$ inside $\tilde I$ which is an interval containing
both $I$ and $I_0$.

The {\it statistics operator} $\varepsilon := u^*\r(u) =
u^*\r_{\tilde I}(u) $ belongs to $\Hom(\r^2_{\tilde I},\r^2_{\tilde
I})$.  An elementary computation shows that it gives rise to a
presentation of  the Artin  braid group
$$\epsilon_i\epsilon_{i+1}\epsilon_i =
\epsilon_{i+1}\epsilon_i\epsilon_{i+1},\qquad
\epsilon_i\epsilon_{i'} =\epsilon_{i'}\epsilon_i
\,\quad{\rm if}\,\, |i-i'|\geq 2,$$ where
$\varepsilon_i=\r^{i-1}(\varepsilon)$. The (unitary equivalence
class of the) representation of the Artin braid group thus obtained is the
{\it statistics} of the superselection sector $\r$.

It turns out the endomorphisms localized in a given interval form a {\em 
braided $C^*$-tensor category} with unitary braiding.

The {\em statistics parameter} $\l_\r$ can be defined in general. In 
particular, assume $\r$ to be localized in $I$ and 
$\r_I\in\text{End}((\A(I))$ to be irreducible with a conditional 
expectation $E: \A(I)\to \r_I(\A(I))$, then 
\[
\l_\r:=E(\epsilon)
\]
depends only on the superselection sector of $\r$.

The {\em statistical dimension} $d_{DHR}(\r)$ and the  {\it univalence}
$\omega_\r$ are then defined by
\[
d_{DHR}(\r) = |\lambda_\r|^{-1}\ ,\qquad \omega_\r = \frac{\lambda_\r}{|\lambda_\r|}\ .
\]

{\sc Refs:} \cite{DHR,FRS,L1,L1',LR}.

\subsection{Index-statistics and spin-statistics relations} Let $\rho$ be an 
endomorphism localized in the interval $I$. A natural connection
between the Jones and DHR theories is realized by the {\em index-statistics 
theorem}
\[
\text{Ind}(\rho)= d_{\text{DHR}}(\rho)^2.
\]
Here $\text{Ind}(\rho)$ is $\text{Ind}(\rho_I))$; 
namely $d((\rho_I))= d_{\text{DHR}}(\rho)$.
We will thus omit the suffix DHR in the dimension. 
Since by duality $\rho(\A(I))\subset\A(I)$ coincides
with $\rho(\A(I))\subset \rho(\A(I'))'$ one may rewrite the above index
formula directly in terms of the representation $\rho$. 

The map $\rho\to\rho_I$ is a 
faithful functor of  $C^*$-tensor categories of endomorphism of $|A$ 
localized in $I$ into $\text{End}(M)$  with $M\equiv\A(I)$. 
Passing to  quotient one obtains a natural embedding 
\[
\text{Superselection
sectors}\longrightarrow\text{Sect}(M).
\]
Restricting to finite-dimensional endomorphisms, the above functor is 
{\em full}, namely, given endomorphisms $\r,\r'$ localized in $I$, if 
$a\in \text{Hom}(\r_I,\r'_I)$ then $a$ intertwines the representations $\r$ 
and $\r'$ (this is obviously true also in the infinite-dimensional 
case if there holds the strong additivity property below, but 
otherwise a non-trivial result).

The {\em conformal spin-statistics theorem} shows that
\[
\omega_\r = e^{i 2\pi L_0(\r)}\ ,
\]
where $L_0(\r)$ is the conformal Hamiltonian (the generator of the 
rotation subgroup) in the representation $\r$. The right hand side in 
the above equality is called the {\em univalence} of $\r$.
\par
{\sc Refs:} \cite{GL2,L1}.

\subsection{Genus 0 $S,T$-matrices}
Next we will recall some of the results of \cite{R2}  and
introduce notations. \par        
Let $\{[\lambda], \lambda\in \L \}$ be a finite set  of
all equivalence classes of irreducible,
covariant, finite-index representations of an irreducible local
conformal net $\A$. 
We will denote the conjugate of $[\lambda]$ by
$[{\bar \lambda}]$
and identity sector (corresponding to the vacuum representation) 
by $[1]$ if no confusion arises, and let
$N_{\lambda\mu}^\nu = \langle [\lambda][\mu], [\nu]\rangle $. 
Here $\langle \mu,\nu\rangle$ denotes the dimension of the
space of intertwiners from $\mu$ to $\nu$ (denoted by $\text {\rm
Hom}(\mu,\nu)$).  We will
denote by $\{T_e\}$ a basis of isometries in $\text {\rm
Hom}(\nu,\lambda\mu)$.
The univalence of $\lambda$ and the statistical dimension of
(cf. \S2  of \cite{GL1}) will be denoted by
$\omega_{\lambda}$ and $d{(\lambda)}$ (or $d_{\lambda})$) respectively. \par
Let $\phi_\lambda$ be the unique minimal
left inverse of $\lambda$, define:   
\begin{equation}\label{Ymatrix}
Y_{\lambda\mu}:= d(\lambda)  d(\mu) \phi_\mu (\epsilon (\mu, \lambda)^*
\epsilon (\lambda, \mu)^*), 
\end{equation} 
where $\epsilon (\mu, \lambda)$ is the unitary braiding operator
 (cf. \cite{GL1} ). \par
We list two properties of $Y_{\lambda \mu}$ (cf. (5.13), (5.14) of \cite{R2}) which
will be used in the following:
\begin{lemma}\label{Yprop}
\begin{equation*}
Y_{\lambda\mu} = Y_{\mu\lambda}  = Y_{\lambda\bar \mu}^* = 
Y_{\bar \lambda \bar \mu}.
\end{equation*} 
\begin{equation*}
Y_{\lambda\mu}  = \sum_k N_{\lambda\mu}^\nu \frac{\omega_\lambda\omega_\mu}
{\omega_\nu} d(\nu) .
\end{equation*}
\end{lemma}
We note that one may take the second equation in the above lemma as the
definition of $Y_{\lambda\mu}$.\par
Define
$a := \sum_i d_{\rho_i}^2 \omega_{\rho_i}^{-1}$.
If the matrix $(Y_{\mu\nu})$ is invertible,
by Proposition on P.351 of \cite{R2} $a$ satisfies
$|a|^2 = \sum_\lambda d(\lambda)^2$.
\begin{definition}\label{c0}
Let $a= |a| \exp(-2\pi i \frac{c_0}{8})$ where  $c_0\in {\mathbb R}$ and $c_0$
is well defined ${\rm mod} \ 8\mathbb Z$.
\end{definition} 
Define matrices
\begin{equation}\label{Smatrix}
S:= |a|^{-1} Y, T:=  C {\rm Diag}(\omega_{\lambda})
\end{equation}
where \[C:= \label{dims} \exp(-2\pi i \frac{c_0}{24}).\]  
Then these matrices satisfy (cf. \cite{R2}):
\begin{lemma}\label{Sprop}
\begin{align*}
SS^{\dag} & = TT^{\dag} ={\rm id},  \\
STS &= T^{-1}ST^{-1},  \\
S^2 & =\hat{C},\\
 T\hat{C} & =\hat{C}T=T, 
\end{align*}

where $\hat{C}_{\lambda\mu} = \delta_{\lambda\bar \mu}$ 
is the conjugation matrix. 
\end{lemma}
Moreover
\begin{equation}\label{Verlinde}
N_{\lambda\mu}^\nu = \sum_\delta \frac{S_{\lambda\delta} S_{\mu\delta} 
S_{\nu\delta}^*}{S_{1\delta}}. \
\end{equation}
is known as Verlinde formula. \par
We will refer the $S,T$ matrices
as defined above  as  {\bf genus 0 modular matrices of ${\A}$} since
they are constructed from the fusion rules, monodromies and minimal 
indices which can be thought as  genus 0 {\bf chiral data} associated to
a Conformal Field Theory. \par
We note that in all cases $c_0-c\in 8\mathbb Z$, where $c$ is the central
charge associated with the projective representations of ${\rm Diff}(S^1)$
of the conformal net $\A$ (cf. \cite{K1} or \cite{LX}). We will prove
in Lemma \ref{central} that $c_0-c\in 4\mathbb Z$ under general 
conditions. \par
The commutative algebra  generated by $\lambda$'s with structure
constants $N_{\lambda\mu}^\nu$ is called {\bf fusion algebra} of 
$\A$. If $Y$ is invertible,  
it follows from Lemma \ref{Sprop}, (\ref{Verlinde})  that any nontrivial 
irreducible representation
of the fusion algebra is of the form
$\lambda\rightarrow \frac{S_{\lambda\mu}}{S_{1\mu}}$ for some
$\mu$. \par  
\subsection{The orbifolds}
Let ${\A}$ be an irreducible conformal net on a Hilbert space
${\H}$ and let $\Gamma$ be a finite group. Let $V:\Gamma\rightarrow U({\H})$
be a  unitary representation of $\Gamma$ on ${\H}$. 
If $V:\Gamma\rightarrow U({\H})$ is not faithful, we set $\Gamma':= \Gamma
/{\rm ker} V$.
\begin{definition} \label{p'}
We say that $\Gamma$ acts properly on ${\A}$ if the following conditions
are satisfied:\par
(1) For each fixed interval $I$ and each $g\in \Gamma$, 
$\alpha_g (a):=V(g)aV(g^*) \in {\A}(I), \forall a\in
{\A}(I)$; \par
(2) For each  $g\in \Gamma$, $V(g)\Omega = \Omega, \forall g\in \Gamma$.\par
\label{Definition 2.1}
\end{definition}
We note that if $\Gamma$ acts properly, then $V(g)$, $g\in\Gamma$ commutes with
the unitary representation $U$ of $\Mob$. \par
Define $\B(I):= \{ a\in \A(I) | \alpha_g (a)=a, \forall g\in \Gamma \}$ and 
${\A}^\Gamma(I):={\B}(I)P_0$ on ${\H}_0$ where $\H_0:=\{ x\in \H| 
V(g)x=x, \forall g\in \Gamma \}$ and $P_0$ is the projection
from $\H$ to $\H_0.$  Then $U$ restricts to
an  unitary
representation (still denoted by $U$) of $\Mob$ on ${\H}_0$.
Then:
\bprop 
The map $I\in {\I}\rightarrow {\A}^{\Gamma}(I)$ on $ {\H}_0$ 
together with the  unitary
representation (still denoted by $U$) of \Mob\ on ${\H}_0$
is an
irreducible M\"{o}bius covariant net.
\label{Prop.2.1}
\eprop
The irreducible  M\"{o}bius covariant net in Prop.  \ref{Prop.2.1} 
will be denoted by
${\A}^\Gamma$ and will be called the {\it orbifold of ${\A}$}
with respect to $\Gamma$. We note that by definition 
${\A}^\Gamma= {\A}^{\Gamma'}$. \par
\subsection{Complete rationality }
\label{complete rationality}
We first recall some definitions from \cite{KLM} .
Recall that   ${\I}$ denotes the set of intervals of $S^1$.
Let $I_1, I_2\in {\I}$. We say that $I_1, I_2$ are disjoint if
$\bar I_1\cap \bar I_2=\emptyset$, where $\bar I$
is the closure of $I$ in $S^1$.  
When $I_1, I_2$ are disjoint, $I_1\cup I_2$
is called a 1-disconnected interval in \cite{Xu6}.  
Denote by ${\I}_2$ the set of unions of disjoint 2 elements
in ${\I}$. Let ${\A}$ be an irreducible M\"{o}bius covariant net
as in \S2.1. For $E=I_1\cup I_2\in{\I}_2$, let
$I_3\cup I_4$ be the interior of the complement of $I_1\cup I_2$ in 
$S^1$ where $I_3, I_4$ are disjoint intervals. 
Let 
$$
{\A}(E):= A(I_1)\vee A(I_2), \quad
\hat {\A}(E):= (A(I_3)\vee A(I_4))'.
$$ Note that ${\A}(E) \subset \hat {\A}(E)$.
Recall that a net ${\A}$ is {\it split} if ${\A}(I_1)\vee
{\A}(I_2)$ is naturally isomorphic to the tensor product of
von Neumann algebras ${\A}(I_1)\otimes
{\A}(I_2)$ for any disjoint intervals $I_1, I_2\in {\I}$.
${\A}$ is {\it strongly additive} if ${\A}(I_1)\vee
{\A}(I_2)= {\A}(I)$ where $I_1\cup I_2$ is obtained
by removing an interior point from $I$.
\bdefin\label{abr}
\cite{KLM}
${\A}$ is said to be completely  rational if
${\A}$ is split, strongly additive, and 
the index $[\hat {\A}(E): {\A}(E)]$ is finite for some
$E\in {\I}_2$ . The value of the index
$[\hat {\A}(E): {\A}(E)]$ (it is independent of 
$E$ by Prop. 5 of \cite{KLM}) is denoted by $\mu_{{\A}}$
and is called the $\mu$-index of ${\A}$. If 
the index $[\hat {\A}(E): {\A}(E)]$ is infinity for some
$E\in {\I}_2$, we define the $\mu$-index of ${\A}$ to be
infinity.
\label{Definition 2.2}
\edefin
A formula for the $\mu$-index of a subnet is proved in 
\cite{KLM}. With the result on 
strong additivity for $\A^{\Gamma}$ in \cite{Xu3}, we have the 
complete rationality in following theorem. 

Note that, by our recent results in \cite{LX}, every irreducible, split, 
local conformal net with finite $\mu$-index is automatically strongly additive.
\btheor\label{orb}
Let ${\A}$ be an irreducible M\"{o}bius covariant net and let $\Gamma$
be a finite group acting properly on ${\A}$. Suppose that 
${\A}$ is completely rational. Then:\par
(1): ${\A}^\Gamma$ is completely rational or $\mu$-rational and  
$\mu_{{\A}^\Gamma}= |\Gamma'|^2 \mu_{{\A}}$; \par
(2): There are only a finite number of irreducible covariant 
representations of ${\A}^\Gamma$ (up to unitary equivalence), 
and they give rise to a unitary modular
category as defined in II.5 of \cite{Tu} by the construction as given in
\S1.7 of \cite{Xu8}. 
\label{Th.2.6}
\etheor
Suppose that ${\A}$ and $\Gamma$ satisfy the assumptions of Th. \ref{orb}. 
Then ${\A}^\Gamma$ has only finite 
finite number of irreducible representations $\dot\lambda$ and
$$
\sum_{\dot\lambda}d(\dot\lambda)^2 = \mu_{{\A}^\Gamma}= |\Gamma'|^2 \mu_{{\A}} .
$$

The set of such $\dot\lambda$'s is closed under conjugation and compositions,
and by Cor. 32 of \cite{KLM}, the $Y$-matrix in (\ref{Ymatrix}) for
${\A}^\Gamma$ is non-degenerate, and we will denote the corresponding
genus $0$ modular matrices by $\dot S, \dot T$. We note that $d(\dot\lambda)$
is conjectured to be related to the asymptotic dimension of Kac-Wakimoto
in \cite{KW}, and one can find a precise statement of the conjecture 
and its consequences in \cite{L3} and in
\S2.3 of \cite{Xu9}.
Denote by $\dot\lambda$ 
(resp. $\mu$) the irreducible covariant representations of ${\A}^\Gamma$
(resp. ${\A}$) with finite index. Denote by
$b_{\mu\dot\lambda}\in {\mathbb N}\cup\{0\}$ 
the multiplicity of representation $\dot\lambda$ which appears
in the restriction of  representation $\mu$ when restricting from 
$ {\A}$ to ${\A}^{\Gamma}$. The $b_{\mu\dot\lambda} $ are also known as 
the {\em branching rules.}  
An irreducible covariant representation  $\dot\lambda$  of ${\A}^\Gamma$
is called
an {\it untwisted} representation if $b_{\mu\dot\lambda}\neq 0$ for 
some representation $\mu$ of $\A$.
These are
representations of ${\A}^\Gamma$ which appear as subrepresentations in the
the restriction of some representation of ${\A}$ to 
${\A}^\Gamma$. A  representation is called {\it twisted} if it is not
untwisted.    
Note that 
$\sum_{\dot\lambda} d(\dot\lambda) b_{\mu\dot\lambda} = d(\mu) |\Gamma'|$, and
$ b_{1\dot\lambda} =d(\dot\lambda)$. So we have
$$
\sum_{\dot\lambda \ \text{\rm untwisted} } d(\dot\lambda)^2 \leq
\sum_\mu({\sum_{\dot\lambda} d(\dot\lambda) b_{\mu\dot\lambda}})^2
= |\Gamma'|+ \sum_{\mu\neq 1} d(\mu)^2 |\Gamma'|^2
< |\Gamma'|^2+ \sum_{\mu\neq 1} d(\mu)^2 |\Gamma'|^2 = \mu_{{\A}^\Gamma}
$$ 
if $\Gamma'$ is not a trivial group, where in the last $=$ we have used 
Th.  \ref{Th.2.6}. 
It follows that the set of twisted representations of 
${\A}^\Gamma$ is not empty. 
This fact has already been observed in a special case in \cite{KLM} under
the assumption that ${\A}^{\Gamma}$ is strongly additive. Note that
this is very different from the case of
cosets, cf. \cite{Xu7} Cor. 3.2 where it was shown that under certain
conditions there are no twisted representations for the coset.\par
\subsection{Restriction to the real line: Solitons}
Denote by $\I_0$ the set of open, 
connected, non-empty, proper subsets of $\mathbb R$, thus $I\in\I_0$ 
iff $I$ is an open interval or half-line (by an interval of $\mathbb 
R$ we shall always mean a non-empty open bounded interval of $\mathbb 
R$).

Given a net $\A$ on $S^1$ we shall denote by $\A_0$ its restriction 
to $\mathbb R = S^1\setminus\{-1\}$. Thus $\A_0$ is an isotone map on
$\I_0$, that we call a \emph{net on $\mathbb R$}. In this paper
we denote by $J_0:=(0,\infty)\subset \mathbb R$.

A representation $\pi$ of $\A_0$ on a Hilbert space $\H$ is a map 
$I\in\I_0\mapsto\pi_I$ that associates to each $I\in\I_0$ a normal 
representation of $\A(I)$ on $B(\H)$ such that
\[
\pi_{\tilde I}\res\A(I)=\pi_I,\quad I\subset\tilde I, \quad 
I,\tilde I\in\I_0\ .
\]
A representation $\pi$ of $\A_0$ is also called a
\emph{soliton}.
As $\A_0$ satisfies half-line duality, namely 
$$
\A_0(-\infty,a)'= \A_0(a,\infty), a\in \mathbb R,
$$
by the usual DHR argument \cite{DHR} $\pi$ is unitarily equivalent to 
a representation $\rho$ which  acts identically on $\A_0(-\infty,0)$, thus
$\rho$ restricts to an endomorphism of $\A(J_0)= \A_0(0,\infty)$.
$\rho$ is said to be localized on $J_0$ and we also refer to $\rho$
as soliton endomorphism.

Clearly a representation $\pi$ of $\A$ restricts to a soliton 
$\pi_0$ of $\A_0$. But a representation $\pi_0$ of $\A_0$ does not 
necessarily extend to a representation of $\A$. 

If $\A$ is strongly additive, and  
a representation $\pi_0$ of $\A_0$ extends to a DHR representation of $\A$,
then it is easy to see that such an extension is unique, and in this 
case we will use the same notation $\pi_0$ to denote the corresponding
DHR  representation of $\A$.  
\subsection{A result on extensions of solitons}
The following proposition will play an important role in proving 
Th.\ref{permu}.
\begin{proposition}\label{normal}
Let $H_1,H_2$ be two subgroups of a compact group $\Gamma$
which acts properly on $\A$
, and let $\pi$ be a soliton of
$\A_0$. Assume that $\A$ is strongly additive. 
Suppose that $\pi\res \A^{H_i}, i=1,2$ are DHR representations.
Then $\pi\res (\A^{H_1}\vee \A^{H_2})$ is also a DHR representation, 
where $\A^{H_1}\vee \A^{H_2}$ is an intermediate net with 
$(\A^{H_1}\vee \A^{H_2})(I)=  \A^{H_1}(I)\vee \A^{H_2}(I), \forall I$.
\end{proposition}
\proof
Let $I$ be an arbitrary interval with $-1\in I$. 
It is sufficient to show that $\pi$ has a normal extension
to $\A^{H_1}(I)\vee \A^{H_2}(I).$
Since $\pi$ is 
a soliton, by choosing a unitary equivalence class of $\pi$ we may
assume that $\pi (x)=x,\forall x\in \A(I')$. 
Let $J\supset I$ be an interval sharing a boundary point with $I$
and let $I_0=J\cap I'$.
Since $\pi\res \A^{H_i}$
is a DHR representation, it is localizable on $I_0$. Denote the corresponding 
DHR representation localized on $I_0$ by $\pi_{i,I_0}$, then we can find
unitary $u_i$ such that 
$u_i \pi_{i,I_0} u_i^*= \pi$ on $\A^{H_i}$. It follows that 
$u_i\in \A^{H_i}(J)$ since $\pi$ is localized on $I$, and we have
$ \pi (x)= u_i x u_i^*, \forall x\in \A^{H_i}(I)$. Note that
$ \A^{H_1}(I) \cap   \A^{H_2}(I)\supset  \A^{\Gamma}(I)$, hence
$u_2^*u_1\in \A^{\Gamma}(I)'\cap \A(J)$. Since
$ \A^{\Gamma}\subset \A$ is a strongly additive pair (cf. \cite{Xu4}), 
it follows that $ \A^{\Gamma}(I)'\cap \A(J)= \A(I_0)$, and
$u_1 x u_1^* = u_2 x u_2^*, \forall x\in \A(I)$. Hence 
$\Ad_{u_1}$ defines a normal extension of $\pi$ from
$ \A^{H_1}(I)$ to $ \A^{H_1}(I)\vee \A^{H_2}(I)$. Such an extension
is also unique by definition. 
\endproof

\section{Induction and restriction for general orbifolds}

Let $\A$ be a M\"obius covariant net and $\B$ a subnet. Given  a 
bounded interval $I_0\in\I_0$ we fix canonical endomorphism 
$\gamma_{I_0}$ associated with $\B(I_0)\subset\A(I_0)$. Then we can choose 
for each $I\subset\I_0$ with $I\supset I_0$ a canonical endomorphism 
$\gamma_{I}$ of $\A(I)$ into $\B(I)$ in such a way that 
$\gamma_{I}\res\A(I_0)=\gamma_{I_0}$ and 
$\l_{I_1}$ is the identity on $\B(I_1)$ if 
$I_1\in\I_0$ is disjoint from $I_0$, where 
$\l_{I}\equiv\gamma_{I}\res\B(I)$.

We then have an endomorphism $\gamma$ of the $C^*$-algebra 
$\mathfrak A\equiv\overline{\cup_{I}\A(I)}$
($I$ bounded interval of $\mathbb R$).

Given a DHR endomorphism $\r$ of $\B$ localized in $I_0$, the 
$\a$-induction $\a_{\r}$ 
of $\r$ is the endomorphism of $\mathfrak A$ given by
\[
\a_{\r}\equiv \gamma^{-1}\cdot\Ad\e(\r,\l)\cdot\r\cdot\gamma\ ,
\]
where $\e$ denotes the right braiding unitary symmetry (there is 
another choice for $\a$ associated with the left braiding).
$\a_{\r}$ is localized in a right half-line containing $I_0$, namely 
$\a_\r$ is the identity on $\A(I)$ if $I$ is a bounded interval 
contained in the left complement of $I_0$ in $\mathbb R$. Up to 
unitarily equivalence, $\a_\r$ is localizable in any right half-line 
thus $\a_\r$ is normal on left half-lines, that is to say, for every 
$a\in\mathbb R$, 
$\a_\r$ is normal on the $C^*$-algebra 
$\mathfrak A(-\infty,a)\equiv\overline{\cup_{I\subset 
(-\infty,a)}\A(I)}$
($I$ bounded interval of $\mathbb R$), namely $\a_\r\res\mathfrak 
A(-\infty,a)$ extends to a normal morphism of $\A(-\infty,a)$. 
We have the following Prop. 3.1 of \cite{LX}:
\begin{proposition}\label{sg}
$\a_\r$ is a soliton endomorphism of $\A_0$.
\end{proposition}

\subsection{Solitons as endomorphisms}
Let $\A$ be a conformal net and $\Gamma$ a finite group acting properly 
on $\A$ (cf. (\ref{p'}). We will assume that $\A$ is strongly additive.
Let $\pi$ be an irreducible
soliton of $\A_0$ localized on $J_0=(0,\infty).$ Note that the 
restriction of $\pi$ to $\A(J_0)$ is an endomorphism and we denote this
restriction by $\pi$ when no confusion arises.  
Let $\pi_{\A^\Gamma}$ be a soliton of $\A^\Gamma_0$ localized on $J_0$
and unitarily equivalent to $\pi\res \A^\Gamma.$
Let $\rho_1$ be an endomorphism of $\A(J_0)$ such that 
$\rho_1(\A(J_0))= \A^\Gamma(J_0)$ and $\rho_1\bar\rho_1= \gamma$, 
where $\gamma$ is the canonical endomorphism from $\A(J_0)$ to 
$\A^\Gamma(J_0)$. Note that $[\gamma]=\bigoplus_{g\in \Gamma'}[g']$, 
where for simplicity we have used $[g]$ to denote the sector 
of  $\A(J_0)$
induced by the automorphism
$\beta_g$.  
By \cite{LR} as sectors of $\A^\Gamma(J_0)$ 
we have $[\pi_{\A^\Gamma}]= 
[\gamma \pi\res \A^\Gamma (J_0)]$.
\begin{definition}\label{gpi}
Define $\Gamma_\pi:=\{ h\in \Gamma| [h \pi h^{-1}]=[\pi] \}$. 
Note that ${\rm ker} V$ (cf. definition before 
(\ref{p'})) is a normal subgroup of $\Gamma_\pi$
and let $\Gamma_\pi':=\Gamma_\pi/{\rm ker} V$.
\end{definition}  
Note that 
\[
\Hom(\pi_{\A^\Gamma}, 
\pi_{\A^\Gamma })_{\A^\Gamma(J_0)}
\simeq \Hom(\bar\rho_1 \pi\rho_1, 
\bar\rho_1 \pi\rho_1)_{\A(J_0)}
.\]
By  Frobenius duality we have
\[
\langle \pi_{\A^\Gamma}, 
\pi_{\A^\Gamma}\rangle
=\langle \lambda, 
 \gamma\lambda\gamma\rangle  
.\]
\begin{lemma}\label{gh}
(1) If $g\neq h$, then $\langle \pi, g\pi h^{-1}\rangle =0$;\par
(2) $ \langle \pi, \gamma  \pi \gamma \rangle= |\Gamma_\pi'| 
= \langle \gamma\pi\res\A^\Gamma,\gamma\pi\res\A^\Gamma\rangle 
$
where $\Gamma_\pi'=\Gamma_\pi/{\rm ker} V;$\par
(3) 
$$
\langle \gamma\pi\res\A^\Gamma,\gamma\pi\res\A^\Gamma\rangle 
= \langle \gamma_1\pi\res\A^{\Gamma_\pi},\gamma_1\pi\res\A^{\Gamma_\pi}\rangle 
$$
where $\gamma_1$ is the canonical endomorphism from $\A(J_0)$ to 
$\A^{\Gamma_\pi}(J_0)$; \par
(4)
Every irreducible summand of $\pi\res\A^{\Gamma_\pi}_0$ (as soliton of 
$\A^{\Gamma_\pi}_0$) remains irreducible when restricting to 
$\A^{\Gamma}_0$. 
\end{lemma}
\proof
Note that $g\pi h^{-1}= g\pi g^{-1} g h^{-1}$, and $ g\pi g^{-1}$ is
a soliton equivalent to $\pi g^{-1}$ but localized on $J_0$. By
Lemma 8.5 of \cite{LX} 
we have proved (1). (2),(3) follows from (1) and the definition of
$\Gamma_\pi$. (4) follows from (3). 
\endproof

\begin{proposition}\label{solitonres}
Let $\pi_1,\pi_2$ be two irreducible solitons of $\A_0$  .  
If there is $g\in \Gamma'$ such that
$[\pi_1]=[ g \pi_2 g^{-1}]$, then
$[\gamma\pi_1\res\A^\Gamma]=[ \gamma\pi_2\res\A^\Gamma]$. Otherwise
$\langle \gamma\pi_1\res\A^\Gamma, \gamma\pi_2\res\A^\Gamma \rangle=0$.
\end{proposition}
\proof
By Frobenius duality  and Lemma 8.5 of \cite{LX} we have
$$
\langle \gamma\pi_1\res\A^\Gamma, \gamma\pi_2\res\A^\Gamma \rangle
=\sum_{g\in \Gamma'}\langle \pi_1, g\pi_2g^{-1} \rangle
$$
Hence $
\langle \gamma\pi_1\res\A^\Gamma, \gamma\pi_2\res\A^\Gamma \rangle=0$ 
if there is 
no $g\in \Gamma'$ such that
$[\pi_1]=[ g \pi_2 g^{-1}]$.  
If there is $g\in \Gamma$ such that
$[\pi_1]=[ g \pi_2 g^{-1}]$, then 
$$
\langle \gamma\pi_1\res\A^\Gamma, \gamma\pi_2\res\A^\Gamma \rangle
=\sum_{h\in \Gamma_{\pi_1}'}\langle 
\pi_1, h g\pi_2 g^{-1} h^{-1} \rangle
=  |\Gamma_{\pi_1}'|
$$
By exchanging $\pi_1$ and $\pi_2$ we get
$$
\langle \gamma\pi_1\res\A^\Gamma, \gamma\pi_2\res\A^\Gamma \rangle= 
\langle \gamma\pi_1\res\A^\Gamma, \gamma\pi_1\res\A^\Gamma \rangle=  
\langle \gamma\pi_2\res\A^\Gamma, \gamma\pi_2\res\A^\Gamma \rangle
$$
It follows that $[\gamma\pi_1\res\A^\Gamma]=  [\gamma\pi_2\res\A^\Gamma]$.
\endproof
\begin{theorem}\label{solitonind}
Assume that $\pi$ is irreducible with finite index and
$[\beta]= [\gamma \pi\res\A^\Gamma] = \bigoplus_j m_j [\beta_j]$. Then
$[\alpha_{\beta_j}]= m_j (\bigoplus_i [h_i \pi h_i^{-1}])$
where $h_i$ are representatives of $\Gamma/ \Gamma_\pi$.
In particular $d(\beta_j)= m_j d(\pi) \frac{|\Gamma|}{|\Gamma_\pi|}$,
and 
$\sum_j d(\beta_j)^2=  \frac{|\Gamma|^2}{|\Gamma_\pi|}  d(\pi)^2.$
\end{theorem} 
\proof
By the definition we have $
[\gamma \alpha_{\beta}]= [\beta \gamma]=[ \gamma \pi\gamma]$.
So we have 
$
\langle \gamma \alpha_{\beta}, \pi\rangle= \langle  \gamma \pi\gamma,\pi
\rangle=|\Gamma_\pi|.$
By Lemma 8.5 of \cite{LX} we have 
$\langle \gamma \alpha_{\beta}, \pi\rangle= \langle \alpha_{\beta}, \pi\rangle
$, and therefore
$\alpha_{\beta}\succ |\Gamma_\pi|\pi$. By Lemma 8.1 of \cite{LX} we have
$[h_i\alpha_{\beta}h_i^{-1}]= [\alpha_{\beta}]$, so
$ [\alpha_{\beta}] \succ |\Gamma_\pi'| \bigoplus_i [h_i \pi h_i^{-1}]$. 
On the other hand $d(\alpha_{\beta})=d(\beta)= |\Gamma'| d(\pi) = 
|\Gamma_\pi'|\sum_i d (h_i \pi h_i^{-1})$. It follows that
$$ [\alpha_{\beta}] = |\Gamma_\pi'| \big( \bigoplus_i [h_i \pi 
h_i^{-1}]\big).$$
Note that by Lemma 8.1 of \cite{LX} $[h_i^{-1}\alpha_{\beta_j} h_i]
=[\alpha_{\beta_j}]$, hence 
$$\langle \alpha_{\beta_j}, h_i \pi h_i^{-1}\rangle =
\langle h_i^{-1}\alpha_{\beta_j} h_i,  \pi \rangle
=\langle \alpha_{\beta_j},  \pi \rangle.$$
So we must have $[\alpha_{\beta_j}]= k_j (\bigoplus_i [h_i \pi 
h_i^{-1}])$
for some positive integer $k_j$. We note that
$k_j=\langle \alpha_{\beta_j}, \pi\rangle \leq \langle \beta_j, 
\gamma \pi\res\A^\Gamma\rangle = m_j$ by definitions and Frobenius duality. 
On the other hand
$\sum_j m_jk_j = |\Gamma_\pi'|= \sum_j m_j^2$, and we conclude that
$k_j=m_j$. Since by definition $\frac{|\Gamma|}{|\Gamma_\pi|}= 
\frac{|\Gamma'|}{|\Gamma_\pi'|}$, the proof of the theorem follows.
\endproof
\subsection{Solitons as representations}
In this section we use $\hat{\pi}$ to denote  an irreducible soliton of 
$\A_0$ on a Hilbert space $\H_\pi$. Let $\pi$ be a soliton 
unitarily equivalent to $\hat\pi$ but localized on $J_0$ as in the 
previous section. 
The restriction of $\hat{\pi}$ to $\A^\Gamma_0$, denoted by
$\hat{\pi}\res \A^\Gamma_0$ is also a soliton. Define 
$\Hom(\hat{\pi}\res \A^\Gamma_0,\hat{\pi}\res \A^\Gamma_0):=
\{ x\in B(\H_\pi)|   
x \hat{\pi}(a)= \hat{\pi}(a) x, \forall x\in \A^\Gamma_0 \},$ and
let $\langle \hat{\pi}\res \A^\Gamma_0,\hat{\pi}\res \A^\Gamma_0\rangle =
\dim \Hom(\hat{\pi}\res \A^\Gamma_0,\hat{\pi}\res \A^\Gamma_0)$.
\begin{lemma}\label{lg}
(1) 
$$
\langle \hat{\pi}\res \A^\Gamma_0,\hat{\pi}\res \A^\Gamma_0\rangle
=\langle \gamma\pi\res \A^\Gamma_0,\gamma\pi\res \A^\Gamma_0\rangle 
;$$\par
(2) $h\in \Gamma_\pi$ if and only if $\hat{\pi}\cdot\Ad_h\simeq \hat{\pi}$ as 
representations of $\A_0$.
\end{lemma}
\proof
By \cite{LR} $\hat{\pi}\res \A^\Gamma_0$ and 
$\gamma\pi\res \A^\Gamma_0$ are
unitarily equivalent as solitons of $\A^\Gamma_0$. Note that 
$\gamma\pi\res \A^\Gamma_0$ is localized on $J_0$,   and (1) follows
directly. As for (2), we note that 
$h^{-1} \pi h$ is localized on $J_0$ and unitarily equivalent to
$\hat{\pi}\cdot \Ad_h$, and (2) now follows from the definition (\ref{gpi}).
\endproof
From (2) of Lemma \ref{lg} we have for any $h\in \Gamma_\pi'$, 
there is a unitary operator denoted by $ \hat{\pi}(h)$ on $\H_\pi$ such that
$\Ad_{\hat{\pi}(h)^*}\cdot\hat{\pi} =\hat{\pi}\cdot\Ad_h$ as solitons of $\A_0.$ 
Since $\hat{\pi}$ is irreducible, there is 
a $U(1)$ valued cocycle $c_\pi(h_1,h_2)$ on $\Gamma_\pi$ such that 
$\hat{\pi}(h_1)\hat{\pi} (h_2) =c_\pi(h_1,h_2) \hat{\pi}(h_1 h_2)$. 
We note that $c_\pi(h_1,h_2)$ is fixed up to coboundaries (cf. \S2 of
\cite{Kar}).
Hence $h\rightarrow \hat{\pi}(h)$
is a projective unitary representation of $\Gamma_\pi'$ on $\H_\pi$
with cocycle $c_\pi$. Assume that 
$$\H_\pi= \bigoplus_{\sigma\in E} M_\sigma\otimes
V_\sigma
$$
where $E$ is a subset of irreducible projective representations of 
 $\Gamma_\pi'$ with cocycle $c_\pi$, and $M_\sigma$ is the multiplicity
space of representation $V_\sigma$ of  $\Gamma_\pi$. Then by definition 
each $M_\sigma$ is a representation of $\A^\Gamma_0$.   
\begin{lemma}\label{fullspectrum1} 
Fix an interval $I$. 
Assume that  
$\hat{\pi}$ is  a representation of $\A(I)$ (resp. a projective representation 
of $\Gamma'$ with cocycle $c_\pi$) 
on a Hilbert space $\H$ such that 
$\hat{\pi} (\beta_h(x))= \hat{\pi}(h)\hat{\pi} (x) \hat{\pi}(h)^*, 
\forall x\in \A(I)$. 
Let $ \sigma_1 \in \hat \Gamma'$ where $\hat  \Gamma'$ denotes the 
set of irreducible representations of $\Gamma'$, 
and  $\sigma_2$ be an irreducible
summand of the representation $\hat{\pi} $ of $\Gamma'$. Then:\par
(1) any irreducible
summand $\sigma$ of  $\sigma_1 \otimes \sigma_2$ appears as an irreducible
summand in the  projective representation $ \hat\pi$ of $\Gamma'$
with cocycle $c_\pi$.  In particular if
$\sigma_2$ is the trivial representation of $\Gamma'$ then all elements of
$\hat \Gamma'$ appear as irreducible
summand of the representation $ \pi$ of $\Gamma'$.\par
(2) Every irreducible projective representation of $\Gamma'$ with 
cocycle $ c_\pi$ appear as an irreducible summand of $\hat{\pi} $, and
$$
\sum_{\sigma, \sigma \ {\rm has} \  {\rm cocycle} \ c_\pi} 
\dim(\sigma)^2 = |\Gamma'|^2.
$$
\end{lemma}
\proof
Ad(1):
Since the action of $\Gamma'$ on $\A$ is proper, and $
\A^\Gamma(I)$ is a type III factor, for any $\sigma_1 \in \hat \Gamma'$,
by Page 48 of \cite{ILP} we can find a basis $V(\sigma_1)_i, 1\leq i\leq 
{\rm dim} \sigma_1$ in $\A(I)$
such that $V(\sigma_1)_i^* V(\sigma_1)_j=\delta_{ij}, $ 
and the linear span of 
$V(\sigma_1)_i, 1\leq i\leq 
{\rm dim}\sigma_1 $ forms the irreducible representation $\sigma_1$ of 
$\Gamma'$. Let $W(\sigma_2)_i \in \H, 1\leq i\leq 
{\rm dim} \sigma_2$ be an orthogonal basis of representation
$\sigma_2$. We claim that the vectors $ \pi( V(\sigma_1)_i)W(\sigma_2)_j,
1\leq i\leq 
{\rm dim} \sigma_1,  1\leq j\leq 
{\rm dim} \sigma_2$ in $ \H$ are linearly independent. If 
$\sum_{ij} C_{ij}\pi( V(\sigma_1)_i)W(\sigma_2)_j=0$ for some complex
numbers $C_{ij}$, multiply both sides by $\hat{\pi} ( V(\sigma_1)_i^*)$ and use the
orthogonal property of $ V(\sigma_1)_j$'s above we have
$\sum_{j} C_{ij}W(\sigma_2)_j=0$, and hence $C_{ij}=0$ since
$W(\sigma_2)_j$'s are linearly independent. It follows that the 
linear span of $\hat{\pi}( V(\sigma_1)_i)W(\sigma_2)_j, 
1\leq i\leq 
{\rm dim} \sigma_1,  1\leq j\leq 
{\rm dim} \sigma_2$ gives a tensor product representation of 
$\Gamma'$ on a subspace of $\H$, and the lemma follows.\par
Ad (2):
Let $\sigma_3$ be an irreducible summand of $\pi$, and let $\sigma_4$
be an arbitrary irreducible projective representation of $\Gamma'$ with 
cocycle $c_\pi$. By definition $\bar\sigma_3\otimes \sigma_4$ is 
a representation of $\Gamma'$ ($\bar \sigma_3$ stands for the conjugate
of $\sigma_3$), and hence $\bar\sigma_3\otimes  \sigma_4\succ\sigma_5$
for some $\sigma_5\in \hat{\Gamma'}$, and it follows that 
$\sigma_4$ appears as an irreducible summand of $\sigma_3\otimes  \sigma_5$,
and so by (2) every irreducible projective representation of 
$\Gamma'$ with cocycle $c_\pi$ appears as an irreducible summand of 
$\pi$. Note the twisted group algebra $\mathbb C^{c_\pi}[\Gamma']$
with cocycle $c_\pi$
(cf. P. 85 of \cite{Kar}) is semisimple, and the equality in (2) follows. 
\endproof
\begin{theorem}\label{projrep}
(1) $\Hom(\hat{\pi}\res\A^\Gamma,\hat{\pi}\res\A^\Gamma)=\bigoplus_{\sigma\in E} 
{\rm Mat}(\dim 
(\sigma)) $ where $E$ is the set of irreducible projective 
representations of $\Gamma_\pi'$ with the  cocycle $c_\pi$; \par   
(2) $\sum_{\sigma\in E} \dim (\sigma)^2 = |\Gamma_\pi'|$;\par
(3) $M_\sigma$ as defined before Lemma \ref{fullspectrum1} 
is an irreducible representation of 
$\A^{\Gamma}_0$, and $M_\sigma$ is not unitarily equivalent to $M_{\sigma'}$
if $\sigma\neq \sigma'.$ 
\end{theorem} 
\proof
(1) (2) follows directly from Lemma \ref{fullspectrum1}. As for (3),
note that by  (2) of Lemma \ref{gh} and (1) of Lemma \ref{lg} we have
$  \langle \hat{\pi}\res \A^{\Gamma}, 
\hat{\pi}\res \A^{\Gamma}\rangle
=|\Gamma'|.
$
On the other hand 
$\langle\pi\res \A^{\Gamma}, 
\pi\res \A^{\Gamma}\rangle \geq \sum_{\sigma\in E} {\rm dim} V_\sigma^2
$
with equality iff  $M_\sigma$ as above is an irreducible representation of 
$\A^{\Gamma}_0$, and $M_\sigma$ is not unitarily equivalent to $M_{\sigma'}$
if $\sigma\neq \sigma'.$ Since we have equality by (2),
(3) is proved.
\endproof
Since for cyclic group
$H^2( {\mathbb Z}_k, U(1))=0$, we have proved the following corollary which
generalizes Lemma 2.1 of  \cite{Xu5}. 
\begin{corollary}\label{cyclic1}
If $\Gamma_\pi' = {\mathbb Z}_k$ for some positive integer $k$, then 
$\Hom(\hat{\pi}\res \A^\Gamma, \hat{\pi}\res \A^\Gamma)$ is isomorphic
to the group algebra of ${\mathbb Z}_k$, and $ \hat{\pi}\res \A^\Gamma$ 
decomposes into
$k$ distinct irreducible pieces.
\end{corollary} 
\section{Solitons in Affine Orbifold}
\subsection{Conformal nets associated with the affine algebras}
Let $G$ be a compact Lie group of the form $G:=G^0\times G^1\times\cdots
\times G^s$ where $G^0=U(1)^r$, 
and $G^j, j=1,...,s,$ are simple 
simply-connected groups. Let $\g^j$ denoted the Lie algebra of 
$G^j,  j=0,...,s$ and let $L:=\{ \omega\in \g^0 | e^{2\pi i \omega}=1 \}.$
Note that $G^0=U(1)^r= {\mathbb R}^r/L$. 
We assume that $\g =\bigoplus_j \g^j$ is equipped with a symmetric  
even negative definite invariant bilinear form. This means that 
the length square of any $\omega\in i\g^j$($j=0,...,s$) such that
$e^{2\pi i \omega}=1$ is an even integer. Note that our condition
on the bilinear form is slightly stronger than the condition on P. 61
of \cite{KT} to ensure locality of our nets (cf. Remark 1.1 of 
\cite{KT}). When restricted to a simple $\g^j$, the even property 
means that the bilinear form is equal to $k_j(v|v'),$ where 
$k_j\in{\mathbb N}$
will be identified with the level of the affine Kac-Moody algebra
$\hat\g^j$ and $$
(v|v')= \frac{1}{2 g_j^{\vee}}{{\rm Tr}_{\g^j}(\Ad_v\Ad_{v'})}
$$
($g_j^{\vee}$ is the dual Coxeter number of $\g^j$). We will
fix $k_0=1$. \par
We will denote by $\tilde LG$ the central extension of $LG$ whose
Lie algebra is the (smooth) affine Kac-Moody algebra $\hat\g$. 
For an interval $I\subset S^1$, we denote by 
$\tilde L_IG:\{ f\in \tilde LG| f(t)=e, \forall t\in I'$
where $e$ id the identity element in $G$, and
$\tilde L_I\g:\{ p\in \tilde L_\g| p(t)=0, \forall t\in I'.$
We will write elements of $\tilde L\g$ as $(f,c)$ where 
$f\in L\g, c\in \mathbb C$ and $(0,c)$ is in the center of $\tilde L\g$. 
Denote by $\A_{G_k}$ the conformal net associated with representations of
$\tilde LG$ at level $k=(k_0,...,k_s)$. The following lemma follows from
\cite{TL}: 
\begin{lemma}\label{strong}
$\A_{G_k}$ is strongly additive. 
\end{lemma}
\par
For simplicity we will denote $\A_{G_k}$ by $\A$ in this chapter. 
Let $Z_j\subset G^j$ denote the center of $G^j, j=1,...,s,$ and let
$Z^0=L^*/L$ where $L^*:= \{\mu\in \g^0|(\mu|\omega)\in {\mathbb Z}, \forall
\omega \in L\}$. 

The following finite subgroup of $G$ will 
play an important role:
$$
Z(G):= Z^0\times Z^1\times\cdots\times Z^s
$$
Recall from \S4.2 of \cite{KT} that an element $g\in G$ is called 
{\bf non-exceptional} if there exists $\beta(g)\in i\g$ such that
$g=e^{2\pi i\beta(g)}$ and the centralizer $G_g:=\{ b'\in G| b'g{b'}^{-1}=g \}
$ of 
$g$ is the same as 
$G_{\beta(g)}: \{ b'\in G| b'\beta(g) {b'}^{-1}=\beta(g)\}$, 
the centralizer of   
$\beta(g)$. \par
Let $\Gamma$ be a finite subgroup of $G$. Then it follows by 
definition that $\Gamma$ acts properly on $\A$. We will be interested in
the irreducible representations of $\A^\Gamma$. 
Note that $Z(G)$ acts on $\A$ trivially.
Hence $\A^\Gamma=\A^{\langle \Gamma,Z(G)\rangle}$ where 
$\langle \Gamma,Z(G)\rangle$ is the subgroup of $G$ 
generated by $\Gamma, Z(G).$
Without losing generality, we will always assume that 
$\Gamma\supset Z(G)$. By definition before (\ref{p'}) we have 
$\Gamma'=\Gamma/Z(G).$ \par 
The following definition is definition 4.1 of \cite{KT}:
\begin{definition}\label{non-excep}
A group $\Gamma$ is called a {\bf non-exceptional} subgroup of
$G$ if for any $g\in \Gamma$ there exists $\zeta \in Z(G)$ such that
$\zeta g$ is a non-exceptional element.
\end{definition}
Recall from \cite{KT} that every element of $Z$ can be written in the form 
$$
\zeta=(\zeta_{j_0}^{(0)},..., \zeta_{j_s}^{(s)}) \in 
Z_0\times\cdots\times Z_s,  \zeta_{j_\nu}^{(\nu)}= e^{2\pi i
\Lambda_j^{(\nu)}}.
$$
Here $\Lambda_j^{(0)}$ generate the finite abelian group $L^*/L$; 
for each 
simple component $\g$ the fundamental weight $\Lambda_j$ belongs to the set $J$ (1.33) of 
\cite{KT}. If both $g$ and $\zeta_j g$ are non-exceptional,
we can write
\begin{equation}\label{conn}
\beta(\zeta_j g) = \beta(g)+ \Lambda_j +m, 
[\beta(g), \beta(\zeta_j g)]=0, e^{2\pi im}=1.
\end{equation}
Now we define the action of $\zeta_j$ on $\Lambda$. 
By Lemma 4.1 of \cite{KT} the phase factor
\begin{equation}\label{phase}
\sigma_j(b')= e^{2\pi i(k \Lambda_j+km|\beta')}, b'=e^{2\pi i \beta'}\in 
\Gamma_g, [\beta', \Lambda_j+m]=0
\end{equation}
gives a 1-dimensional representation of $\sigma_j$ of $\Gamma_g$. The
transformation 
$\Lambda\rightarrow \zeta_j(\Lambda)$ of a lattice weight
$\Lambda\in L^*$ is given by $ \zeta_j(\Lambda)=(\Lambda+ \Lambda_j) {\rm mod}
L$. If $\g$ is a simple rank $l$ Lie algebra and $\Lambda$ is an 
integral weight at level $k$, then $\zeta_j(\Lambda):=
k\Lambda_j + w_j \Lambda$ where $w_j$ is the unique element of the 
Weyl group of $\g$ that permutes the 
set $\{ -\theta, \alpha_1,...,\alpha_l \}$ and satisfies 
$-w_j \theta = \alpha_j$. 
\begin{definition}\label{centeract}
\cite{KT}
For any $\zeta\in Z$, $\Lambda=\sum_\nu \Lambda^\nu$, we define:
$$\zeta(\Lambda)= 
\sum_\nu( w_{j_\nu}\Lambda^{\nu} + k_\nu\Lambda_{j_\nu})
.$$
We will use $\pi_\Lambda$ to denote the irreducible representations of
$\tilde LG$ on a Hilbert space $\H_\Lambda$ with highest weight $\Lambda$. 
\end{definition}
Note that $\pi_\Lambda$ gives an irreducible representation of $\A_{G_k}$
by \S3 of \cite{FG} on $\H_\Lambda$.
We will write $\zeta=e^{2\pi i \beta(\zeta)}$ with 
$\beta(\zeta)=(\beta(\zeta_{j_0}^{(0)}),...,\beta(\zeta_{j_s}^{(s)}))$
and $\beta(\zeta_{j_\nu}^{(\nu)})= \Lambda_j^{(\nu)} +m $
where $m$ is as in (\ref{phase}).
Let $P_g: [0,1]\rightarrow G$ be a map with $P_g(\theta)=
e^{ 2\pi i\beta(g) \theta}, 0\leq \theta\leq 1$, 
 $P_{\zeta g}: [0,1]\rightarrow G$ be a map with $P_{\zeta g}(\theta)=
e^{ 2\pi i\beta(\zeta g) \theta}, 0\leq \theta\leq 2\pi$,
and $P_\zeta: [0,1]\rightarrow G$ be a map with $P_\zeta(\theta)=
e^{ 2\pi i\beta(\zeta) \theta}, 0\leq \theta\leq 1$. We note that 
$\Ad_{P_\zeta}$ is an automorphism of $LG$ since $\zeta$ is in the center
of $G$.
\begin{lemma}\label{nexcp}
(1) If $g$ is non-exceptional then  $P_g \in Z(G_g)$;
\par
(2) If $\zeta g$, $g$ are non-exceptional  then
$P_{\zeta g} P_{g}^{-1} =P_\zeta.$
\end{lemma}
\proof
If $h\in G_g$, since $g$ is non-exceptional, 
it follows that 
$$ h e^{2\pi i \theta \beta(g)}h^{-1} =e^{2\pi i \theta \beta(g)},  
0\leq \theta \leq 1$$ and (1) is proved. \par 
Since $\zeta g$, $g$ are non-exceptional ,  
by (\ref{conn}) $[\beta(\zeta g), \beta]=0$ and 
(2) follows immediately.   
\endproof
\begin{lemma}\label{cocycle}
 If $\zeta g$, $g$ are non-exceptional in ,  
and with notations as above, we have:
\par
(1) $\Ad_{P_\zeta}$ lifts to an automorphism  
denoted by $Ad_\zeta$ of $\tilde LG$;\par
(2) The induced action of $\Ad_\zeta$ on $\tilde Lg$ is given by
$$
\Ad_\zeta (f,c) = (\Ad_{P_\zeta}.f, k(\zeta|f)+c)
;$$\par
(3) There is an unitary $U: \H_{\zeta(\Lambda)}\rightarrow H_\Lambda$ such that
$U^*\pi_{\zeta(\Lambda)} (Ad_\zeta)U=\pi_\Lambda $ 
as representations of $\tilde LG;$\par
(4) 
$ U^*\pi_{\zeta(\Lambda)}(h) \sigma_\zeta(h) U=
\pi_\Lambda(h) 
$ 
for any $h\in \Gamma_g$, where $\sigma_\zeta=\otimes_\nu \sigma_{j_\nu}$
with $\sigma_{j_\nu}$ as defined in (\ref{phase}). 

\end{lemma}
\proof
We note that the path $ P_{\zeta g}P_g^*$ is an element of $
L(G/Z(G))$. When $G$ is semisimple, (1),(2) follows from 
Lemma 4.6.5 and equation (4.6.4) of \cite{PS}. The proof in \S4.6 of \cite{PS}
also generalizes easily  to the proof of (1) and (2)  
when $G=G_0=U(1)^r$. 
As for (3), 
first note that  $\pi_{\zeta(\Lambda)} (  \Ad_\zeta)$ 
is an irreducible representation of $\tilde LG$,
since such irreducible representations are classified (cf. \cite{PS} and 
\cite{K1}), we just 
have to identify it with the known representations. By using 
Th. 4.2 of \cite{KT} for the special case when the group
$\Gamma$ is trivial, 
we conclude that the character  of 
 $\pi_{\zeta(\Lambda)}\cdot\Ad_\zeta$ 
is the same as that of $\pi_\Lambda (\tilde LG)$, and it follows that they
are unitarily equivalent as representations of $\tilde LG$. \par
For any $h=e^{2\pi i \beta'}\in G_g \subset \tilde LG$, by (2) we have
$$\pi_\Lambda (Ad_\zeta(h))= \pi_\Lambda (h) 
\prod_{\nu} e^{2\pi i k(\beta'|\Lambda_{j_\nu} +m)}
=\pi_\Lambda (h) \sigma_{\zeta} (h)
$$
Using  (3) we have 
$$
U^*\pi_{\zeta(\Lambda)} (h) \sigma_\zeta (h) U=\pi_\Lambda (h) 
$$ 
\endproof
\subsection{Constructions of solitons}
Let $\pi_\Lambda$ be an irreducible representation of $\tilde LG$ with
highest integral weight $\Lambda$. We will denote the net $\A_{G_k}$ simply
by $\A$ in this section.
For $g\in G$,
let $\beta(g)$ be an element in the Lie algebra of $G$ such that 
$e^{2\pi i \beta(g)}= g$. Define $P_g(\theta):= e^{2\pi i \theta \beta(g)},
0\leq \theta\leq 1.$ 
Identify $\mathbb R$ with the open interval $(0,1)$ via a smooth  map
$\phi:  (-\infty,+\infty)\rightarrow (0,1), \phi(t) = \frac{1}{\pi}( 
{\rm tan}^{-1}(t)+ \frac{\pi}{2})$. 
For any $I\subset \mathbb R$, 
Let $P_{g,I}\in L_I G$ be a loop localized on $I$ such that 
$P_{g,I} (t)= P_g (\phi(t)), \forall t\in I$.
\begin{definition}\label{affinesl}
For any $x\in \A(I)$, define
$\hat{\pi}_{\Lambda,g,I} (x):= \pi_\Lambda(P_{g,I} x P_{g,I}^*)$.
\end{definition}
We note that the above definition is independent of the choice of
$P_{g,I}$: if $\tilde P_{g,I}$ is another loop such that
$\tilde P_{g,I}(t)= P_{g,I}(t),\forall t\in I$, then
$ \tilde P_{g,I}(t)  P_{g,I}^{-1}$ is a loop with support in $I'$, and
so
$\pi_\Lambda(P_{g,I} x P_{g,I}^*)= 
\pi_\Lambda(\tilde P_{g,I} x \tilde P_{g,I}^*),  \forall x\in \A(I).$
One checks easily that definition (\ref{affinesl}) defines a soliton,
and we denote it by $\hat{\pi}_{\Lambda,g}$. \par
Fix $J_0:=(0,\infty)\subset \mathbb R$. To obtain a soliton equivalent to
$\hat{\pi}_{\Lambda,g}$ but localized on $J_0$, 
we choose a smooth path $P_g^{J_0}
\in C^\infty({\mathbb R}, G)$ which satisfies the following 
 boundary conditions: 
$P_g^{J_0}(t) = e, {\rm if } -\infty < t \leq 0 $ and
$P_g^{J_0}(t) = g, {\rm if } \ 1\leq t < \infty$. 
For any interval $I\subset {\mathbb 
R}$, we choose a loop $P^{J_0}_{g,I} \in LG$ such that
$P^{J_0}_{g,I}(t)= P^{J_0}_g(t), \forall t\in I$.
\begin{definition}\label{affinesll}
For any $x\in\A(I)$, define $\pi_{\Lambda,g,I}:= \Lambda(
 P^{J_0}_{g,I}x {P^{J_0}_{g,I}}^*)$ where we use $\Lambda$ to denote a 
representation unitarily equivalent to $\pi_\Lambda$ but  localized on $J_0$.
\end{definition}
We denote the soliton in the above definition as $\pi_{\Lambda,g}$.
\begin{proposition}
The unitary equivalence class of  $\pi_{\Lambda,g}$ is 
independent of the choice of the path $P_g^{J_0}$ as long as it satisfies
the boundary conditions given as above, and  $\pi_{\Lambda,g}$ is
localized on $J_0$. Moreover $\pi_{\Lambda,g}$ is unitarily equivalent 
to $\hat{\pi}_{\Lambda,g}$, and 
$\pi_{\Lambda,g}$ restricts to a DHR representation
of $\A^{\langle g\rangle}$ where $\langle g\rangle$ denotes the closed
subgroup of $G$ generated by $g$.
\end{proposition}
\proof
If  $\tilde P_g^{J_0}$ is another path which satisfies the same boundary
condition as  $P_g^{J_0}$, then $\tilde P_g^{J_0} {P_g^{J_0}}^{-1} \in LG,$
and the first statement of the proposition follows by definition.  By 
definition $\pi_{\Lambda, g, J_0'} (x)=x, \forall x\in \A(J_0')$ since 
$ P_g^{J_0}(t)=e,$ if $ -\infty < t \leq 0$, and so  $\pi_{\Lambda,g}$ is
localized on $J_0$. Since $ P_g^{J_0}P_g^{-1}$ extends to an element in $LG$,
it follows that   $\pi_{\Lambda,g}$ is unitarily equivalent 
to $\hat{\pi}_{\Lambda,g}$. 
To prove the last statement, let $I$ be an interval 
with $-1\in I$. It is sufficient to show that   $\pi_{\Lambda,g}$
has a normal entension to $\A^{\langle g \rangle}(I)$.  Recall from \S3.6
that we identify $\mathbb R= S^1 \setminus \{-1\}$ and $J_0=(0,\infty)\subset 
\mathbb R$. Since the  net $\A$ is strongly additive by Lemma \ref{strong},
and so $\A^{\langle g \rangle}$ is strongly additive by \cite{Xu4},  
we can assume that 
$\A^{\langle g \rangle} (I)= \A^{\langle g \rangle}(-\infty, a)\vee 
\A^{\langle g \rangle} (b, \infty)$ where $a<b$. 
Let us assume that $P^{J_0}_{g, (-\infty,a)}$ and 
$P^{J_0}_{g, (b,\infty)}$ are the elements in $LG$ such that 
$\Ad_{\Lambda(P^{J_0}_{g, (-\infty,a)})}= \pi_{\Lambda,g, (-\infty,a)}$ and 
$\Ad_{\Lambda(P^{J_0}_{g, (b,\infty)})}= \pi_{\Lambda,g, (b,\infty)}$ as in Definition
\ref{affinesll}. Choose an element $\tilde P\in LG$ so that 
$\tilde P(t)= g P^{J_0}_{g, (-\infty,a)} (t),  -\infty < t< a$
and $ \tilde P(t)= P^{J_0}_{g, (b,\infty)} (t),  b < t< \infty$. 
Then by definition $\Ad{\Lambda(\tilde P)} (x)= \pi_{\Lambda,g}(x), \forall
x\in \A^{\langle g \rangle}(-\infty, a)\vee 
\A^{\langle g \rangle} (b, \infty)$, and hence $\Ad\Lambda({\tilde P})$ 
defines the
normal extension of $\pi_{\Lambda,g}$ to $\A^{\langle g \rangle} (I)$.
\endproof 

\begin{proposition}\label{afusion}
As sectors of $\A(J_0)$ we have: \par
(1) $[\pi_{\Lambda,g}]= [\Lambda \pi_{1,g}];$\par
(2) $[\pi_{1,g_1} \pi_{1,g_2}]=  [\pi_{1,g_1 g_2}], 
[h \pi_{\Lambda,g} h^{-1}]= [\pi_{\Lambda,hgh^{-1}}];$\par
(3) Assume that $\Lambda,\mu$ are irreducible DHR representations of $\A.$ 
Then
$\langle \Lambda, \mu\pi_{1,g} h\rangle =1$ if and only if
$h\in Z(G), g\in Z(g)$ and $\Lambda= g^{-1}(\mu)$ 
where the action of the center is as in (\ref{centeract}).
In all other cases
$ \langle \Lambda, \mu\pi_{1,g} h\rangle =0$;\par
(4) If $\Lambda_1,\Lambda_2$ are irreducible DHR representations
of $\A$ , then
$\langle \pi_{\Lambda_1,g_1}, \pi_{\Lambda_2,g_2}h\rangle
=1$ if and only if 
$h\in Z(G)$ and there exists a $g\in Z(G)$ such that 
$g_2=gg_1$ and $\Lambda_2= g^{-1}(\Lambda_1).$ In all other cases
$ \langle \pi_{\Lambda_1,g_1}, \pi_{\Lambda_2,g_2}h  \rangle=0;$\par
(5) The stabilizer $\Gamma_{\Lambda,g} $of $\pi_{\Lambda,g}$ 
(cf. (\ref{gpi})) is given by $\Gamma_{\Lambda,g}=\{ h\in \Gamma
   | hgh^{-1}= g_1g, g_1(\Lambda)=\Lambda, g_1\in Z(G) \}.$
\end{proposition}
\proof
(1) and (2) follows directly from definition \ref{affinesll}. Now
assume that $\langle \Lambda, \mu\pi_{1,g} h\rangle =1$. By lemma 8.5
of \cite{LX} we conclude that $[h]=[1]$ and so $h\in Z(G)$, hence
$[\Lambda]=[ \mu\pi_{1,g}]$, and it follows that $ \mu\pi_{1,g}$ is
a DHR representation of the net $\A$. In particular $\mu\pi_{1,g}$ is
normal on $\A(-\infty,0)\vee \A(1,\infty)$. Choose $\mu$ to be localized on
$\A(0,1)$. Since 
$\A(-\infty,0)\vee \A(1,\infty)$ is a type III von Neumann algebra, there
is a unitary $u$ such that 
$\pi_{1,g}(x)= uxu^*, \forall x\in \A(-\infty,0)\vee \A(1,\infty)$. 
Since $ \pi_{1,g}=id$ on $\A(-\infty,0)$ and  $ \pi_{1,g}= {\rm Ad}_g$ on 
$\A(1,\infty)$, we have $u\in \A(-\infty,0)'\cap  \A^{G}(1,\infty)'$. 
By (2) of Lemma 3.6 in \cite{Xu4} the
pair $\A^G \subset \A$ is strongly additive (cf. Definition 3.2 of 
\cite{Xu4} ) since $\A$
is strongly additive by Lemma \ref{strong}, 
and so  $ \A(-\infty,0)'\cap  \A^{\Gamma}(1,\infty)'
= \A(0,1).$ Therefore
$u\in \A(0,1), \Ad_g (x)=x, \forall x\in \A(1,\infty),$ and so 
$g\in Z(G)$. Hence we have $\langle \Lambda, \mu\pi_{1,g}\rangle =1$.
By (3) of Lemma \ref{cocycle} and definition of $\pi_{1,g}$ we have 
$\Lambda=g^{-1}(\mu)$ where the action of the center is defined in 
(\ref{centeract}).\par
As for (4), by (1) and (2) we have
\begin{align}
\langle \pi_{\Lambda_1,g_1}, \pi_{\Lambda_2,g_2}h\rangle = 
\langle \Lambda_1\pi_{1,g_1}, \Lambda_2\pi_{1,g_2}h\rangle 
= \langle\Lambda_1, \Lambda_2 \pi_{1,g_2}h \pi_{\Lambda_1,g_1^{-1}}
\rangle \\
=\langle \Lambda_1, \Lambda_2\pi_{1,g_2 h g_1^{-1} h^{-1}} h
\rangle
\end{align}
and (4) follows from the above equation and (3). (5) follows from
definitions and (4). 
\endproof
\subsubsection{Comparing solitons with ``twisted representations''}
Let $e^{2\pi i \beta}=g$ and choose the Cartan subalgebra of $\g$ which
contains $\beta$. 
In the definition (\ref{affinesl}), if we choose $x=\pi_1(y), y\in \tilde 
L_IG$, 
then ${\hat{\pi}}_{\Lambda,g,I}(\pi_1(y))= 
\pi_\Lambda(P_{g,I}y P_{g,I}^*)$. Note that $\Ad_{P_{g,I}}$ is
an automorphism of $\tilde L_IG$, and induces an automorphism on
$\tilde L_I \g$. By Prop. 4.3.2 of \cite{PS}, if we write 
elements of $\tilde L_I \g$ as $(f,c)$, where $f\in C^\infty(S^1,\g)$
with support in $I$, and $c\in \mathbb C$, then
\begin{equation}\label{compare}
\Ad_{P_{g,I}}(f,c)= (\Ad_{P_{g,I}}.f, c+k(\beta|f))
\end{equation}
Let us  check  that  (\ref{compare}) agrees with the definition 
of twisted representation 2.11-2.14
of \cite{KT} on $\tilde L_I \g, \forall I\subset \mathbb R.$ 
Let $E^\alpha$ be a raising or lowering operator as on Page 64 of
\cite{KT}. Let $f_1\in C^\infty(S^1,\mathbb R)$ be a smooth map such that
$f_1(t)=0, \forall t\in I'$. By the commutation relation
$[E^\alpha, \beta]= -(\alpha|\beta) E^\alpha$
we have $\Ad_{P_{g,I}}.f= z^{-(\alpha|\beta)}E^\alpha f_1$ 
where $z^{-(\alpha|\beta)}:= e^{-2\pi i \theta (\alpha|\beta)}$ as 
a function on $[0,1]$, and $(\beta|f_1E^\alpha)=0$ by definition. 
By (\ref{compare}) we have
$$
\Ad_{P_{g,I}}(f_1 E^\alpha ,c)= (z^{-(\alpha|\beta)}E^\alpha f_1 , c)
$$
which is the restriction of (2.11) of \cite{KT} to
 $\tilde L_I \g.$ Similarly one can check that 
 (\ref{compare}) agrees with the definition 
of twisted representation 2.12-2.14
of \cite{KT} on $\tilde L_I \g, \forall I\subset \mathbb R.$  
Hence our soliton representations
in Definition \ref{affinesl} can be regarded as ``exponentiated'' version of
the twisted representations in \S2 of \cite{KT}. In the next section we 
shall see that these soliton representations are important in constructing
 irreducible DHR representations of $\A^\Gamma$. Motivated by the
above observations, we have the following conjecture:
\begin{conjecture}\label{conjecture}
There is a natural one to one correspondence between the 
set of  irreducible DHR representations of $\A^\Gamma$ and the
set of  irreducible representations of the orbifold chiral algebra 
as defined on Page 74 of \cite{KT} with gauge group $\Gamma$.
\end{conjecture}
We note that this conjecture, together with the results of \S5.4 and 
\S5.5, give a prediction on the  the
set of  irreducible representations of the orbifold chiral algebra 
as defined on Page 74 of \cite{KT} with non-exceptional 
gauge group $\Gamma$.  
\subsection{Completely rational case}
Assume that the net $\A$ associated to $G$ has the property that
\begin{equation}\label{ra}
\mu_\A= \sum_\Lambda d(\Lambda)^2
\end{equation}
where the sum is over all 
irreducible projective representations of $LG$ of a fixed level. 
When $G=SU(N)$ this property is proved by \cite{Xu6}. 
We show that all irreducible DHR representations
of $\A^\Gamma$ are  obtained from decomposing the 
restriction of solitons $\pi_{\Lambda,g}$ to $\A^\Gamma$, answering one
of the motivating questions for this paper. 
By Prop. \ref{solitonres}
$\pi_{\Lambda_1,g_1}\res\A^\Gamma\simeq \pi_{\Lambda_2,g_2}\res\A^\Gamma$ 
iff there exists
$h\in \Gamma$ such that $ [h \pi_{\Lambda_1,g_1} h^{-1}]
= [\pi_{\Lambda_2,g_2}.$ By (2) and (4) of Prop. \ref{afusion}
this is true if there is a $g_3\in Z(G)$ such that 
$\Lambda_2= g_3^{-1} (\Lambda_1)$
and $g_2= h g_3g_1h^{-1}$. Define an action of group $Z(G)\times \Gamma$
on the set $(\Lambda,g)$ by $(g_3, h).(\Lambda,g)= ({g_3}^{-1}(\Lambda),
h g_3g_1h^{-1})$. 
Denote the orbit of $(\Lambda,g)$ by $\{ \Lambda,g\}$. 
Note that the stabilizer of  $(\Lambda,g)$ has the same order as
the stabilizer $\Gamma_{\Lambda,g}$ of $\pi_{\Lambda,g}$ by (5) of Prop.
\ref{afusion}.  
Hence 
the orbit   
$\{ \Lambda,g\}$ contains 
$\frac{|Z(G)\times \Gamma|}{|\Gamma_{\Lambda,g}|}$ elements.  
Let $[\gamma\pi_{\Lambda,g}]= \sum_{i} m_i[\beta_i]$ where $\beta_i$
are irreducible DHR representations of $\A^\Gamma$.

By Th. \ref{solitonind} 
$\sum_{i} d(\beta_i)^2 = \frac{|\Gamma'|^2}{ |\Gamma_{\Lambda,g}'|}
d(\Lambda)^2$. By Prop. \ref{solitonres} 
we get the sum of index of all different
irreducible DHR representations  of $\A^\Gamma$ 
coming from decomposing the restriction 
of $\pi_{\Lambda,g}$ to  $\A^\Gamma$ is given by
$$
\sum_{\{\Lambda,g\}}  \frac{|\Gamma'|^2}{ |\Gamma_{\Lambda,g}'|}
d(\Lambda)^2. 
$$
Since  the orbit   
$\{ \Lambda,g\}$ contains $\frac{|Z(G)\times \Gamma|}{|\Gamma_{\Lambda,g}|}$ 
elements, the above sum is equal to
$$
\sum_{\Lambda,g}   \frac{|\Gamma'|^2}{ |\Gamma|}d(\Lambda)^2
= |\Gamma'|^2 \mu_{\A} =\mu_{\A^\Gamma}
$$
where in the last $=$ we have used Th. \ref{orb}. By Th. 33
of \cite{KLM}  we have proved the following:
\begin{theorem}\label{allaffine}
If  equation (\ref{ra}) holds, then 
every irreducible DHR representation of $\A^\Gamma$ is contained in
the restriction of $\pi_{\Lambda,g}$ to  $\A^\Gamma$ for some 
$\Lambda,g$ where $\pi_{\Lambda,g}$ is defined as in (\ref{affinesll}).
\end{theorem}
Let $G=SU(N_1)\times SU(N_2)\times\cdots\times SU(N_m)$ and let level
$k=(k_1,...,k_m)$. 
Since $\A_{G_k}$ verifies  equation (\ref{ra}) by \cite{Xu6}, we have
the following:
\begin{corollary}\label{allaa}
Let $\Gamma\subset G=SU(N_1)\times SU(N_2)\times\cdots\times SU(N_m) $ be 
a finite subgroup. Then 
every irreducible DHR representation of $\A_{G_k}^\Gamma$ is contained in
the restriction of $\pi_{\Lambda,g}$ to  $\A_{G_k}^\Gamma$ for some 
$\Lambda,g\in \Gamma$ where $\pi_{\Lambda,g}$ is defined as in definition 
(\ref{affinesll}) and  $\A_{G_k}^\Gamma$ is the conformal
net associated with the projective representation of $LG$ at level 
$k=(k_1,...,k_m) $.
\end{corollary} 
\subsection{Identifying representations of $\A^\Gamma$ for 
non-exceptional $\Gamma$}
In this section  we assume that $\Gamma$ is a non-exceptional finite subgroup 
of $G$ (cf. \ref{non-excep}).
Assume that $g\in \Gamma$ is a non-exceptional element in $\Gamma$ with
$g=e^{2\pi i\beta}$ and $G_g =G_\beta.$ 
We will choose the path $P_g$ as 
$P_g(\theta)=e^{2\pi i \theta \beta}, 0\leq \theta\leq 1.$
Let $\sigma$ be an irreducible character of the group 
$\Gamma_\beta:=\Gamma\cap G_\beta=\Gamma_g$. Let 
\begin{equation}\label{projector}
P_{\Lambda, \sigma}:=\frac{\sigma(1)}{|\Gamma_g|}\sum_{h\in \Gamma_\beta}
\sigma^*(h)\pi_{\Lambda}(h)
\end{equation} 
By Lemma \ref{commu} 
$P_{\Lambda, \sigma}\pi_{\Lambda,g}$ is a direct sum of 
$\sigma(1)$ copies of a DHR representation of $\A^{\Gamma}$ (on 
$P_{\Lambda, \sigma}\H_{\Lambda}$) which we denote by 
${\pi}_{\Lambda,g, \sigma}$. We have:
\begin{proposition}\label{id1}
Let $h\in N_G(\Gamma_g):=\{ b\in G| b\Gamma_g b^{-1}= \Gamma_g \}$. 
Then as representation of $\A^{\Gamma_g}$ we have
$$
\pi_{\Lambda,g, \sigma}\cdot\Ad_{h^{-1}}\simeq \pi_{\Lambda, hgh^{-1}, \sigma^h}
$$
where $\sigma^h$ is an irreducible representation of $\Gamma_{hgh^{-1}}$
defined by $\sigma^h(b)= \sigma (h^{-1}bh)$.
\end{proposition}
\proof
By definition (\ref{affinesl}) $\forall x\in \A(I), I\subset \mathbb R$ 
we have
\begin{align}
\hat{\pi}_{\Lambda,g}( \Ad_{h^{-1}} x) = \pi_{\Lambda}( P_{g,I} h^{-1} x 
h P_{b,I}^*)
= \pi_{\Lambda}(g)^* \pi_{\Lambda} (h  P_{g,I} h^{-1} x h P_{g,I}^* h^{-1}) 
\pi_{\Lambda}(h)\\
= \pi_{\Lambda}(h)^* \hat{\pi}_{\Lambda, hgh^{-1}} (x) \pi_{\Lambda}(h)   
\end{align}
On the other hand from the definition (\ref{projector}) one checks that
$$
\pi_{\Lambda}(h)^* P_{\Lambda, \sigma^h} \pi_{\Lambda}(h)
=  P_{\Lambda, \sigma}
$$
It follows that $\forall y\in \A^{\Gamma_g}(I)$ 
$$
\pi_{\Lambda,g,\sigma}\cdot\Ad_{h^{-1}} y =
\pi_{\Lambda}(g)^* \pi_{\Lambda, hgh^{-1},\sigma^h } (y) \pi_{\Lambda}(g)
.$$ 
\endproof
\begin{proposition}\label{id2}
For the pair of non-exceptional triples
$X=(\Lambda, g, \sigma)$ 
and 
$$
\zeta(X):= (\sum_\nu( w_{j_\nu}\Lambda^{\nu} + k_\nu\Lambda_{j_\nu},
\zeta g, \sigma\otimes(\otimes_\nu\sigma_{j_\nu}))
$$ 
where $\sigma_{j_\nu}$ is defined as in (\ref{phase}), 
we have
$\pi_X \simeq \pi_{\zeta(X)}$ as DHR representations of $\A^{\Gamma_g}_0$. 
\end{proposition}
\proof
For any $a\in \A(I)$ we have:
$$
\hat{\pi}_{\zeta(\Lambda), \zeta g}(a) = 
\pi_{\zeta(\Lambda)} (P_{\zeta g}P_g^* P_g a P_g^* 
P_g  P_{\zeta g}^*)
= \pi_{\zeta(\Lambda)} (P_\zeta P_g x P_g^*P_\zeta^* ) 
$$
where  we have used (2) of Lemma \ref{nexcp}. By (3) of Lemma \ref{cocycle}
There exists a unitary $U$ such that 
$$
\pi_{\zeta(\Lambda)} (P_\zeta P_g a P_g^*P_\zeta^* ) 
=  U \pi_{\Lambda,g} ( a  ) U^*.
$$ 
By (4) of Lemma \ref{cocycle}
$$
\pi_{\zeta(\Lambda)} (h) = U  \pi_{\Lambda} (h) \sigma_\zeta(h) U^*
,$$ and it follows by definition (\ref{projector}) 
$$
P_{\zeta(\Lambda),\sigma\otimes \sigma_\zeta} = U 
P_{\Lambda,\sigma}U^*
,$$
hence the proposition is proved by definition.
\endproof
\subsection{Details on decomposing solitons: fixed point resolutions}  
Assume that $g\in \Gamma$ is a non-exceptional element with
$g=e^{2\pi i\beta}$ and $G_g=G_\beta$. We will choose the path $P_g$ as 
$P_g(\theta)=e^{2\pi i \theta \beta}, 0\leq \theta\leq 1.$
Let $\hat{\pi}_{\Lambda,g}\res \A^\Gamma \simeq \sum_{i} m_i\beta_i$ 
where $\beta_i$
are irreducible DHR representations of $\A^\Gamma$. Define
$\Gamma_g:=\{h\in \Gamma| hg=gh \}$. Note that $\Gamma_g$ is a normal
subgroup of $\Gamma_{\Lambda,g}$ and $\Gamma_{\Lambda,g}/\Gamma_g
=\{ h\in Z(G)| h\Lambda=\Lambda \}$ is an abelian group (cf. (5) of
Lemma \ref{afusion}). 
\begin{lemma}\label{commu}
For all $x\in \A (I), h\in \Gamma_g$,
$$
\pi_\Lambda(h) \hat{\pi}_{\Lambda,g}(x) \pi_\Lambda(h)^* =
\hat{\pi}_{\Lambda,g}(hxh^*).
$$
\end{lemma}
\proof
Since $\pi_1(\tilde L_IG)$ generates $\A(I)$, it is sufficient to check the
equation for $x=\pi_1(y), y\in L_IG$. As elements in
$LG$ we have
$$
h P_{g,I} y  P_{g,I}^{-1} h^{-1} =
P_{g,I}  h y  h^{-1}P_{g,I}^{-1}
$$
where we have used $h P_{g}h^{-1}= P_{g}$ by (1) of Lemma \ref{nexcp}. 
It follows by definition (\ref{affinesl}) that
$$
\pi_\Lambda(h) \hat{\pi}_{\Lambda,g}(x) \pi_\Lambda(h)^* =
 \hat{\pi}_{\Lambda,g}(hxh^*).
$$
\endproof
Assume that when restricting to $\A^{\Gamma_g}$, 
$\H_\Lambda=\bigoplus_{\sigma\in E } M_\sigma \otimes 
V_\sigma$ where
$V_\sigma$ are irreducible representation spaces of $\Gamma_g$, $E\subset
{\rm Irr}\Gamma_g$  and 
$M_\sigma$ the
corresponding multiplicity spaces. By Th. 4.1 of \cite{KT}, 
$\sigma$ appears in the above decomposition iff 
$\sigma|Z(G)= \Lambda|Z(G)$.   
Applying Th. \ref{projrep} to the pair $\A^{\Gamma_g}\subset \A$, each
$ M_\sigma$ with $\sigma|Z(G)= \Lambda|Z(G)$ is an irreducible  
DHR representation of $\A^{\Gamma_g}$. We will denote 
$ M_\sigma$ by $\pi_{\Lambda,g,\sigma}$. 
When $\Gamma_{\Lambda,g}/\Gamma_g$ is nontrivial, 
the next question is how $\pi_{\Lambda,g,\sigma} $ 
decomposes when restricting to
$\A^{\Gamma_{\Lambda,g}}$. This is the issue of ``fixed point resolutions'',
since the action of the center has a nontrivial fixed point on the 
quadruples as described on Page 78 of \cite{KT}, and 
the question about the nature of 
how $\pi_{\Lambda,g,\sigma} $ 
decomposes as representation of $\A^\Gamma$ is implicitly raised. 
Assume that $\Gamma_{\Lambda,g}/\Gamma_g
=\{ h\in Z(G)| h\Lambda=\Lambda \}$. Then 
$\A^{\Gamma_{\Lambda,g}}\subset \A^{\Gamma_g}$ is the fixed point
subnet under the action of $\Gamma_{\Lambda,g}/\Gamma_g $. 
Note that $\Gamma_{\Lambda,g}/\Gamma_g
\simeq \{ \zeta\in Z(G)| \zeta\Lambda=\Lambda \} $
and denote the isomorphism by $h\rightarrow \zeta(h)$.
Then we have:
\begin{theorem}\label{fix}
(1):$$
\langle \pi_{\Lambda,g,\sigma}\res \A^\Gamma, 
\pi_{\Lambda,g,\sigma}\res \A^\Gamma \rangle
=|\{ h\in \Gamma_{\Lambda,g}/\Gamma_g|
\sigma_{\zeta(h)}\simeq \sigma\otimes \sigma_\zeta \}|
$$
where $\sigma_{\zeta(h)}$ is as defined in (4) of Lemma \ref{cocycle};
\par
(2): $\pi_{\Lambda,g,\sigma}\res \A^\Gamma$ decomposes into
irreducible representations of $\A^\Gamma$ which are in one-to-one
correspondence with all irreducible projective representation of the 
group $\Gamma_{\Lambda,g}/\Gamma_g$ with a fixed cocycle.
\end{theorem}
\proof
Ad (1): 
$\A^{\Gamma_{\Lambda,g}}\subset \A^{\Gamma_g}$ is the fixed point
subnet under the action of $\Gamma_{\Lambda,g}/\Gamma_g $. 
Applying Lemma \ref{gh} to the pair 
$\A^{\Gamma_{\Lambda,g}}\subset \A^{\Gamma_g}$, 
$\langle \pi_{\Lambda,g,\sigma}\res \A^{\Gamma_g}, 
\pi_{\Lambda,g,\sigma}\res \A^{\Gamma_g} \rangle$ is equal to the number of
elements $h \in\Gamma_{\Lambda,g}/\Gamma_g $ such that
$ \pi_{\Lambda, g,\sigma} \simeq  \pi_{\Lambda, g,\sigma}(\Ad.h)$
as representations of $\A^{\Gamma_g}$. 
By Prop.\ref{id1}  $\pi_{\Lambda, g,\sigma}(\Ad_h) \simeq 
\pi_{\Lambda, hgh^{-1},\sigma^h} = \pi_{\Lambda, \zeta(h)g,\sigma^h}
,$ 
and by Prop.\label{s2} $\pi_{\Lambda, g,\sigma}\simeq 
\pi_{\Lambda, \zeta(h)g,\sigma \otimes\sigma_{\zeta(h)}}$ as 
representations of  $\A^{\Gamma_g}$. It follows that
$ \pi_{\Lambda, g,\sigma} \simeq  \pi_{\Lambda, g,\sigma}(\Ad.h)$
as representations of  $\A^{\Gamma_g}$
iff  $\sigma^h\simeq \sigma\otimes \sigma_{\zeta(h)}.$ 
Hence 
$$
\langle \pi_{\Lambda,g,\sigma}\res \A^{\Gamma_{\Lambda,g}}, 
\pi_{\Lambda,g,\sigma}\res \A^{\Gamma_{\Lambda,g}} \rangle
=|\{ h\in \Gamma_{\Lambda,g}/\Gamma_g|
\sigma_{\zeta(h)}\simeq \sigma\otimes \sigma_\zeta \}|
$$
By
(4) of Lemma \ref{gh}  (1)  is proved. (2) follows by applying 
Th. \ref{projrep} to the pair $\A^{\Gamma_{\Lambda,g}}\subset \A^{\Gamma_g}$
and (4) of Lemma \ref{gh}. 
\endproof
Combine the above theorem  with   Cor. \ref{cyclic1} we immediately have:
\begin{corollary}\label{cyclic2}
If the group $\{ h\in \Gamma_{\Lambda,g}/\Gamma_g|
\sigma^h\simeq \sigma\otimes \sigma_\zeta(h)\}$ is cyclic of order $m$, then
$ \pi_{\Lambda,g,\sigma}\res \A^\Gamma$ decomposes into $m$ irreducible
pieces.
\end{corollary}
\subsection{An example}
Here we illustrate Cor. \ref{cyclic2} in the example 6.4
of \cite{KT}. We keep the same notation of \cite{KT}. Set $G=SU(2)$ and
$\Gamma= H_8$ the quaternion group. $H_8$ has 8 elements,$\{1,\epsilon,q_i,
\epsilon,q_i,i=1,2,3, \};$ they obey the multiplication rules
$q_i^2 =\epsilon, q_iq_jq_i^{-1}= \epsilon q_j= q_j^{-1}, i\neq j$.
We note that $q_i, \epsilon q_i$ are non-exceptional elements of
$SU(2)$. 
The centralizer of $\Gamma_{q_i}\simeq {\mathbb Z}_4$, and we will label
its irreducible representations by the exponents $\sigma=0,+1,-1,2$. 
There are 5 irreducible representations of $H_8$, $\{ \alpha_0, \alpha_1,
\alpha_2,\alpha_3,\alpha_4 \}$ with dimensions $1,1,2,1,1$ respectively. The
characters of these representations are given on Page 94 of \cite{KT}.\par
Consider the net $\A_{SU(2)_{2k_1}}$. The irreducible DHR representations
of $\A_{SU(2)_{2k_1}}$ are labeled by irreducible representations of 
$\tilde LSU(2)$ at level $2k_1$, and we will use integers $0,1,...,2k_1$
to label these representations such that $0$ is the vacuum representation.
The only representation which is fixed by the action of the center is
$k_1$. We note that $\sigma_{\epsilon} = 2k_1({\rm mod}) 4$.
When $k_1$ is odd, consider the DHR representation $\pi_{k_1,q_j,1}$.
We have $\Gamma_{k_1,q_j} = H_8$.   
We note that $\sigma_{\epsilon} = 2k_1({\rm mod}) 4$, 
and so the stabilizer of 
$\pi_{k_1,q_j,\pm 1}$ is $\{ h\in  H_8/{\Bbb Z_4}| \sigma^h \simeq 
\sigma_{\zeta(h)} \} \simeq {\Bbb Z}_2$. Hence by Cor. \ref{cyclic2}, 
 $\pi_{k_1,q_j,\pm 1}$ decomposes into two distinct 
irreducible DHR representations 
of $\A_{SU(2)_{2k_1}}^{H_8}$. When $k_1=1$ this is first observed in 
\cite{KT} by identifying  $\A_{SU(2)_{2k_1}}^{H_8}$ with the tensor 
products of three ``Ising Models'' (cf. Page 99 of \cite{KT}).\par
When $k_1$ is even, consider the DHR representation $\pi_{k_1,q_j,0}$ or
 $\pi_{k_1,q_j,2}$. Similar as above the stabilizer of 
 $\pi_{k_1,q_j,0}$ or
 $\pi_{k_1,q_j,2}$ is ${\Bbb Z}_2$, and by using  Cor. \ref{cyclic2}
again we conclude that  $\pi_{k_1,q_j,0}$ or
$\pi_{k_1,q_j,2}$ decomposes into two distinct  
irreducible DHR representations 
of $\A_{SU(2)_{2k_1}}^{H_8}$.
\section{Constructions of solitons for permutation orbifolds}
\subsection{Preliminaries on cyclic orbifolds}
In the rest of this paper we assume that $\A$ is completely rational. 
$\D:= \A\otimes\A...\otimes\A$ ($n$-fold tensor product)
and $\B:=\D^{\mathbb Z_n}$ 
(resp. $\D^{\mathbb P_n}$
where $\mathbb P_n$ is the permutation group on $n$ letters) 
is the fixed point subnet 
of $\D$ under the action of cyclic permutations (resp. permutations).
Recall that $J_0=(0,\infty)\subset \mathbb R$. Note that the action of
$\mathbb Z_n$ (resp. $\mathbb P_n$) on $\D$ is faithful and proper. 
Let $v\in \D(J_0)$ be a unitary such that $\beta_g(v)=e^{\frac{2\pi i}{n}}
v$ (such $v$ exists by P. 48 of \cite{ILP}) where $g$ is 
the generator of the cyclic group $\mathbb Z_n$ and $\beta_g$ stands for
the action of $g$ on $\D$. Note that $\sigma:=\Ad_v$ is a DHR representation
of $\B$ localized on $J_0$.
Let $\gamma: \D(J_0)\rightarrow \B(J_0)$ be the 
canonical endomorphism from $ \D(J_0)$ to $\B(J_0)$ and let $\gamma_{\B}:=
\gamma\res \B(J_0)$. Note $[\gamma]=[1]+[g]+...+[g^{n-1}]$ as sectors of
$\D(J_0)$ and $[\gamma_{\B}]=[1]+[\sigma]+...+[\sigma^{n-1}]$ as sectors of
$\B(J_0)$. Here $[g^i]$ denotes the sector of $\D(J_0)$ which is the 
automorphism 
induced by $g^i.$ 
All the sectors considered
in the rest of  this paper will be sectors of  $ \D(J_0)$ or  $ \B(J_0)$ as
should be clear from their definitions. 
All DHR representations will be assumed to be localized
on $J_0$ and have finite statistical dimensions unless noted otherwise. 
For simplicity of notations,
for a DHR representation
$\sigma_0$ of $\D$ or $\B$ localized on $J_0$,
we will use the same notation $\sigma_0$ to denote its restriction to
$ \D(J_0)$ or  $ \B(J_0)$
and we will make no distinction
between local and global intertwiners  
for DHR representations localized on $J_0$ since they are the same by the
strong additivity of $\D$ and $\B$. The following is Lemma 8.3 of \cite{LX}:
\begin{lemma}\label{grading}
Let $\mu$ be an irreducible DHR representation of $\B$. Let
$i$ be any integer. Then:\par 
(1) $G(\mu,\sigma^i):=\e (\mu,\sigma^i) \e (\sigma^i,\mu) \in {\mathbb C},
$ 
$G(\mu,\sigma)^i=G(\mu,\sigma^i) $. Moreover $G(\mu,\sigma)^n=1$;\par
(2) If $\mu_1\prec \mu_2\mu_3$ with
$\mu_1,\mu_2,\mu_3$ irreducible, then $G(\mu_1,\sigma^i)=G(\mu_2, \sigma^i)
G(\mu_3,\sigma^i)$;\par
(3)   $\mu$ is untwisted if and only if  $G(\mu,\sigma)=1;$ \par
(4) $G(\bar \mu, \sigma^i)=\bar G(\mu, \sigma^i).$
\end{lemma}

\subsection{ One cycle case}
First we recall the construction of solitons  for permutation orbifolds
in \S6 of \cite{LX}. 
Let $h:S^1\setminus\{-1\}\simeq\mathbb R\to S^1$ be a smooth, 
orientation preserving, injective map 
which is smooth also at $\pm\infty$, namely the left and right limits
$\lim_{z\to -1^{\pm}}\frac{{\rm d}^n h}{{\rm d}z^n}$ exist for all 
$n$.

The range $h(S^1\setminus\{-1\})$ is either $S^1$ minus a 
point or a (proper) interval of $S^1$.

With $I\in\I$, $-1\notin I$, we set
\[
\Phi_{h,I}\equiv \Ad U(k)\ ,
\]
where $k\in\Diff(S^1)$ and $k(z)=h(z)$ for all $z\in I$ and $U$ is 
the projective unitary representation of $\Diff(S^1)$ associated 
with $\A$.
Then $\Phi_{h,I}$ does not depend on the choice of 
$k\in\Diff(S^1)$ and 
\[
\Phi_{h}:I\mapsto \Phi_{h,I}
\]
is a well defined soliton of $\A_0\equiv\A\restriction\mathbb R$.

Clearly $\Phi_h(\A_0(\mathbb R))''=\A(h(S^1\setminus\{-1\}))''$, 
thus $\Phi_h$ is irreducible if the range of $h$ is  
dense, otherwise it is a type III factor representation.  It is easy 
to see that, in the last case, $\Phi_h$ does not depend on $h$ 
up to unitary equivalence.

Let now $f:S^1\to S^1$ be the degree $n$ map $f(z)\equiv z^n$. 
There are $n$ right inverses 
$h_i$, $i=0,1,\dots n-1$, for $f$ ($n$-roots); namely there are $n$ injective 
smooth maps $h_i:S^1\setminus\{-1\}\to S^1$ such that $f(h_i(z))=z$, 
$z\in S^1\setminus\{-1\}$. The $h_i$'s are smooth also at $\pm\infty$.

Note that the ranges $h_i(S^1\setminus\{-1\})$ are $n$ pairwise 
disjoint intervals of $S^1$, thus we may fix the labels of 
the $h_i$'s so that 
these intervals are counterclockwise ordered, namely
we have $h_0(1)<h_1(1)<\dots<h_{n-1}(1)<h_0(1)$, and
we choose $h_j= e^{\frac{2\pi ij}{n}}h_0, 0\leq j\leq n-1.$ 
\par
For any interval $I$ of $\mathbb R$, we set
\begin{equation}\label{ts}
\pi_{1,\{0,1...n-1 \},I}\equiv \chi_I\cdot 
(\Phi_{h_{0},I}\otimes\Phi_{h_1,I}\otimes\cdots\otimes\Phi_{h_{n-1},I})\ ,
\end{equation}
where $\chi_I$ is 
the natural isomorphism from $\A(I_0)\otimes\cdots\otimes\A(I_{n-1})$ to 
$\A(I_0)\vee\cdots\vee\A(I_{n-1})$ given by the split property, 
with $I_k\equiv h_k(I)$. 
Clearly $\pi_{1,\{0,1...n-1 \}} $ is a soliton of 
$\D_0\equiv\A_0\otimes\A_0\otimes\cdots\otimes\A_0$ ($n$-fold 
tensor product).
Let $p\in \mathbb P_n$. We set
\begin{equation}\label{1cycle1}
\pi_{1,\{p(0),p(1),...,p(n-1) \}}=\pi_{1,\{0,1...,n-1\}} \cdot\b_{p^{-1}}
\end{equation}
where $\b$ is the natural action of $\mathbb P_n$ on $\D$, and 
$\pi_{1,\{0,1...,n-1\}} $ is as in (\ref{ts}). 
The following is part of Prop. 6.1 in \cite{LX}:
\begin{proposition}\label{pind1}
\label{pi}
$(1)$: ${\rm Index}(\pi_{1,\{0,1...,n-1\}} ) =\mu_{\A}^{n-1}$.

$(2)$: The conjugate of $\pi_{1,\{0,1...,n-1\}} $ is  
 $\pi_{1,\{0,n-1,n-2,...,1\}} $. 
\end{proposition}
Let $\lambda$ be a DHR representation of $\A$. Given an 
interval $I\subset S^1\setminus\{-1\}$,  we set
\begin{definition}\label{1cycle}
\[
\pi_{\lambda,\{ p(0),p(1),...,p(n-1)\}, I}(x)= 
\pi_{\lambda,J}(\pi_{1, \{ p(0),p(1),...,p(n-1)\}, I}
(x))\ 
,\quad x\in\D(I)\ ,
\]
where $ \pi_{1, \{ p(0),p(1),...,p(n-1)\}, I} $ is defined as in 
(\ref{1cycle1}), 
and $J$ is any interval which contains $I_0\cup I_1\cup...\cup I_{n-1}$. 
Denote the corresponding soliton by $\pi_{\lambda,\{ p(0),p(1),...,p(n-1)\}}.$
When $p$ is the identity element in $\mathbb P_n$, we will denote the
corresponding soliton by $\pi_{\lambda, n}$. 
\end{definition}
The following follows from  
Prop. 6.4 of \cite{LX}:
\begin{proposition}\label{non-v}
The above definition is independent of the choice of $J$, thus 
$\pi_{\lambda, \{p(0),p(1),...p(n-1)\},I}$ is a well defined soliton of $\D$. 

We can localize $\pi_{1,\{p(0),p(1),...p(n-1)\}}$, 
$\pi_{\lambda, \{p(0),p(1),...p(n-1)\}} $ and $\lambda$  on
$J_0$. Denote by $\tilde \pi, \tilde \pi_{\lambda}$ and 
$(\l,1,1,...,1):=\l\otimes \iota \otimes \iota\cdots \otimes 
\iota\restriction\D(J_0)$ respectively the 
corresponding endomorphisms of  $\D(I)$. 
Then  as sectors of $\D(J_0)$ we have
\[
[\tilde\pi_{\lambda}]= 
[\tilde\pi\cdot(\l,1,1,...1)\ ].
\]
In particular 
${\rm Index}(\pi_{\lambda, \{p(0),p(1),...p(n-1)\}})= d(\lambda)^2 
\mu_{\A}^{n-1}.$ 
\end{proposition}
\subsection{General case}
Let $\psi: \{ 0,1,...,n-1 \}\rightarrow \L$ where $\L$ is the set of
all irreducible DHR representations of $\D$. For any $p\in \mathbb P_n$ 
we set $p.\psi (i):= \psi( p^{-1}.i), i=0,...,n-1$ where $\mathbb P_n$ acts 
via permutation on the $n$ numbers $ \{ 0,1,...,n-1 \}$. Assume that 
$p.\psi=\psi$, and $p=c_1...c_k$ is a product of disjoint cycles. 
Since
$p.\psi=\psi$, $\psi$ takes the same value denoted by $\psi(c_j)$ 
on the elements $\{a_1,a_2,... a_l \}$ of each cycle $c_j=(a_1...a_l)$.
A {\it presentation} $f_j$ of the cycle $c_j=(a_1...a_l)$ is a 
list of numbers $\{b_1,...,b_l\}$ such that $(b_1...b_l)=c_j$ as cycles.
The length $l(f_j)$ of $f_j$ is $l$. We note that for a cycle of 
length $l$ there are $l$ different presentations.  
For each element $x=x_0\otimes x_1\otimes \cdots \otimes x_{n-1} \in \D$,
and each cycle $c=(a_1...a_l)$ with a fixed presentation $f=
\{b_1,...b_l\}$, we define 
$x_{c,f} = x_{b_1}\otimes x_{b_2}\otimes \cdots \otimes x_{b_l}$. Now 
we are ready to define solitons for permutation orbifolds:
\begin{definition}\label{sp}
Assume that $p.\psi=\psi$ and 
$p=c_1...c_k$ is a product of disjoint cycles 
as above. For each $c_j$ we fix a presentation $f_j$. 
Then for any 
$x= x_0\otimes x_1\otimes \cdots \otimes x_{n-1} \in \D(I), 
I\subset S^1\setminus {-1}=\mathbb R$, 
$$\pi_{\psi,p}\equiv\pi_{\psi, c_1c_2...,c_k,f_1,...f_k} (x)= \pi_{\lambda_1,l(f_1)}
(x_{c_1,f_1})
\otimes  \pi_{\lambda_2,l(f_2)}(x_{c_2,f_2})
\otimes \cdots \otimes \pi_{\lambda_k,l(f_k)}(x_{c_k,f_k})
$$
on $\H_{\psi(c_1)}\otimes\H_{\psi(c_2)}\otimes\cdots\otimes \H_{\psi(c_k)}$
where $\pi_{\lambda_j,l(f_j)}$ is as in Def. \ref{1cycle}.
\end{definition}
Here and in the following, to simplify notations, 
we do not put the interval suffix $I$ in a 
representation, if no confusion arises.
\begin{lemma}\label{equv}
The unitary equivalence class of $\pi_{\psi, p}$ in Definition \ref{sp}
depends only on $ p \in \mathbb P_n$.
\end{lemma}
\proof
We have to check that the unitary equivalence class of $\pi_{\psi, p}$ 
in Definition \ref{sp} is independent of the order $c_1,...c_k$ and the 
presentation of $c_j$. The first case is obvious, and second case follows
from (a) of Prop. 6.2 in \cite{LX}.
\endproof
Due to the above lemma, for each $p\in \mathbb P_n$ we will fix a choice
of the order $c_1,...c_k$ and presentations of $ c_1,...c_k$. For
simplicity we will denote the corresponding soliton simply by
$\pi_{\psi, p}$.\par 
\begin{proposition}\label{conj}
$\pi_{h.\psi, hph^{-1}} \simeq \pi_{\psi, p}\cdot\beta_{h^{-1}}$
as solitons of $\D_0$, $p,h\in\mathbb P_n$.
\end{proposition}
\proof
Let $p=c_1...c_k$ be a product of disjoint cycles with $c_j=(a_1...a_l)$.  
Then
$ hph^{-1}= hc_1h^{-1}... hc_kh^{-1}$ with $ hc_jh^{-1}=(h(a_1)...h(a_l))$.
Note that $h.\psi(h(a_1))=\psi(a_1)= \psi(c_j)$, and
$\beta_{h^{-1}}(x_0\otimes x_1\otimes\cdots\otimes x_{n-1})
= x_{h(0)}\otimes x_{h(1)}\otimes\cdots\otimes x_{h(n-1)}$, 
$\forall x_0\otimes x_1\otimes\cdots\otimes x_{n-1}\in \D(I)
.$
The proposition now follows directly from definition (\ref{sp}).
\endproof
\section{Identifying solitons in the permutation orbifolds}
The goal in this section is to prove the following:
\begin{theorem}\label{permu}
Let $\pi_{\psi_1, p_1}$,  $\pi_{\psi_2, p_2}$ be two solitons
as given in definition (\ref{sp}). Then
$\pi_{\psi_1, p_1} \simeq \pi_{\psi_2, p_2}$ as 
solitons of $\D_0$ if and only 
if $\psi_1=\psi_2, p_1=p_2$.
\end{theorem}
We note that even for the first nontrivial case $n=3$ we do not know
a direct proof of the theorem. Our proof is indirect and is divided into
the following steps:
\subsection{Identifying solitons: Cyclic case}
We will first prove Th. \ref{permu} for the case when
both $p_1,p_2$ are one cycle. In this case $\psi_1$ (resp.  $\psi_2$)
is a constant function with value denoted by $\lambda_1$ (resp. $\lambda_2$). 
We will denote $\psi_1$ (resp.  $\psi_2$)  simply by $\lambda_1$ 
(resp. $\lambda_2$). If $g\in\Gamma$, we will denote by $\D^{\langle 
g\rangle}$ the fixed-point subnet of $\D$ under the subgroup generated 
by $g$.
\begin{proposition}\label{cycle1}
(1). Let $g_1=(01,...n-1)$ and $g_2=g_1^m$ with $(m,n)=1$. Then
$\pi_{\lambda_1,g_1}\simeq  \pi_{\lambda_1,g_2}$ if and only if
$\lambda_1= \lambda_2, g_1=g_2$;\par
(2). If $\pi_{\lambda_1,g_1}$ restricts to a DHR representation 
a subnet $\B$ with $\D^{\langle g_1\rangle}\subset \B \subset\D$, then $\B=
\D^{\langle g_1\rangle}$.
\end{proposition}
\proof
Ad (1):
It is sufficient to show that if $\pi_{\lambda_1,g_1}\simeq  \pi_{\lambda_1,g_2}$, then $\lambda_1= \lambda_2, g_1=g_2$. \par
Since $(m,n)=1$, there exists $h\in \mathbb P_n$ such that 
$h g_1h^{-1}= g_2$. By Prop. \ref{conj}, we can assume that
$ \pi_{\lambda_1,g_1}\simeq  \pi_{\lambda_2,g_1}\cdot\Ad_h$. 
As in \S8.3 of \cite{LX}, we denote the $n$ irreducible DHR
representations of $\D^{\langle g_1\rangle}$ of 
$\pi_{\lambda_1,g_1}$ by $\tau_{\lambda_1}^{(0)},...,
\tau_{\lambda_1}^{(n-1)}$. Since 
$ \pi_{\lambda_1,g_1}\simeq  \pi_{\lambda_2,g_1}\cdot\Ad_h$, we must have
$\tau_{\lambda_1}^{(0)}\simeq \tau_{\lambda_2}^{(i)}\cdot\Ad_h$ for
some $0\leq i\leq n-1$. By (48) of \cite{LX} we have that 
$[\tau_{\lambda_1}^{(0)}]\prec [(\lambda,1,...,1)\res \D^{\langle g_1\rangle} 
\tau^{(0)}]$, and (2) and (3) of 
Lemma \ref{grading} we have 
$G(\tau_{\lambda_1}^{(0)}, \sigma^{k(1)})=G(\tau^{(0)},  \sigma^{k(1)})  
= e^{\frac{2\pi i}{n}}$, where $1\leq k(1)\leq n-1$ and $(k(1),n)=1$ (cf.
the paragraph after (47)). 
Similarly $G(\tau_{\lambda_2}^{(i)}, \sigma^{k(1)})=e^{\frac{2\pi i}{n}}.$
On the other hand note that by 
definition 
$$
G(\tau_{\lambda_2}^{(i)}\cdot\Ad_h),   \sigma^{k(1)}\cdot\Ad_h)=
G(\tau_{\lambda_2}^{(i)}, \sigma^{k(1)})   
= e^{\frac{2\pi i}{n}}
$$
Since $[\Ad h.g] =[g^m]$, we have  $\sigma\cdot\Ad_h\simeq\sigma^{m}$, 
and so we have 
$$G(\tau_{\lambda_2}^{(i)}\cdot\Ad_h,   \sigma^{k(1)}\cdot\Ad_h)
=G(\tau_{\lambda_2}^{(i)}\cdot\Ad_h,   \sigma^{mk(1)})   
=G(\tau_{\lambda_1}^{(0)},   \sigma^{k(1)}))^m
= e^{\frac{2\pi i}{n}}
$$
where in the second = we have used (1) of Lemma \ref{grading}.
Hence 
$ e^{\frac{2\pi i}{n}}= e^{\frac{2\pi im}{n}}$
and it follows that $m=1$ since $(m,n)=1$. So we have 
$\pi_{\lambda_1, g_1}\simeq \pi_{\lambda_2, g_1}$, and by (2) of 
Th. 8.6 in \cite{LX} we have $\lambda_1=\lambda_2$.
\par
Ad (2): First we note that the subnet $\B=\D^{\langle g_1^l\rangle}$ 
for some
$1\leq l\leq n, n=ll_1$ by the Galois correspondence (cf. \cite{ILP}). 
Also the vacuum representation of $\D^{\langle g_1^l\rangle}$ restricts to 
$\bigoplus_{ 1\leq i \leq l} \sigma^{l_1i}$ of 
$\D^{\langle g_1 \rangle}$
. If 
$\pi_{\lambda_1, g_1}$ restricts to a DHR representation of
$\D^{\langle g_1^l\rangle}$ , by applying (3) of Lemma \ref{grading}
to the pair $\D^{\langle g_1 \rangle}\subset  \D^{\langle g_1^l\rangle}$
we conclude that 
$G(\tau_{\lambda_1}^{(0)}, \sigma^{l_1})= 1$. Since
$G(\tau_{\lambda_1}^{(0)}, \sigma^{k(1)})=e^{\frac{2\pi i}{n}}$, by using
(1) of   Lemma \ref{grading} we have
$$
G(\tau_{\lambda_1}^{(0)}, \sigma^{l_1k(1)})= 
G(\tau_{\lambda_1}^{(0)}, \sigma^{l_1})^{k(1)}=1
= G(\tau_{\lambda_1}^{(0)}, \sigma^{k(1)})^{l_1}
=e^{\frac{2\pi l_1i}{n}} 
$$
Hence $n|l_1$ and we conclude that $l_1=n$, $\B=\D^{\langle g_1\rangle}$.  
\endproof
\begin{proposition}\label{cycle2}
Let $g_1$ (resp. $g_2$) be one cycle of length $n$. Then
$\pi_{\lambda_1, g_1}\simeq \pi_{\lambda_2, g_2}$ if and only if
$\lambda_1= \lambda_2, g_1=g_2$.
\end{proposition}
\proof
It is  sufficient to show that if $\pi_{\lambda_1, g_1}\simeq \pi_{\lambda_2, g_2}$, then $\lambda_1= \lambda_2, g_1=g_2$.\par
Note that $\pi_{\lambda_1, g_1}$ (resp. $\pi_{\lambda_2, g_2}$)
restricts to a DHR representation of
$\D^{\langle g_1\rangle}$ (resp. $\D^{\langle g_2\rangle}$), it follows that
$\pi_{\lambda_1, g_1}$ restricts to a DHR representation of 
$\D^{\langle g_2\rangle}$. By Prop.\ref{normal} $\pi_{\lambda_1, g_1}$ 
restricts to a DHR representation of 
$\D^{\langle g_1\rangle}\vee \D^{\langle g_2\rangle}$, 
and by (2) of \ref{cycle1} we must have
$\D^{\langle g_1\rangle}\vee \D^{\langle g_2\rangle}= 
\D^{\langle g_1\rangle}
.$
It follows that $\D^{\langle g_2\rangle}\subset 
\D^{\langle g_1\rangle}$ and by Galois correspondence again (cf. \cite{ILP})
we have $\langle g_2\rangle\subset \langle g_1\rangle$. Exchanging
$ g_1$ and $g_2$ we conclude that
$\langle g_2\rangle= \langle g_1\rangle$. Hence 
$g_2=g_1^m$ for some integer $m$ with $(m,n)=1$. By (1) of 
Prop.\ref{cycle1} we have proved that
$g_1=g_2, \lambda_1=\lambda_2.$
\endproof
\subsection{Proof of Th.\ref{permu} for general case and its corollary}
Assume that $g_1=c_1c_2...c_k$ and $g_2=c_1'...c_l'$ 
where $c_j$ (resp. $c_i'$) are disjoint cycles.
Fix $1\leq j\leq k$ and 
let $c_j=(a_1...a_m)$. 
Let us first show that $a_1,...,a_m$ must appear in one cycle of $g_2$.
Let $U$ be the 
unitary such that $\pi_{\psi_1,g_1}=\Ad_U\cdot\pi_{\psi_2,g_2}$. 
Choose $x=x_0\otimes x_1\otimes\cdot\otimes x_{n}\in \D_0$ such that
$x_i=1$ if $i\neq a_j, j=1,...,m$, and no
other constraints.  Denote by $\D_{0,c_j}$ the subalgebra
of $\D_0$ generated by such elements.
We note that
$ \pi_{\psi_1,g_1}(\D_{0,{c_j}})$ is $B(\H_{c_j})$, a type I factor by 
strong additivity. If
 $a_1,...,a_m$  appear in more than one cycle of $g_2$, then by definition
(\ref{ts}) 
$ \pi_{\psi_2,p_2}(\D_{0,{c_j}})$ will be tensor products of factors 
of the form $\pi_\lambda( \A_{J})$, where $J$ is a union of intervals 
of $S^1$, but $\bar J\neq S^1$, and so 
$ \pi_{\psi_1,p_1}(\D_{0,{c_j}})$ will be tensor 
products of type III factors, contradicting $
\pi_{\psi_1,p_1}(\D_{0,{c_j}})= U  \pi_{\psi_2,p_2}(\D_{0,{c_j}}) U^*$. 
By exchanging the role of $g_1$ and $g_2$ we conclude that 
$a_1,...,a_m$ must be exactly the elements  in one cycle $c_i'$ of $g_2$ for
some $1\leq i\leq l$,  
and we have 
$
\pi_{\psi_1,p_1}(\D_{0,c_j})= U  \pi_{\psi_2,p_2}(\D_{0,c_i'}) U^*
.$
Let $\H=\H_{\psi_1(c_j)}\otimes \H_{r}$. 
We have $UB(\H_{\psi_1(c_j)})U^* =B(\H_{\psi_2(c_i')})$. 
Since every automorphism of a type I factor  is inner, there a unitary
$U_1 \H_{\psi_2(c_i')}\rightarrow  \H_{\psi_1(c_j)}$ such that 
$ UB(\H_{\psi_1(c_j)})U^*= U_1B(\H_{\psi_2(c_i') })U_1^*$. Hence 
$\pi_{\psi_1(c_j),l(c_j)}=U_1 \pi_{\psi_2(c_i'),l(c_i')}U_1^*$ on 
$\D_{0,{c_j}}$, and by Prop.\ref{cycle2} we conclude that 
$c_j=c_i', \psi_1(c_j) = \psi_2(c_i')$. Since $j$ is arbitrary, exchanging
the roles of $g_1$ and $g_2$ we have proved $g_1=g_2, \psi_1=\psi_2$.   
\begin{proposition}\label{indexp}
Assume that $p=c_1...c_k$ where $c_i$ are disjoint cycles. Let
$\psi$ be such that $p.\psi=\psi$. Then:\par
(1): The centralizer (cf. (\ref{gpi}))
of $\pi_{\psi,p}$ in $\mathbb P_n$ is  
$\Gamma_{\psi,}=\{ h\in \mathbb P_n| h.\psi=\psi,
hph^{-1}=p;$ \par
(2) If $p\in \mathbb Z_n$, the  centralizer (cf. (\ref{gpi}))
of $\pi_{\psi,p}$ in $\mathbb Z_n$ is  
$\Gamma_{\psi,p}=\{ h\in \mathbb Z_n| h.\psi=\psi,
hph^{-1}=p;$ \par
(3):
$d(\pi_{\psi,p})^2= \prod_{1\leq i\leq k}d(\psi(c_i))^2 \mu_{\A}^{n-k}.$
\end{proposition}
\proof
(1) (2) follows from Prop. \ref{conj} and Th. \ref{permu}. 
Assume that each cycle $c_i$ has length $m_i, 1\leq i\leq k.$ Then
$\sum_{1\leq i\leq k} m_i=n$. By definition (\ref{sp}), we have
$d(\pi_{\psi,p})^2= \prod_{1\leq i\leq k} d(\pi_{\psi(c_i)})^2$.
By Prop.\ref{non-v} and (1) of Prop. \ref{pind1} 
$ d(\pi_{\psi(c_i)})^2= d(\psi(c_i))^2 \mu_{\A}^{m_i-1},$
hence
$$
d(\pi_{\psi,p})^2= \prod_{1\leq i\leq k}d(\psi(c_i))^2 \mu_{\A}^{m_i-1}
=\prod_{1\leq i\leq k}d(\psi(c_i))^2 \mu_{\A}^{n-k}
.$$
\section{Identifying all the irreducible representations of the
permutation orbifolds}
\subsection{Cyclic orbifold case}
\begin{theorem}\label{allcyclic}
Let $g=(01...n-1)$. Then
every irreducible DHR representation of $\D^{\mathbb Z_n}$ appears as 
an irreducible summand of $\pi_{\psi,g^i}$ for some $\psi,g^i$.
\end{theorem}
\proof
By Prop. \ref{solitonres} 
$\pi_{\psi_1,g^{i_1}}\res \B\simeq \pi_{\psi_2,g^{i_2}}\res\B$ iff 
there exists $h\in \mathbb Z_n$ such that 
$\pi_{\psi_1,g^{i_1}}(\beta_{h^{-1}}) \simeq 
\pi_{\psi_2,g^{i_2}}
,$ and by Prop. \ref{conj} and (2) of Prop. \ref{indexp} we have 
$h.\psi_1=\psi_2,  hg^{i_1}h^{-1}= g^{i_2}$. Denote the orbit 
of $\pi_{\psi_1,g^{i_1}}$ under the action of $\mathbb Z_n$ 
by $\{\psi_1,g^{i_1} \}$. Note that the orbit   $\{\psi_1,g^{i_1} \}$
has length $\frac{n}{ |\Gamma_{\psi_1,g^{i_1}}|}$.  
By Th. \ref{solitonind}
the sum of index of the irreducible summands of $\pi_{\lambda,g^i}$
is $\frac{n^2}{|\Gamma_{\psi,g^{i}}| }d(\pi_{\lambda,g})^2$. 
Hence the sum of index of 
distinct irreducible summands of  $\pi_{\lambda,g}$ for all 
$\psi, g\in \mathbb Z_n$ is given by$
\sum_{\{\psi, g^i \}} \frac{n^2}{|\Gamma_{\psi,g^{i}}| }  
d(\pi_{\lambda,g^i})^2$ where the sum
is over different orbits.  Assume that $g^i= c_1...c_k$. Then 
$k=(n,i)$ (the greatest common divisor of $n$ and $i$) 
and each cycle $c_i$ has length $\frac{n}{(n,i)}$.  
For each element $\psi_2,g^{i_2}$
in the orbit $\{\psi, g^i \}$, by (3) of Prop. \ref{indexp} 
$d(\pi_{\psi_2,g^{i_2}})^2= d(\pi_{\lambda,g^i})^2= 
\prod_{1\leq j\leq (n,i)} d(\psi(c_j)^{2} \mu_\A^{n-(n,i)}   
$
Hence
\begin{multline}
\sum_{\{\psi, g^i \}} \frac{n^2}{|\Gamma_{\psi,g^{i}}| }  d(\pi_{\psi,g^i})^2
= \sum_{\lambda, 0\leq i \leq n }  \frac{1}{\frac{n}
{ |\Gamma_{\psi,g^{i}}|}} \frac{n^2}{|\Gamma_{\psi,g^{i}}| }
d(\pi_{\psi,g^i})^2\\
=n \sum_{\psi, 0\leq i \leq n }  \prod_{1\leq j\leq (n,i)} d(\psi(c_j)^{2}
 \mu_\A^{n-(n,i)}
= n^2 \mu_\A^n = \mu_{\D^{\mathbb Z_n}}
\end{multline}
where in the last $=$ we have used Th. \ref{orb}.
The theorem now follows from Th. 30 of  \cite{KLM}.
\endproof   
Let us now decompose $\pi_{\lambda,g}$ into irreducible pieces. 
In this case $\Gamma_{\lambda,g}= \mathbb Z_n$ since $g=(012...n-1)$
(cf. (2) of Prop. \ref{indexp}). By definition (\ref{ts})
$\forall x_0\otimes x_1\otimes\cdots\otimes x_{n-1} \in \D(I)$, 
\begin{align}
\pi_{\lambda,g}\cdot\Ad_{g^{-1}}(x_0\otimes x_1\otimes\cdots\otimes x_{n-1})
= \pi_{\lambda,g} (  x_1\otimes x_2\otimes\cdots\otimes x_{0})
\\
= \pi_{\lambda} (R(\frac{2\pi}{n})) 
\pi_{\lambda,g} (x_0\otimes x_1\otimes\cdots...\otimes x_{n-1})
\pi_{\lambda} (R(\frac{2\pi}{n}))^*
\end{align}
Here $\pi_{\l}(R(\cdot))$ denotes the unitary one-parameter rotation 
subgroup in the representation $\l$.
Note that $ \pi_{\lambda} (R(\frac{2\pi}{n}))^n =  \pi_{\lambda} 
(R(2\pi)) = C_\lambda {\rm id} $ for some complex number $C_\lambda, 
|C_\lambda|=1 $. 
Let $\Omega_\lambda\in \H_\lambda$ be a unit vector such that
$\pi_{\lambda} (R(\frac{2\pi}{n}) \Omega_\lambda =  C_\lambda'  
\Omega_\lambda$ with $(C_\lambda')^n =C_\lambda$. 
\begin{definition}\label{act1}
$\pi_{\lambda,g}(g):={C_\lambda'}^{-1} \pi_{\lambda} 
(R(\frac{2\pi}{n}))$, 
and $\pi_{\lambda,g}(g^i):= \pi_\lambda(g)^i$.
\end{definition} 
Then it follows that $g^i\rightarrow \pi_{\lambda,g}(g^i)$ 
gives a representation
of $\mathbb Z_n$ on $\H_\lambda$, and $ \pi_{\lambda,g}(g^i). \Omega_\lambda
=\Omega_\lambda$. So $\Omega_\lambda$ affords a trivial representation of 
of  $\mathbb Z_n$ on $\H_\lambda$.
It follows from Lemma \ref{fullspectrum1} 
that all irreducible representations of 
$\mathbb Z_n$ appear in the representation $\pi_{\lambda}$.
It follows by Th. \ref{projrep} that 
$\pi_{\lambda,g, i}, i\in \hat{\mathbb Z}_n$ are distinct irreducible
representations. \par
Note that 
$g^i= c_1...c_k$ is a product of $k=(n,i)$ disjoint cycles of the
same length $\frac{n}{(n,i)}.$ Let $h\in \Gamma_{\psi,g^i}$. Then 
$\Ad_h$ induces a permutation among the cycles $c_1,...,c_k$. 
We define an element $h'\in \mathbb P_k$ by the formula
$ hc_ih^{-1}=c_{h'(i)}, i=1,...,k$. We note that in the definition
of $\pi_{\psi,g}$ a presentation of $g$ has been fixed. Assume that  
$hf_{{h'}^{-1}(i)}h^{-1}= h''(i). f_i$ where $h''(i)$ is an element in the
cyclic group generated by $c_i$. 
Define
\begin{definition}\label{act2}
$$\pi_{\psi,g^i}(h):= h' \big(\pi_{\psi(c_1),c_1}(h''(1)) \otimes\cdots
\otimes  \pi_{\psi(c_k),c_k}(h''(k)\big)
$$
where the action of $h'\in {\mathbb P}_k$ on $\H_{\psi(c_1)} \otimes\cdots\H_{\psi(c_k)}$ 
is by permutation of the tensor factors, and 
$ \pi_{\psi(c_i),c_i}(h''(i))$ is as defined in definition (\ref{act1}).
\end{definition}
One checks easily that Definition \ref{act2} gives a representation of
$ \Gamma_{\psi,g^i},$
$
\Ad_{\pi_{\psi,g^i}(h)} \pi_{\psi,g^i}= \pi_{\psi,g^i}\Ad_h,$ and the vector 
$\Omega_{\psi(c_1)}\otimes\cdots \otimes \Omega_{\psi(c_k)}$
is fixed by $\pi_{\psi,g^i}(\Gamma_{\psi,g^i})$. It follows by 
lemma \ref{fullspectrum1} and Th. \ref{projrep} 
that we have proved the following: 
\begin{theorem}\label{decomposecy}
$ \pi_{\psi,g^i, \sigma\in \hat\Gamma_{\psi,g^i} }$ gives all the 
irreducible summands of  $\pi_{\psi,g^i}\res \D^{\mathbb Z_n}$.
\end{theorem}
We note that Th. \ref{allcyclic} and 
Th. \ref{decomposecy} generalizes the considerations of
\S8 of \cite{LX} for the case $n=2,3,4.$\par
\subsection{Permutation orbifold case}
\begin{theorem}\label{allpermu}
Every irreducible DHR representation of $\D^{\mathbb P_n}$ appears as 
an irreducible summand of $\pi_{\psi,p}$ for some $\psi,p\in \mathbb P_n$.
\end{theorem}
\proof
The proof is similar to the proof of Th. \ref{allcyclic} with 
small modifications. 
By Prop. \ref{solitonres} and  Th. \ref{permu}
$\pi_{\psi_1,p_1}\res \D^{\mathbb P_n} \simeq \pi_{\psi_2,p_2}\res 
\D^{\mathbb P_n} $ iff 
there exists $h\in \mathbb P_n$ such that 
$h.\psi_1=\psi_2,  hp_1h^{-1}= p_2$. Denote the orbit 
of $\pi_{\psi_1,p_1}$ under the action of $\mathbb P_n$ 
by $\{\psi_1,p_1 \}$. Note that the orbit   $\{\psi_1,p_1 \}$
has length $\frac{n}{ |\Gamma_{\psi_1,p_1}|}$.  
By Prop. \ref{solitonind} 
the sum of index of the irreducible summands of $\pi_{\lambda,p}$
is $\frac{n!^2}{|\Gamma_{\psi,p}| }d(\pi_{\psi,p})^2$. 
Hence the sum of index of 
distinct irreducible summands of  $\pi_{\psi,p}$ for all 
$\psi, p\in \mathbb P_n$ is given by$
\sum_{\{\psi, p \}} n d(\pi_{\lambda,p})^2$ where the sum
is over different orbits.  Assume that $p= c_1...c_k$ is a product of disjoint
cycles. 
For each element $\psi_2,p_2$
in the orbit $\{\psi, p\}$, by Prop. \ref{indexp} 
$d(\pi_{\psi_2,p_2})^2= d(\pi_{\psi,p})^2= 
\prod_{1\leq j\leq k} d(\psi(c_j)^{2} \mu_\A^{n-k}   
$
Hence
$$
\sum_{\{\psi, p \}} \frac{n!^2}{|\Gamma_{\psi,p}|}  d(\pi_{\psi,p})^2
=n! \sum_{\psi, p }  \prod_{1\leq j\leq k} d(\psi(c_j)^{2}
 \mu_\A^{n-k}
= n!^2 \mu_\A^n = \mu_{\D^{\mathbb P_n}}
$$
where in the last $=$ we have used Th. \ref{orb}. 
The theorem now follows from Th. 30 of \cite{KLM}.
\endproof   
Let
$p= c_1...c_k$ be  a product of $k$ disjoint cycles . 
Let $h\in \Gamma_{\psi,p}$. Then 
$\Ad_h$ induces a permutation among the cycles $c_1,...,c_k$. 
We define an element $h'\in \mathbb P_k$ by the formula
$ hc_ih^{-1}=c_{h'(i)}, i=1,...,k$. We note that in the definition
of $\pi_{\psi,g}$ a presentation of $g$ has been fixed. Assume that  
$hf_{{h'}^{-1}(i)}h^{-1}= h''(i). f_i$ where $h''(i)$ is an element in the
cyclic group generated by $c_i$. 
Define
\begin{definition}\label{act3}
$$\pi_{\psi,p}(h):= h'\big(\pi_{\psi(c_1),c_1}(h''(1)) \otimes\cdots
\otimes  \pi_{\psi(c_k),c_k}(h''(k)\big)
$$
where the action of $h'\in {\mathbb P}_k$ on $\H_{\psi(c_1)} \otimes\cdots\H_{\psi(c_k)}$ is by permutation of the tensor factors, and 
$ \pi_{\psi(c_i),c_i}(h''(i))$ is as defined in definition (\ref{act1}).
\end{definition}
One checks easily that Definition \ref{act3} gives a representation of
$ \Gamma_{\psi,p},$
$
\Ad_{\pi_{\psi,p}(h)} \pi_{\psi,p}= \pi_{\psi,p}\Ad_h
,$ and the vector 
$\Omega_{\psi(c_1)}\otimes\cdots \otimes \Omega_{\psi(c_k)}$
is fixed by $\pi_{\psi,p}(\Gamma_{\psi,p})$. It follows by 
Lemma \ref{fullspectrum1}, Th. \ref{projrep} that 
$ \pi_{\psi,p, \sigma\in \hat\Gamma_{\psi,p} }$ gives all the 
irreducible summands of  $\pi_{\psi,p}\res \D^{\mathbb P_n}$, 
and we have proved:\par
\begin{theorem}\label{decomposepermu}
$ \pi_{\psi,p, \sigma\in \hat\Gamma_{\psi,p} }$ gives all the 
irreducible summands of  $\pi_{\psi,p}\res \D^{\mathbb P_n}$.
\end{theorem}
Note that by Prop. \ref{decomposepermu} and Th. \ref{allpermu} 
the irreducible DHR representations
of $\D^{\mathbb P_n}$ are labeled by triples
$(\psi,p,\sigma)$ with $p.\psi=\psi, \sigma\in \hat\Gamma_{\psi,p}$
with equivalence relation $\sim$,  
$(\psi,p,\sigma)\sim (\psi_1,p_1,\sigma_1)$ iff there is $h\in \mathbb P_n$
such that $\psi_1=h.\psi, p_1= hph^{-1},\sigma_1= \sigma^h$. 
In \cite{Ba}, based on heuristic argument it is claimed that the 
irreducible representations of $\D^{\mathbb P_n}$ should be given by 
the set of pairs $(\psi, \phi)$ where $\phi$ is an irreducible 
representation of the double $\D(F_\psi)$ of the stabilizer
$F_\psi=\{ p\in \mathbb P_n| p.\psi=\psi \}$ with equivalence relation
$(\psi, \phi)\sim (\psi_1, \phi_1)$  iff there is $h\in \mathbb P_n$
such that $ \psi_1=h.\psi,   \phi_1= \phi^h$. We note that the irreducible
representation of  the double $\D(F_\psi)$ are labeled by 
$(g,\pi)/F_\psi $, where $g\in F_\psi$, $\pi$ is an irreducible representation
of the centralizer of $g$ in $F_\psi$, and the action of 
$F_\psi$ on $(g,\pi)$ is given by 
$h.(g,\pi)=(hgh^{-1}, \pi^h)$. Hence the labels \cite{Ba} are exactly
the same as the labels we described above, and we have confirmed this
claim of \cite{Ba}.
\section{Examples of fusion rules}
\subsection{Some properties of S matrix for general orbifolds}
Let $\A$ be a completely rational conformal net and let $\Gamma$ be a finite 
group acting properly on $\A$. By Th. \ref{orb} $\A^\Gamma$ has 
only finitely many irreducible representations. We use $\dot\lambda$ 
(resp. $\mu$) to
label representations of $\A^\Gamma$ (resp. $\A$). 
We will denote the corresponding
genus $0$ modular matrices by $\dot S, \dot T$ (cf. (\ref{Smatrix}). 
Denote by $\dot\lambda$ 
(resp. $\mu$) the irreducible covariant representations of ${\A}^\Gamma$
(resp. ${\A}$) with finite index. Recall that
$b_{\mu\dot\lambda}\in {\Bbb Z}$ denote
the multiplicity of representation $\dot\lambda$ which appears
in the restriction of  representation $\mu$ when restricting from 
$ {\A}$ to ${\A}^\Gamma$. 
$ b_{\mu\dot\lambda} $ is also known as the branching
rules. 
\begin{lemma}\label{Smatrix1}
(1)
If $\tau$ is an automorphism (i.e., $d(\tau)=1$) then
$S_{\tau(\lambda) \mu} = G(\tau,\mu)^* S_{\lambda\mu}$ where 
$\tau(\lambda):= \tau\lambda, G(\tau,\mu)=\epsilon(\tau,\mu)
\epsilon(\mu,\tau) ;$\par
(2) For any $h\in \Gamma$, let $h(\lambda)$ be the DHR representation 
$\lambda\cdot\Ad_{h^{-1}}$. Then
$ S_{\lambda\mu}=  S_{h(\lambda)h(\mu)}$;\par
(3) If $\lambda\rightarrow z(\lambda) \frac{ S_{\lambda\mu}}{S_{1\mu}}
$ gives a representation of the fusion algebra of $\A$ where $z(\lambda)$ is
a complex-valued function, $z(1)=1$, 
then there exists an automorphism  $\tau$ 
such that $z(\lambda)= \frac{S_{\lambda\tau}}{S_{\lambda 1}}$;\par
(4) If $[\alpha_{\dot \lambda}] = [\mu \alpha_{\dot \delta}]$, then for
any $\dot\lambda_1, \mu_1$ with $b_{\dot\lambda_1\mu_1}\neq 0$ we have
$
\frac{S_{\dot \lambda \dot\lambda_1}}{S_{\dot 1\dot\lambda_1}}
= \frac{S_{\mu\mu_1}}{S_{ 1\mu_1}}
\frac{S_{\dot \delta \dot\lambda_1}}{S_{\dot 1\dot\lambda_1}}
$;\par
\end{lemma}
\proof
Ad (1):
Since $\lambda\rightarrow \frac{S_{\lambda \mu}}{S_{1\mu}}$ is a 
representation of the fusion algebras, it follows that
$$
\frac{S_{\tau(\lambda) \mu}}{S_{1\mu}}=
\frac{S_{\lambda \mu}}{S_{1\mu}}\frac{S_{\tau \mu}}{S_{1\mu}}
$$
On the other hand
$$
\frac{S_{\tau \mu}}{S_{1\mu}}= \frac{\omega_{\tau}
\omega_{\mu}}{\omega_{\tau\mu}} = G(\tau,\mu)^*
$$
where the last equation follows from the monodromy equation (cf.\cite{R2})
and (1) is proved. \par
Ad (2)
By Lemma \ref{Yprop}, it is sufficient to show that 
$N_{h(\lambda)h(\mu)}^{h(\delta)}=N_{\lambda\mu}^{\delta}$
and $\omega_{h(\lambda)}= \omega_{\lambda}$. The first equation 
follows from the definition. For the second one, we note that
$\omega_{\lambda}= \pi_\lambda(R(2\pi))$.  
Since $h$ commutes with the vacuum unitary representation of $\Mob$,
it follows that $\omega_{h(\lambda)}= \omega_{\lambda}$.    
\par
Ad (3): 
By assumption $\lambda\rightarrow z(\lambda) \frac{S_{\lambda1}}{S_{11}}$
is a non-trivial representation of the fusion algebra, and so 
there exists $\tau$ such that
$ z(\lambda) \frac{S_{\lambda1}}{S_{11}}= \frac{S_{\lambda\tau}}{S_{1\tau}}
,\forall \lambda.$ Hence $ |z(\lambda)|\leq 1$. From
\begin{align}
 z(\lambda_1) \frac{S_{\lambda_1\mu}}{S_{1\mu}}
z(\lambda_2) \frac{S_{\lambda_2\mu}}{S_{1\mu}}
=  \sum_{\lambda_3} N_{\lambda_1\lambda_2}^{\lambda_3}
z(\lambda_1) z(\lambda_2) \frac{S_{\lambda_3\mu}}{S_{1\mu}}
\\
= \sum_{\lambda_3} N_{\lambda_1\lambda_2}^{\lambda_3}
z(\lambda_3) \frac{S_{\lambda_3\mu}}{S_{1\mu}}
\end{align}
we have
$$\sum_{\lambda_3} N_{\lambda_1\lambda_2}^{\lambda_3}
(z(\lambda_1) z(\lambda_2)-z(\lambda_3))  
\frac{S_{\lambda_3\mu}}{S_{1\mu}}=0
.$$
Using $ N_{\lambda_1\lambda_2}^{\lambda_3}= \sum_\delta
\frac{S_{\lambda_1\delta} S_{\lambda_2\delta} S_{\lambda_3\delta}^*}
{S_{1\delta}}$ and the orthogonal property of $S$ matrix in Lemma \ref{Sprop}
we have
$$
N_{\lambda_1\lambda_2}^{\lambda_3}(z(\lambda_1) z(\lambda_2)-z(\lambda_3))=0  
.$$
Since $N_{\lambda_1\bar\lambda_1}^{1}=1$ we have
$z(\lambda_1) z(\bar\lambda_1)=1$. So we conclude that
$ |z(\lambda)|=1, \forall \lambda$, and
$$
|\frac{1}{S_{1\tau}}|^2
\sum_\lambda |\frac{S_{\lambda\tau}}{S_{1\tau}}|^2= 
\sum_\lambda |\frac{S_{\lambda 1}}{S_{11}}|^2
=|\frac{1}{S_{11}}|^2
$$
Hence $S_{1\tau}=S_{11}$ and $d(\tau)=1$, i.e., $\tau$ is an automorphism.\par
Ad (4):
By \cite{Xu4} or \cite{BEK} 
there is a unit vector $\psi$ in the vector space spanned by 
the irreducible components of $\alpha_{\dot \lambda_2}, \forall \lambda_2$
such that 
$$
\alpha_{\dot\lambda} \psi= \frac{S_{\dot\lambda\dot\lambda_1}}
{S_{\dot 1\dot\lambda_1}} \psi, 
\mu \psi = \frac{S_{\mu\mu_1}}{S_{1\mu_1}},
\alpha_{\dot\delta} \psi= \frac{S_{\dot\delta\dot\lambda_1}}
{S_{\dot 1\dot\lambda_1}} \psi
$$ 
and (4) follows immediately.
\endproof
\subsection{Fusions of solitons in cyclic orbifolds}
Let $\B\subset \D$ be as in \S6.1. 
Set $i=0$ in  Th. \ref{decomposecy}. In this $\psi$ is a constant function, 
and we denote it by its value $\lambda$. For simplicity we will
label the representation $\pi_{\lambda,g^j,i}$ ($g=(01...n-1)$) by
$(\lambda,g^j,i)$. 
Define $(\lambda i):=(\lambda,1,i)$ where $i\in \widehat{\mathbb Z_n}\simeq 
\mathbb Z_n$.   
\begin{lemma}\label{sumN}
If $(k,n)=1$, then
$$
\sum_{ 0\leq j\leq n-1} e^{\frac{2\pi i kj}{n}}
N_{(\lambda 0)(\mu 0)}^{(\delta 0)} = N_{\lambda \mu}^{\delta}
$$
\end{lemma}   
\proof
Let $V:= {\rm Hom}(\delta,\lambda\mu) \subset \A(J_0)$. 
Note that $\mathbb Z_n$ acts on $W:=V\otimes V\otimes\cdots\otimes V$(n-tensor
factors) by permutations.
Let $W_j:= \{ w\in W| \beta_g(w)= e^{-2\pi ij}{n} w\}$. Note that if
$w\in W_j$, then $w v^{j} \in {\rm Hom} (v^{-j} \delta^{\otimes^n} v^j, 
\lambda^{\otimes^n}\cdot\mu^{\otimes^n}) \cap \D^{\mathbb Z_n}(J_0),
$ where $v$ is defined as before Lemma \ref{grading}.
Hence we have an injective map $w\in W_j\rightarrow w v^{j} \in {\rm Hom}
((\delta j), (\lambda 0)(\mu 0))$. By definition the map is also
surjective. So we have
$$
\sum_{ 0\leq j\leq n-1} e^{-\frac{2\pi i kj}{n}}
N_{(\lambda 0)(\mu 0)}^{(\delta 0)} =\sum_{ 0\leq j\leq n-1} 
e^{-\frac{2\pi i kj}{n}} {\rm dim} W_j 
= \Tr_W (\beta_{g^k})
$$ 
When $(k,n)=1$, $g^k$ is one cycle, and it follows that
$ Tr_W (\beta_{g^k})= {\rm dim} V =  N_{\lambda \mu}^{\delta}.$
Take the complex conjugate of both sides we have proved the lemma.
\endproof
\begin{lemma}\label{N}
Let $ f_\mu:=(\mu,g,0)$. Then: \par
(1) $G(\sigma,f_\mu)= e^{\frac{2\pi l_1 i}{n}}$ for some integer $l_1$
with $(l_1,n)=1$;\par
(2) $\lambda\rightarrow \frac{S_{(\lambda 0)f_\mu}}{S_{(10) f_\mu}}$
is a representation of the fusion algebra of $\A$; \par
(3) There exists an automorphism $\tau, [\tau^2]=[1]$ such that 
$
\frac{S_{(\lambda 0)f_\mu}}{S_{(10)f_\mu}}=
\frac{S_{\lambda \tau(\mu)}}{S_{1 \tau(\mu)}} 
$
\end{lemma}
\proof
Ad(1): By the paragraph after (47) in \cite{LX} we have
$G(\sigma^{k(1)},f_\mu)=  e^{\frac{2\pi i}{n}}$ where $(k(1),n)=1$.
By (1) of Lemma \ref{grading} we have 
$G(\sigma,f_\mu)^{k(1)}=  e^{\frac{2\pi i}{n}}.$ Choose $l_1$ such that
$l_1k(1)=1 {\rm mod} n$ we have 
$G(\sigma,f_\mu)= e^{\frac{2\pi l_1 i}{n}}$ for some integer $l_1$
with $(l_1,n)=1$. \par
As for (2) and (3), first we note that by  Lemma \ref{grading}, 
if $\delta \prec (\lambda 0) (\mu 0)$,
then $\delta$ is untwisted. Suppose that $\delta$ is an irreducible 
component of the restriction of $(\delta_1,...,\delta_n)$ to
$\D^{\mathbb Z_n}$. We claim that 
$S_{\delta f_\mu}=0$ if $\delta_i\neq \delta_j$ for some $i\neq j$. In fact
if  $\delta_i\neq \delta_j$ for some $i\neq j$, then the stabilizer of
$(\delta_1,...,\delta_n)$ under the action of $\mathbb Z_n$ is a proper
subgroup of  $\mathbb Z_n$, and by Th. \ref{solitonind}  
$\alpha_\delta $ is reducible,
and $[\sigma^k \delta]=[\delta]$ for some $1\leq k\leq n-1 $ . By
(1) of Lemma \ref{Smatrix} we have
$S_{\delta f_\mu}=  S_{\sigma^k (\delta) f_\mu}= 
G(\sigma^k, f_\mu)^*S_{\delta f_\mu}
$
Since $G(\sigma^k, f_\mu)= e^{\frac{2\pi kl_1}{n}} $ 
with $(l_1,n)=1$ by (1),  $G(\sigma^k, f_\mu)^*\neq 1$, hence
$S_{\delta f_\mu}=0$. So we have
\begin{align}
\frac{S_{(\lambda_1 0)f_\mu}}{S_{(10)f_\mu}} 
\frac{S_{(\lambda_2 0)f_\mu}}{S_{(10)f_\mu}}
=\sum_{\lambda_3, 0\leq j\leq n-1}
N_{(\lambda_1 0)(\lambda_2 0)}^{((\lambda_3 j)} 
\frac{S_{(\lambda_3 j)f_\mu}}{S_{(10)f_\mu}}
\\
= \sum_{\lambda_3, 0\leq j\leq n-1}
N_{(\lambda_1 0)(\lambda_2 0)}^{(\lambda_3 j)} e^{\frac{2\pi i j l_1}{n}} 
\frac{S_{(\lambda_3 0)f_\mu}}{S_{(10)f_\mu}}\\
= N_{\lambda_1\lambda_2}^{\lambda_3} 
\frac{S_{(\lambda_3 0)f_\mu}}{S_{(10)f_\mu}}
\end{align}
where we have used (1) of Lemma \ref{Smatrix1} and 
Lemma \ref{sumN} in the second $=$ and third $=$ respectively. \par 
Ad (2):
Since $\alpha_{f_\mu}= (\mu,1,...,1)\alpha_{f_1}$ by (48) of \cite{LX},  
by (4) of Lemma \ref{Smatrix1} we have
$$
\frac{S_{f_\mu (\lambda 0)}}{S_{(10)(\lambda 0)}}= 
\frac{S_{\mu \lambda}}{S_{1\lambda}}  
\frac{S_{f_1 (\lambda 0)}}{S_{(10)(\lambda 0)}}
$$
Combined with (1) it follows that there exists $\tau$ such that
the map 
$$
\lambda\rightarrow \frac{S_{\lambda\tau}}{d(\lambda) S_{1\tau}} 
\frac{S_{\lambda\mu}}{S_{1\mu}} 
$$
gives a representation of the fusion algebra of $\A$. By (3) of lemma
\ref{Smatrix1} we have that $\tau$ is an automorphism and
$$
\frac{S_{(\lambda  0)f_1}}{S_{(10)f_1}}
= \frac{S_{\lambda \tau}}{S_{1 \tau}}
$$
Let $h \in \mathbb P_n$ such that $h gh^{-1}= g^{-1}$. By definition
$h((\lambda  0)) =(\lambda 0)$. By Prop. \ref{pind1} 
$[h(f_1)]=[\sigma^j(\bar f_1)]$ for
some $1\leq j\leq n$,  and it follows
by Lemma \ref{Smatrix1} 
that 
$$
S_{(\lambda  0)f_1}=S_{h((\lambda  0))h(f_1)}=
S_{(\lambda  0) \sigma^j(\bar f_1)}=S_{(\lambda  0) \bar f_1}
=S^*_{(\lambda  0) f_1}, 
$$
hence 
$\frac{S_{\lambda \tau}}{S_{1 \tau}}= \frac{S_{\lambda \bar\tau}}{S_{1 \tau}}
, \forall \lambda$, and so
$[\tau]=[\bar \tau].$
\endproof
We conjecture that $[\tau]=[1]$ in the above lemma.  \par
Let $f_1:=(1,g,0)$ where $0$ stands for the trivial representation of 
$\mathbb Z_n$.  In \cite{LX} the questions about the nature of
$[\alpha_{f_1}^k]= [\pi_{1,g}^k]$ (cf. (44) of \cite{LX}) 
where $k$ is an integer is raised. 
\begin{proposition}\label{degreen}
When $n$ is even we have
$$
[\pi_{1,g}^n]=\bigoplus_{\lambda_1,...,\lambda_n} 
M_{\lambda_1,...,\lambda_n} [(\lambda_1,...,\lambda_n)]
$$ where
$M_{\lambda_1,...,\lambda_n}:=\sum_{\lambda} S_{1\lambda}^{2-2g}
\prod_{1\leq i\leq n} \frac{S_{\lambda_i\lambda}}{ S_{1\lambda}}
$ 
with $g=\frac{(n-1)(n-2)}{2}$.
\end{proposition}  
\proof
We note that by Lemma \ref{grading}  $\pi_{1,g,0}^{n-1}$ is untwisted,
and must be sum of irreducible untwisted representations. It follows that 
by Cor. 8.4 of \cite{LX} that 
$$[\alpha_f^n]=\bigoplus_{\lambda_1,...,\lambda_n} 
M_{\lambda_1,...,\lambda_n} [(\lambda_1,...,\lambda_n)]
$$
with $M_{\lambda_1,...,\lambda_n}$ non-negative integers. 
Let $\mu$ be any irreducible subsector of $\alpha_f^{n-1}$. By the equation
above $\mu\alpha_f\succ (\lambda_1,...,\lambda_n)$ for 
some $(\lambda_1,...,\lambda_n)$, and by Frobenius duality 
$\mu\prec (\lambda_1,...,\lambda_n) \bar\alpha_f$. By (46) of \cite{LX}
$[(\lambda_1,...,\lambda_n) \bar\alpha_f]
= \sum_{\lambda} \langle\lambda_1\cdots \lambda_n, \lambda\rangle
[(\lambda,1,1,...,1)\bar\alpha_f]
$ 
and by (48) of \cite{LX} each $(\lambda,1,1,...,1)\bar\alpha_f$ 
is irreducible. Hence $[\mu]= [(\lambda,1,1,...,1)\bar\alpha_f]$
for some $\lambda$. 
Hence 
$$
[\alpha_{f_1}^{n-1}]=
[\pi_{1,g}^{n-1}]=\bigoplus_{\lambda} m_\lambda [(\lambda,1,...,1) \bar 
\pi_{1,g}]
$$
with $m_\lambda$  non-negative integers.   
By (4) of Lemma \ref{Smatrix1} we have
$$
(\frac{S_{f_1(\mu 0)}}{ S_{(10)(\mu 0)}})^{n-1}
=\sum_{\lambda} m_\lambda \frac{S_{\lambda \mu}}{S_{1 \mu}}  
\frac{S^*_{f_1(\mu 0)}}{ S_{(10)(\mu 0)}}
$$
Note that 
$ \frac{1}{S_{(10)(10)}^2}= \mu_{\D^{\mathbb Z_n}}= 
n^2 \mu_\D = n^2  \frac{1}{S_{11}^{2n}},$ hence
$ S_{(10)(10)} = \frac{S_{11}^n}{n}$. From 
$
\frac{S_{(\lambda 0)(10)}}{S_{(1 0)(10)}}= d((\lambda 0))= n d(\lambda)^n
$ we have
$
S_{(\lambda 0)(10)} = S_{\lambda 1}^n
.$ 
By (2) of Lemma \ref{N} and our assumption that $n$ is even and
hence $[\tau^n]=[1]$, we have
$$
\frac{1}{S_{1\mu}^{(n-1)(n-2)}}
= \sum_\lambda m_\lambda \frac{S_{\lambda\mu}^*}{S_{1\lambda}}
$$
By orthogonal property of $S$ matrix in Lemma \ref{Sprop} we have
$$
m_\lambda = \sum_{\mu} \frac{S_{\lambda\mu}}{S_{1\mu}^{ n^2-3n+1}}
$$
Combine this with (46) of \cite{LX} and (\ref{Verlinde}) 
the proposition follows.
\endproof
We remark $M_{\lambda_1,...,\lambda_n}$ is the dimension of 
genus $\frac{(n-1)(n-2)}{2}$ conformal
blocks with the insertion of representations $\lambda_1,...,\lambda_n$. 
Also note that  $\frac{(n-1)(n-2)}{2}$ is the genus of an algebraic 
curve with degree $n$. It may be interesting to give a geometric 
interpretation of Prop.\ref{degreen}.\par
Note that the conjecture $[\tau]=[1]$ at the end of previous section is
true, then the above proposition is also true for odd $n$. 
\subsection{n=2 case}
In this section we consider the fusions rules for the simplest non-trivial
case $n=2$. Partial results have been obtained in \S8 of \cite{LX}. We will
confirm the results in \S4.6 of \cite{BHS}. Let us first simply our
notations by introducing similar notations in \cite{BHS}.
Let $\widehat{(\lambda 0)}: =(\lambda, -1, 0),\widehat{(\lambda 1)}: 
=(\lambda, -1, 1)$. Note that by \S2 of \cite{LX} we can choose
\begin{equation}\label{conformaldt}
\omega_{\widehat {(\lambda 0)}}= e^{\frac{2\pi i (\Delta_\lambda+\frac{c}{8})}
{2}}, 
\omega_{\widehat {(\lambda 1)}}= e^{\frac{2\pi i (\Delta_\lambda+1+\frac{c}{8})}
{2}}  
\end{equation}
where $c$ is the central charge.  
We also note by definitions
\begin{equation}\label{conformaldu}
\omega_{(\lambda 0)}= e^{4 \pi i \Delta_\lambda}
=\omega_{(\lambda 1)}^2, \
\omega_{(\lambda_1\lambda_2)}= e^{2\pi i (\Delta_{\lambda_1}+
\Delta_{\lambda_2}) }, \lambda_1\neq \lambda_2
\end{equation}
\begin{lemma}\label{tt}
(1)
$$
N_{\widehat{(\lambda_1 \epsilon_1)} \widehat{(\lambda_2 \epsilon_2)} }^{
(\lambda_30)}+N_{(\lambda_1 \epsilon_1)(\lambda_2 \epsilon_2)}^{(\lambda_31)} 
=
\sum_{\mu} \frac{S_{\bar\lambda_3 \mu}^2 S_{\lambda_1 \mu}S_{\lambda_2 \mu}}
{S_{1\mu}^2}
$$ where $\epsilon_1, \epsilon_2= 0 \ {\rm or} \ 1$;\par
(2)
$$
N_{(\lambda_1 \epsilon_1)(\lambda_2 \epsilon_2) }^{
(\lambda_4 \lambda_5)}
= \sum_{\mu}
\frac{S_{\bar\lambda_4 \mu} S_{\bar\lambda_5 \mu} S_{\lambda_1 \mu}
S_{\lambda_2 \mu}}
{S_{1\mu}^2};
$$\par
(3)
$$
\sum_{\mu} \frac{S_{\lambda_3 \mu}^2 S_{\lambda_1 \mu}S_{\lambda_2 \mu}}
{S_{1\mu}^2}  \frac{1}{\omega_{\lambda_3}^2} d((\lambda_30))
+ \sum_{\mu, \lambda_4 \neq \lambda_5}
\frac{S_{\lambda_4 \mu} S_{\lambda_5 \mu} S_{\lambda_1 \mu}S_{\lambda_2 \mu}}
{S_{1\mu}} \frac{1}{\omega_{\lambda_3}^2}  d(\lambda_4) d (\lambda_5))
=\frac{ S_{\lambda_1 \mu}   T_\mu^2 S_{\lambda_2 \mu}}{S_{11}^2}
e^{\frac{-2\pi i c_0}{6}}
$$
where $c_0$ is defined as in (\ref{c0}).
\end{lemma}
\proof
Ad (1):
By (48) and (3) of Prop. 8.8 in \cite{LX} we have
$$
\langle \alpha_{\widehat{(\lambda_1 \epsilon_1)}}  
\alpha_{\widehat{(\lambda_2 \epsilon_2)}}, 
(\lambda_3,\lambda_3) \rangle
= \langle \lambda_1\lambda_2 \bar\lambda_3\bar\lambda_3,1\rangle
$$
Note that 
$$
\langle \alpha_{(\lambda_1 \epsilon_1)}  \alpha_{(\lambda_2 \epsilon_2)}, 
(\lambda_3,\lambda_3) \rangle= 
N_{(\lambda_1 \epsilon_1)(\lambda_2 \epsilon_2) }^{
(\lambda_30)}+N_{(\lambda_1 \epsilon_1)(\lambda_2 \epsilon_2)}^{(\lambda_31)}
$$
and by (\ref{Verlinde})
(1) is proved. (2) is proved in a similar way.\par
Ad (3):
\begin{align}
\sum_{\mu} \frac{S_{\lambda_3 \mu}^2 S_{\lambda_1 \mu}S_{\lambda_2 \mu}}
{S_{1\mu}^2}  \frac{1}{\omega_{\lambda_3}^2} d((\lambda_30))
+ \sum_{\mu, \lambda_4 \neq \lambda_5}
\frac{S_{\lambda_4 \mu} S_{\lambda_5 \mu} S_{\lambda_1 \mu}S_{\lambda_2 \mu}}
{S_{1\mu}^2} \frac{1}{\omega_{\lambda_3}^2}  d(\lambda_4) d (\lambda_5))
\\
=\sum_{\mu} \frac{S_{\lambda_1 \mu}S_{\lambda_2 \mu}}
{S_{1\mu}^2} (\sum_{\lambda_3}  
\frac{S_{\bar\lambda_3 \mu} S_{\bar\lambda_3 1}}{
\omega_{\lambda_3}})^2 
\end{align}

From Lemma \ref{Sprop} we have
$
S^* T^{-1} S^*= TS^*T
$
and so 
$$
\sum_{\lambda_3} S_{\bar\lambda_3\mu} S_{\lambda_3 1} 
\frac{1}{\omega_{\lambda_3}}
=S_{1\mu} e^{\frac{-\pi i c_0}{12}}  T_\mu
$$
Substitute into the equations above we have proved (3).
\endproof
Define matrices  $ T^{\frac{1}{2}}$ such that $ T^{\frac{1}{2}}_{\lambda\mu}
= \delta_{\lambda\mu} e^{\pi i (\Delta_\lambda- \frac{c_0}{24})}$ and 
\begin{definition}\label{Pmatrix}
$$P:=T^{\frac{1}{2}} ST^2 S T^{\frac{1}{2}}, 
\widetilde{P}= e^{\frac{2\pi i(c-c_0)}{8}} P
.$$
\end{definition} 
It follows by (\ref{Ymatrix}) that 
\begin{align}
Y_{\widehat{(\lambda_1 \epsilon_1)}\widehat{(\lambda_2 \epsilon_2)} }
=\omega_{(\widehat{\lambda_1 \epsilon_1)}} \omega_{(\widehat{\lambda_2 \epsilon_2)}}
\times
(\sum_{\lambda_3} (N_{(\lambda_1 \epsilon_1)(\lambda_2 \epsilon_2) }^{
(\lambda_30)}+N_{(\lambda_1 \epsilon_1)(\lambda_2 \epsilon_2) }^{
(\lambda_31)}) \frac{1}{\omega_{\lambda_3}^2} d((\lambda_30))
+  \\
\frac{1}{2} 
\sum_{\lambda_4 \neq \lambda_5} (N_{(\lambda_1 \epsilon_1)(\lambda_2 \epsilon_2) }^{
(\lambda_4 \lambda_5)}+N_{(\lambda_1 \epsilon_1)(\lambda_2 \epsilon_2) }^{
(\lambda_31)}) \frac{1}{\omega_{\lambda_4} \omega_{\lambda_5}}  d((\lambda_4 \lambda_5))  
) \\
= \omega_{(\lambda_1 \epsilon_1)} \omega_{(\lambda_2 \epsilon_2)}
\times (\sum_{\mu} \frac{S_{\lambda_3 \mu}^2 S_{\lambda_1 \mu}S_{\lambda_2 \mu}}
{S_{1\mu}^2}  \frac{1}{\omega_{\lambda_3}^2} d((\lambda_30))
+ \\
\sum_{\mu, \lambda_4 \neq \lambda_5}
\frac{S_{\lambda_4 \mu} S_{\lambda_5 \mu} S_{\lambda_1 \mu}S_{\lambda_2 \mu}}
{S_{1\mu}^2} \frac{1}{\omega_{\lambda_4}\omega_{\lambda_5} }  d(\lambda_4) d (\lambda_5))
=  e^{\pi i (\epsilon_1+  \epsilon_2)} 
\frac{\widetilde{P}_{\lambda_1 \lambda_2}}{S_{11}^2}
\end{align}
where in the last $=$ we have used (3) of Lemma \ref{tt}. 
Note that $S_{(10)(10)}^2= \frac{1}{4\mu_D}= \frac{1}{4\mu_\A^2}=
\frac{1}{4} S_{11}^4$, and so $S_{(10)(10)}=\frac{1}{2} S_{11}^2$.  
It follows  by (\ref{Smatrix}) that 
$$
S_{\widehat{(\lambda_1\epsilon_1)}\widehat{(\lambda_2\epsilon_2)}}
= e^{\pi i (\epsilon_1+  \epsilon_2)} \frac{1}{2}  
\widetilde{P}_{\lambda_1 \lambda_2}
$$
Note that by Lemma \ref{N} we have
$$
S_{(\lambda 0)\widehat{(\mu 0)}} 
= S_{(\lambda 0)\widehat{(\mu 0)}} \frac{S_{\lambda\mu}}{S_{1\mu}}
\times \frac{S_{\tau\lambda}}{S_{1\lambda}} 
$$
Since $[\tau^2]=[1]$, $\frac{S_{\tau\lambda}}{S_{1\lambda}}^2=1$, 
and so $\frac{S_{\tau\lambda}}{S_{1\lambda}}=\pm 1$. By (1) of 
Lemma \ref{Smatrix1} we can choose our labeling $\widetilde{(\lambda0)}, 
\widetilde {(\lambda 1)}$ such that as a set $\{ \widetilde {(\lambda0)}, 
\widetilde {(\lambda 1)} \}$ is the same as 
 $\{  (\lambda0), 
 (\lambda 1) \}$ and
$$
S_{\widetilde{(\lambda \epsilon)}\widehat{(\mu 0)}} 
= e^{\pi i \epsilon} 
S_{(1 0)\widehat{(1 0)}} \frac{S_{\lambda\mu}}{S_{1\mu}}
$$
From 
$$
[\alpha_{(\lambda_1\lambda_2)}]=[(\lambda_1,\lambda_2)] +
[(\lambda_2,\lambda_1)] 
$$
and (4) of Lemma \ref{Smatrix1} we have
$$
\frac{S_{(\lambda_1\lambda_2)(\lambda \epsilon)}}
{S_{(10)(\lambda \epsilon)}} 
= \frac{S_{\lambda_1\lambda}}
{S_{1\lambda}^2} +
 \frac{S_{\lambda_1\lambda}}
{S_{1\lambda}^2},
\frac{S_{(\lambda_1\lambda_2)(\lambda \mu)}}
{S_{(10)(\lambda \mu
)}} 
= \frac{S_{\lambda_1\lambda}S_{\lambda_1\mu} }
{S_{1\lambda} S_{1\mu}} +
\frac{S_{\lambda_2\lambda} S_{\lambda_2\mu} }
{S_{1\lambda} S_{1\mu}}
$$
Since  
$$
S_{(10)(10)} = \frac{1}{2} S_{11}^2
,$$  we get the following on the entries of $S$-matrix of $\D^{\mathbb Z_2}$: 
\begin{align}\label{Smatrix2}
S_{(\lambda\mu)(\lambda_1\mu_1)}
=S_{\lambda\lambda_1} S_{\mu\mu_1}+ 
S_{\lambda\mu_1} S_{\mu\lambda_1}, 
S_{(\lambda\mu)}{\widetilde{(\lambda_1\epsilon)}}=
S_{\lambda\lambda_1} S_{\mu\lambda_1} \\
S_{(\lambda\mu}{\widehat{(\lambda_1\mu_1)}}= 0, 
S_{\widetilde{(\lambda,\epsilon)} \widetilde{(\lambda_1,\epsilon_1)}}= 
\frac{1}{2} S_{\lambda\lambda_1}^2 \\
S_{\widetilde{(\lambda,\epsilon)} \widehat{(\lambda_1,\epsilon_1)}}  
=\frac{1}{2} e^{\pi i \epsilon} S_{\lambda\lambda_1}, 
S_{\widehat{(\lambda,\epsilon)} \widehat{(\lambda_1,\epsilon_1)}}  
=\frac{1}{2} e^{\pi i (\epsilon+\epsilon_1)} P_{\lambda\lambda_1}, \epsilon=0,1
\end{align}
Denote by $\dot c_0$ the number (well defined ${\rm mod} 8\mathbb Z$) of
$\D^{\mathbb Z_2}$ (cf. (\ref{c0}). 
\begin{lemma}\label{central}
(1) $\dot c_0-2 c_0 \in  8\mathbb Z$; \par
(2) $c_0-c \in 4\mathbb Z.$
\end{lemma}
\proof
By Lemma \ref{Sprop} we have 
\begin{equation}\label{sts}
 S T S= { T}^{-1}  S { T}^{-1} 
\end{equation} 
First let us compare the $(10)(10)$ entry of both sides for $S,T$ matrix
of $\D^{\mathbb Z_2}.$  By using the
formula before the lemma we have:
$$(\sum_{\lambda} S_{1\lambda}^2 e^{2\pi i \Delta_\lambda})^2=
e^{\frac{2\pi i \dot c_0}{8}} S_{11}^2.
$$
On the other hand comparing the entry $11$ of (\ref{sts}) for  
$S$-matrix of $\D$ we have
$$
\sum_{\lambda} S_{1\lambda}^2 e^{2\pi i \Delta_\lambda}= e^{\frac{2\pi i 
c_0}{8}} S_{11}
,$$ 
and (1) follows by combining the two equations. \par
As for (2), we compare the entry $\widehat{(\lambda 0)} (10)$ of
both sides of (\ref{sts}). By using the equations  before the 
lemma the $\widehat{(\lambda 0)} (10)$ entry of the left hand side
of  (\ref{sts}) is given by $
e^{\frac{2\pi i(c-c_0)}{8}}e^{\frac{\pi i(c_0)}{24}}
$ multiplied by the $\lambda 1$ entry of the matrix $PT^{\frac{1}{2}}S.$
By applying (\ref{sts}) to $S,T$ matrix of $\D$ we have
$$
PT^{\frac{1}{2}}S= T^{\frac{1}{2}} ST^2 STS= T^{\frac{-1}{2}} S T^{-2}
.$$
Using these equations to compare with the 
$\widehat{(\lambda 0)} (10)$ entry of right hand side
of  (\ref{sts}) we have 
$ e^{\frac{2\pi i(c-c_0)}{4}}=1$ and (2) if proved.
\endproof
By (\ref{Verlinde}) and (2) of Lemma \ref{central} 
we immediately obtain the following fusion rules:
\begin{align}\label{Fu}
N_{(\lambda\mu)(\lambda_1\mu_1)}^{ (\lambda_2\mu_2)}
= N_{\lambda\lambda_1}^{\lambda_2}  N_{\mu\mu_1}^{\mu_2}+ 
 N_{\lambda\lambda_1}^{\mu_2}  N_{\mu\mu_1}^{\lambda_2}+ 
 N_{\lambda\mu_1}^{\mu_2}  N_{\mu\lambda_1}^{\lambda_2}+
 N_{\lambda\mu_1}^{\lambda_2}  N_{\mu\mu_2}^{\lambda_1} \\
N_{(\lambda\mu)(\lambda_1\mu_1)}^{ \widetilde{(\lambda_2\epsilon)}}
=  N_{\lambda\lambda_1}^{\lambda_2}  N_{\mu\mu_1}^{\lambda_2}
+  N_{\lambda\mu_1}^{\lambda_2}  N_{\mu\lambda_1}^{\lambda_2}\\
N_{\widetilde{(\lambda\epsilon)}\widetilde{(\lambda_1\epsilon_1)}}^{ 
\widetilde{(\lambda_2\epsilon)}} = 
\frac{1}{2} N_{\lambda\lambda_1}^{\lambda_2}( 
N_{\lambda\lambda_1}^{\lambda_2} + e^{\pi i(\epsilon+ \epsilon_1+
\epsilon_2)})\\
N_{\widehat{(\lambda\epsilon)}\widehat{(\lambda_1\epsilon_1)}}^{\lambda_2\mu_2}
= \sum_{\mu} N_{\lambda \lambda_1}^{\mu}   N_{\mu \bar\lambda_2}^{\mu_2}\\
N_{\widehat{(\lambda\epsilon)}\widehat{(\lambda_1\epsilon_1)}}^{\widehat
{(\lambda_2\epsilon_2)}}
= \frac{1}{2}\sum_{\mu}\frac{S_{\lambda\mu}^2 S_{\lambda_1\mu}
S_{\lambda_2\mu}}{ S_{1\mu}^2} 
+ \frac{1}{2}e^{\pi i (\epsilon+\epsilon_1+\epsilon_2)} 
\sum_{\mu} \frac{S_{\lambda\mu} P_{\lambda_1\mu}
P_{\lambda_2\mu}}{ S_{1\mu}^2}
\end{align}
where $\epsilon,\epsilon_1, \epsilon_2= 0 \ {\rm or } \ 1$.
Let us summarize the above equations in the following:   
\begin{theorem}\label{n=2f}
The fusion rules of $\D^{\mathbb Z_2}$ are given by the above equations.
\end{theorem}
From the theorem we  immediately have:
\begin{corollary}\label{integer}
For any completely rational $\A$ 
$$
\frac{1}{2}\sum_{\mu} \frac{S_{\lambda_1\mu}^2 S_{\lambda_2\mu}
S_{\lambda_3\mu}}{ S_{1\mu}^2} 
\pm \frac{1}{2}\sum_{\mu} \frac{S_{\lambda_1\mu} P_{\lambda_2\mu}
P_{\lambda_3\mu}}{ S_{1\mu}^2} 
$$
is a non-negative integer where $P$ is defined in (\ref{Pmatrix}).
\end{corollary}
Cor.  \ref{integer} confirmed a conjecture in \S4.6 of \cite{BHS}.  We 
note that even for known examples the direct confirmation of Cor.  
\ref{integer} seems to be very tedious.\par It will be an interesting 
question to generalize our results to $n>2$ cases.

\end{document}